\documentclass[11pt,reqno]{amsart}
\usepackage[a4paper]{geometry}

\usepackage[utf8]{inputenc}
\usepackage{amsmath,amssymb,amsthm}
\usepackage{mathtools}
\usepackage[english]{babel}
\usepackage[toc]{appendix}
\usepackage{dsfont}
\usepackage{esint}
\usepackage{url}
\usepackage{enumerate}

\usepackage{fullpage}

\usepackage{palatino}

\usepackage{comment}

\usepackage{pgf,tikz,pgfplots}
\pgfplotsset{compat=1.15}
\usetikzlibrary{arrows}

\usepackage{graphicx}

\usepackage[pagebackref]{hyperref}
\hypersetup{
    backref=true, 
    pagebackref=true,
    plainpages=false,
    bookmarks=true,         
    unicode=false,          
    pdftoolbar=true,        
    pdfmenubar=true,        
    pdffitwindow=true,      
    pdftitle={My title},    
    pdfauthor={Author},     
    pdfsubject={Subject},   
    pdfcreator={Creator},   
    pdfproducer={Producer}, 
    pdfkeywords={keywords}, 
    pdfnewwindow=true,      
    colorlinks=true,       
    linkcolor=blue,          
    citecolor=blue,        
    filecolor=blue,      
    urlcolor=blue,          
    urlbordercolor=0 1 1
}

\newcommand{\R}{\mathbb R}
\newcommand{\N}{\mathbb N}

\newcommand{\cW}{\mathcal{W}}
\newcommand{\cM}{\mathcal{M}}

\def\BV{{\mathrm{BV}}}

\def\wpq#1#2{{{\mathrm W}^{#1,#2}}}
\def\wpqO#1#2{{{\mathrm W}_0^{#1,#2}}}
\def\lp#1{{\mathrm L}^{#1}}

\def\cont{{\mathrm{C}}}
\def\Hau{{\mathcal H}}
\def\hausp#1{{\mathcal H}^{{#1}}}

\def\haushd{\hausp{d-1}}
\def\one#1{\mathds{1}_{#1}}

\def\S{\mathcal{S}}
\def\div{\operatorname{div}}
\def\spt{\operatorname{spt}}

\newtheorem{theo}{Theorem}[section]

\newtheorem{lemma}[theo]{Lemma}
\newtheorem*{theo*}{Theorem}

\theoremstyle{remark}
\newtheorem{remk}[theo]{Remark}

\newcommand{\E}{\mathcal{E}}
\DeclareMathOperator*{\argmin}{arg\,min}

\DeclareMathOperator*{\dist}{dist}

\renewcommand{\epsilon}{\varepsilon}
\newcommand{\e}{\epsilon}

\def\ds{\displaystyle}

\newcommand{\res}{\mathop{\hbox{\vrule height 7pt width .5pt depth 0pt
\vrule height .5pt width 6pt depth 0pt}}\nolimits}

\definecolor{Blue}{rgb}{0.,0.,1}

\newcommand{\wto}{\rightharpoonup}

\numberwithin{equation}{section} 

\title{A Cahn--Hilliard--Willmore phase field model\\ for non-oriented interfaces}

\author{Elie Bretin}
\address{INSA Lyon, CNRS, Ecole Centrale de Lyon, Universit\'e Claude Bernard Lyon 1, Universit\'e Jean Monnet, ICJ UMR5208,
69621 Villeurbanne, France \\ elie.bretin@insa-lyon.fr}
\author{Antonin Chambolle}
\address{CEREMADE, CNRS, Université Paris-Dauphine, PSL University, Paris, MOKAPLAN, INRIA, Paris,
France \\ chambolle@ceremade.dauphine.fr.}
\author{Simon Masnou}
\address{Universit\'e Claude Bernard Lyon 1, CNRS, Ecole Centrale de Lyon, INSA Lyon, Universit\'e Jean Monnet, ICJ UMR5208, 69622 Villeurbanne, France \\masnou@math.univ-lyon1.fr}

\makeatletter
\@namedef{subjclassname@2020}{\textup{2020} Mathematics Subject Classification}
\makeatother

\date{\today}

\begin{document}

\begin{abstract}
We investigate a new phase field model for representing non-oriented interfaces, approximating their area and simulating their area-minimizing flow. Our contribution is related to the approach proposed in~\cite{MR4483340} that involves ad hoc neural networks. We show here that, instead of neural networks, similar results can be obtained using a more standard variational approach that combines a Cahn-Hilliard-type functional involving an appropriate non-smooth potential and a Willmore-type stabilization energy.  We give a $\Gamma$-convergence analysis of this phase field model in dimension $1$ and, for radially symmetric functions, in arbitrary dimension. We also propose a simple numerical scheme to approximate its $L^2$-gradient flow. We illustrate numerically that the new flow approximates fairly well the mean curvature flow of codimension $1$ or $2$ interfaces in dimensions $2$ and $3$.
\end{abstract}

\maketitle

\section{Introduction}
A phase field approach was proposed in~\cite{MR4483340} to represent non-oriented interfaces and simulate their area-minimizing flow. The approach is based on ad hoc neural networks that are trained to approximate an evolution by mean curvature flow on a fixed time step. The main contribution of the current paper is to propose an alternative to neural networks based on a new but more standard variational phase field model. This model combines a Cahn-Hilliard-type energy to approximate the area and a Willmore-type energy to stabilize the phase field profile. The Cahn-Hilliard-type energy involves a nonstandard potential with one well at $0$ and an obstacle at $\frac 1 4$. The gradient flow associated with the proposed model is numerically similar to the mean curvature flow. We establish some analytical results which support, in part, the numerical observations.

Given $d\in\N^*$, $\e>0$ and $Q\subset\R^d$ an open and bounded set with Lipschitz boundary, we consider the phase field functional $\E_\e$ defined for every $u\in \wpq{2}{2}(Q)$ such that $u\leq \frac 1 4$ a.e. by
\[ \E_\e(u,Q) = \int_Q \Big(\frac{\e}{2}
|\nabla u|^2 + \frac{1}{\e} F(u)\Big) dx
+ \frac{\sigma_\e}{ 2 \e} \int_Q \left( \e \Delta u - \frac{1}{\e} F'(u) \right)^2dx,\]
where $\sigma_\e$ is any positive function of $\e$ such that $\e^2/\sigma_\e\to 0$ as $\e\to 0^+$, and the potential $F$  is nonsmooth and defined as
\[
F(s)= \begin{cases}
          s^2(\frac{1}{2}-2s)  & \text{if }  s \leq \frac 1 4 \\
          + \infty  &\text{ otherwise }
         \end{cases}
\]

The main theoretical contribution of this work is a convergence analysis of $\E_\e$ as $\e\to 0^+$
in dimension $d=1$, and in the radially symmetric case in higher dimension $d\geq 2$.
We also introduce a numerical scheme
to approximate the $L^2$-gradient flow of $\E_\e$
\[
\begin{cases}
 u_t &=  - \mu +   \sigma_\e  (\Delta \mu - \frac{1}{\e^2} F''(u) \mu), \\
 \mu &=  - \Delta u + \frac{1}{\e^2} F'(u),
\end{cases}
\]
and we illustrate with several numerical simulations that this new model provides a fairly good approximation of the mean curvature flow of non-oriented interfaces of codimension $1$ or $2$ in dimensions $d=2, 3$.

\medskip
In the case where the potential $F$ is replaced with a smooth and classical double-well potential such as $W(s)= \frac{1}{2} s^2 (1-s)^2$, it is well known that the $\Gamma$-limit of $\E_\e$ in dimensions $d=2, 3$ coincide on sets with smooth boundary with a weighted sum of their perimeter and Willmore energy~\cite{belpao93,Roger_schatzle_2006,MR3383330}.

In contrast with $W$, our potential $F$ has a well at $s= 0$ with $F'(0)=0$ and an obstacle at $s=\frac 1 4$ with $F'(\frac 1 4)<0$ (considering  the left-derivative of $F$ at $\frac 1 4$). If we consider only the first term of the functional, i.e.
\[\int_Q \Big(\frac{\e}{2} |\nabla u|^2 + \frac{1}{\e} F(u)\Big) dx,\]
the classical argument of Modica and Mortola~\cite{Modica1977} can be easily adapted to prove the $\Gamma$-convergence to the perimeter functional (up to a constant multiplier). However, adding the Willmore term changes everything, as it penalises the $\frac 1 4$-phase so that any minimizing sequence must converge to $0$ almost everywhere. In the limit, the $\frac 1 4$-phase has zero volume.

Other phase field models, such as in~\cite{Ambrosio_Tortorelli,BouchiDubsSepp96}, yield one phase only in the limit (up to a negligible set). For instance, the approximation of the Mumford-Shah functional by Ambrosio-Tortorelli's model~\cite{Ambrosio_Tortorelli} involves the single-well potential $G(s)=\frac{1}{2} (s-1)^2$ and an auxiliary phase field. In the limit, the $1$-phase connects with the (negligible) $0$-phase through a nonsmooth profile with infinite derivative at $0$. It is a major difference with our potential $F$ which guarantees smooth connecting profiles between phases $0$ and $\frac 1 4$. This smoothness allows us to use the Willmore energy, whose role is to stabilize the optimal phase-field profile, as will be demonstrated by our (partial) theoretical analysis and numerical experiments.

\medskip
As we will illustrate numerically, the main interest of the functional $\E_\e$, together with the proposed scaling of $\e$ and $\sigma_\e$, lies in the fact that the associated $L^2$-gradient flow provides a good numerical approximation of the mean curvature flow of non oriented interfaces, possibly with triple points. In addition, a key advantage of the phase-field approach is its versatility: the flow could be easily coupled with inclusion constraints, as in~\cite{MR4483340}, to approximate solutions such as those of the Steiner or Plateau problems.

\subsection{Classical phase field models for the perimeter and the Willmore energy}
It follows from the results of Modica and Mortola~\cite{Modica1977,Modica1977a} that, for a class of smooth double-well potentials $W$ that includes $W(s)= \frac{1}{2} s^2 (1-s)^2$, the $\Gamma$-limit in $\lp{1}(Q)$ of the so-called Cahn-Hilliard functional
$$P_\varepsilon(u)=\left\{\begin{array}{ll}\ds\int_Q \left(\frac\varepsilon 2|\nabla u|^2+\frac{W(u)}\varepsilon \right) dx&\mbox{if }u\in\wpq{1}{2}(Q)\\
+\infty&\mbox{otherwise in }\lp{1}(Q)\end{array}\right.$$
is $c_0 P(u)$ where
$$P(u)=\left\{\begin{array}{ll}\ds |Du|(Q)&\mbox{if }u\in\BV(Q,\{0,1\})\\
+\infty&\mbox{otherwise in }\lp{1}(\Omega)\end{array}\right.$$
and $c_0=\int_{0}^1\sqrt{2W(s)}ds$. In particular, if  $\Omega \subset\R^d$ has finite perimeter in the bounded set $Q$, thus $u:=\one{\Omega}\in \BV(Q,\{0,1\})$, one can build a sequence of functions $(u_\varepsilon)\in\wpq{1}{2}(Q)$ such that $u_\varepsilon\to u\in\lp{1}(Q)$ and $P_\varepsilon(u_\varepsilon)\to c_0|Du|(\Omega)=c_0P(\Omega,Q)$ as $\e\to 0$, with $P(\Omega,Q)$ the perimeter of $\Omega$ in $Q$. \medskip

A time-dependent smooth domain $\Omega(t) \subset\R^d$ evolves under the classical mean curvature flow if its inner normal velocity at every point $x\in\partial \Omega(t)$ is the scalar mean curvature $H_{\Omega(t)}(x)$ of $\partial \Omega(t)$ at $x$ (with the orientation convention that the scalar mean curvature on the boundary of a convex domain is positive).
This evolution may be understood as the $L^2$-gradient flow of the perimeter $P(\Omega(t))$. Now, $P$ can be approximated by the Cahn-Hilliard energy $P_\e$, and the $L^2$-gradient flow of $P_{\varepsilon}$ leads to the celebrated Allen-Cahn equation \cite{Ambrosio2000} which reads as, up to a time rescaling:

$$\partial_t u_{\e}  = \Delta u_{\e} - \frac{1}{\varepsilon^2}W'(u_{\e}).$$

A formal asymptotic expansion of the solution $u_{\varepsilon}$ to the Allen-Cahn equation near the interface $\partial \Omega_{\varepsilon}(t):=\partial \left\{u_\varepsilon(\cdot,t) \geq \frac{1}{2}\right\}$ gives (see~\cite{BellettiniBook})

\begin{equation*}
    u_{\varepsilon} (x, t) = q\left(\frac{\dist(x,\Omega_{\varepsilon}(t))}{\varepsilon}\right) + O(\varepsilon^2),
\end{equation*}
where
 $\dist(\cdot,\Omega_{\varepsilon}(t))$ is the signed distance function to $\Omega_{\varepsilon}(t)$ (negative inside, positive outside) and $q:\R\to [0,1]$ is the so-called \textit{optimal profile} which minimizes the parameter-free one-dimensional Allen-Cahn energy under some constraints:

$$
    q = \argmin_{p} \left\{ \int_{\R} \left(\frac{1}{2}|p'(s)|^2+{W(p(s))} \right) ds , ~ p\in C^{0,1}(\R), ~ p(-\infty) = 1, ~ p(0)=\frac 1 2,~ p(+\infty) = 0\right\}.
    $$
Furthermore, the velocity $V_\varepsilon$ of the boundary $\partial\Omega_{\varepsilon}(t)$ satisfies
$$
V_{\varepsilon}(t) = H_\varepsilon(t) + O(\varepsilon^2),
$$
where $H_\varepsilon(t)$ denotes the scalar mean curvature on $\partial \Omega_{\varepsilon}(t)$, which suggests that the Allen-Cahn equation approximates the mean curvature flow with an error of order $\varepsilon^2$~\cite{BellettiniBook}.

Rigorous proofs of convergence of the Allen-Cahn flow to the smooth mean curvature flow for short times (in particular before the onset of singularities)
have been presented in \cite{Chen1992, DeMottoniSchatzman1995, bellettini1995quasi}. More precisely, given an initial smooth set $\Omega_0$, its mean curvature flow $\Omega(t)$, an initial condition $u_0=q(\frac{\dist(x,\Omega_0)}\e)$, the solution $u_\e(\cdot,t)$ to the Allen-Cahn equation with $u_\e(\cdot,0)=u_0$, and the evolving set $\Omega_{\varepsilon}(t)=\left\{u_\varepsilon(\cdot,t) \geq \frac{1}{2}\right\}$, a quasi-optimal error on the Hausdorff distance between
$\Omega(t)$ and $\Omega_{\varepsilon}(t)$ is proved in these papers, namely
$$ \text{dist}_{\Hau}(\Omega(t),\Omega_{\varepsilon}(t)) \leq C \epsilon^2 |\log(\varepsilon)|^2,$$
where the constant $C$ depends on the regularity of $\Omega_0$.

The use of less regular potentials has been proposed, see for instance~\cite{Blowey_Elliott_1992},  to guarantee solutions with better physical or mathematical properties.
For example, the convergence analysis of the Allen-Cahn equation using a double obstacle potential
\[
W(s) = \begin{cases}
        s(1-s)  & \text{ if } s \in [0,1] ,\\
        +\infty &\text{ otherwise}
       \end{cases}
\]
is carried out in \cite{MR1290089,Elliott_Schatzle_1996}, where it is proved that solutions,
which are constrained to remain valued in $[0,1]$,
have finite transition zone around the interface $\partial \Omega_\e$ such that $ \text{dist}_{\Hau}(\Omega(t),\Omega_{\varepsilon}(t)) \leq C \epsilon^2$. 

The use of logarithmic potentials in Cahn-Hilliard models for the approximation of surface diffusion flows has been extensively studied in \cite{MR1290089}, allowing also to guarantee the inclusion property of the phase field solution $u_{\e}$ and no singularities.
\medskip

The Willmore energy of  a set $\Omega\subset\R^d$ with smooth boundary in~$Q$ is defined as
$$\cW(\Omega,Q)= \frac 1 2\int_{\partial \Omega \cap Q} |H_\Omega(x)|^2~d\haushd,$$
where $\haushd$ denotes the $(d-1)$-dimensional Hausdorff measure. Based on a conjecture of De Giorgi~\cite{degiorgi-a}, several authors~\cite{belpao93,tone02,belmug04,moser,Roger_schatzle_2006,naga-tone,MR3383330,MR4575399}
 have investigated the phase field approximation of $\cW$ using a double-well potential such as the function $W$ defined in the introduction.
 A classical approximation functional is
\[
\cW_\varepsilon(u)=\left\{\begin{array}{ll}\ds \frac 1{2\varepsilon}\int_Q \left( \varepsilon\Delta u-\frac{W'(u)}{\varepsilon} \right)^2dx&\mbox{if }u\in\lp{1}(Q)\cap\wpq{2}{2}(Q)\\
+\infty&\mbox{otherwise in }\lp{1}(Q).\end{array}\right.
\]

\par Bellettini and Paolini~\cite{belpao93} have proved for this functional a $\Gamma$-limsup property: given a set $\Omega$ with smooth boundary, there exists a suitable family of smooth approximations $u_\e$ of $\one{\Omega}$ such that the limit of $\cW_\varepsilon(u_\varepsilon)$ is,  up to a multiplicative constant, the Willmore energy $\cW(\Omega,Q)$.
\par The $\Gamma$-liminf property is much harder to prove, see the various contributions in~\cite{tone02,belmug04,moser,Roger_schatzle_2006,naga-tone}. In particular, it was proven by R{\"o}ger and Sch{\"a}tzle in~\cite{Roger_schatzle_2006} in space dimensions $N=2,3$ and, independently, by Nagase and Tonegawa~\cite{naga-tone} in dimension $N=2$, that the result holds true for smooth sets. More precisely, given $u=\one{\Omega}$ with $\partial \Omega\in\cont^{2}(Q)$, and $u_\varepsilon$ converging to $u$ in $\lp{1}$ with a uniform control of $P_\varepsilon(u_\varepsilon)$, then
$$c_0 \mathcal{W}(\Omega,Q)\leq\liminf_{\varepsilon\to 0^+} \cW_\varepsilon(u_\varepsilon).$$
The proof is based on a careful control of the discrepancy measure
$$\xi_\varepsilon= \left(\frac\varepsilon 2|\nabla u|^2-\frac{W(u)}\varepsilon\right){\mathcal L}^2,$$
which guarantees that the minimizing sequence has the correct phase field profile
$$u_\e \simeq q(\dist(x,\Omega)/\e).$$
Then, in R{\"o}ger \& Sch{\"a}tzle's proof, the result follows from a representation with varifolds of the limit measure, a lower semicontinuity argument and the locality property of the generalized mean curvature of integral varifolds. Due to dimensional requirements for Sobolev embeddings and for the control of singular terms, the proof works in dimensions $N=2, 3$ only. The result in higher dimension is still open. \medskip

It was shown, first formally in \cite{MR1757511,MR2377282,MR3383330} with asymptotic expansion arguments, then rigorously in~\cite{rigorousWillmoreconvergence}, that the $L^2$-gradient flow of $\cW_{\e}$,
\[
\left\{
\begin{aligned}
\partial_t u_\e &= \Delta \mu_\e  - \frac{1}{\e^2} W''(u_\e) \mu_\e, \\
\mu_\e &=  - \Delta u_\e + \frac{1}{\e^2} W'(u),
\end{aligned}
\right.\]
converges to the Willmore flow, i.e. the $L^2$-gradient flow of $\cW$, which is characterized by the inner normal velocity
$$ V_n =  \Delta_{\partial E} H +  \|A\|^2H  - \frac{H^3}{2},$$
where $H$ is the mean curvature on $\partial E$, $A$ is the second fundamental form and $\|A\|^2$ coincides with the sum of the squared principle curvatures. \medskip

We are not aware of any use of potentials less regular than $W(s)=s^2(1-s)^2$ for the approximation of the Willmore energy and its associated flow, the lack of regularity being an obstacle to any proof of convergence. The idea of the present article is therefore to go a little deeper into this case. We will show the interest of considering a potential with one well and one obstacle to both bound and stabilize the phase associated to the obstacle.

\subsection{Ambrosio-Tortorelli's functional and the case of non oriented interfaces}
Some applications, typically free-discontinuity problems, require the minimization of the area of interfaces that are not domain's boundaries and sometimes even not orientable. In such a case, obviously, the Cahn-Hilliard energy cannot be used to approximate the area. \medskip

A first typical example is in image segmentation with the Mumford-Shah model~\cite{Mumfordshah,David_Guy_book}. Given a gray-scale image $I_0$ defined on a domain $Q\subset\R^2$, one seeks for a
piecewise smooth approximation $v$ of $I_0$ by minimizing the Mumford-Shah functional
$$ \E^{MS}(v,\Gamma)  = \int_{Q} (v - I_0)^2 dx +  \mu \int_{Q\setminus \Gamma} |\nabla v|^2 dx + \eta  \Hau^{1}(\Gamma).$$
 among all pairs $(v,\Gamma)$ such that
 $$\Gamma \subset Q \text{ is closed  and } v \in C^1(Q \setminus \Gamma).$$

Ambrosio and Tortorelli proposed in \cite{Ambrosio_Tortorelli} a phase-field approximation of $\E^{MS}$ of the form
$$  \E^{AT}_{\varepsilon}(v,u) =  \int_{Q} (v - I_0)^2 dx +  \mu \int_{Q} u^2 |\nabla v|^2 dx + \eta \mathcal{F}_{\varepsilon}(u)$$
where the so-called Ambrosio-Tortorelli term $\mathcal{F}_{\varepsilon}$ is defined as
$$\mathcal{F}_{\varepsilon}(u) = \int_Q \left(\varepsilon |\nabla u| + \frac{1}{4 \varepsilon} (1-u)^2 \right)dx.$$
Ambrosio and Tortorelli proved a $\Gamma$-convergence result of $\E^{AT}_{\varepsilon}$  to $\E^{MS}$ in a suitable functional framework.
Intuitively, a phase field minimizer $u_{\varepsilon}$ of $\E^{AT}_{\varepsilon}$ is of the form $u_{\varepsilon} = \varphi(\dist(x,\Gamma)/ \varepsilon)$ (where $\dist$ now denotes the classical distance function), it vanishes on $\Gamma$ and $\mathcal{F}_{\varepsilon}(u_\e)$ approximates the length of $\Gamma$. However, unlike the Cahn-Hilliard energy, the term $\mathcal{F}_{\epsilon}$ considered separately $\Gamma$-converges to $0$ rather than to the length of $\Gamma$. In other words, $\mathcal{F}_{\epsilon}$ needs to be coupled with additional terms and constraints to approximate the length of $\Gamma$, which raises a number of numerical difficulties. \smallskip

A second example is related to the Steiner problem. Recall that, given a collection of points $a_1, \ldots, a_L \in Q$, the associated Steiner problem consists in finding a compact connected set $\Gamma\subset Q$ containing every $a_i$ and having minimal length.
It is equivalent to finding the optimal solution to the problem
	\[
	\min_\Gamma\left\{\Hau^1(\Gamma),~ \Gamma\subset Q, \;\Gamma \text{ connected, } a_i\in \Gamma, ~\forall i=1,\ldots, L \right\},
	\]

Ad hoc variational phase field models for the Steiner problem have recently been introduced, see for instance \cite{MR3337998,MR4011534,BonafiniOrlandiOudet,BonafiniOudet,ChambolleFerrariMerlet2019-1}. The model introduced in \cite{MR3337998}
involves the Ambrosio-Tortorelli term  ${\mathcal F}_{\varepsilon}$
coupled to an additional geodesic term ${\mathcal G}_{\varepsilon}$ that forces the connectedness of$\Gamma$:
$$
{\mathcal G}_{\varepsilon}(u) = \frac{1}{\varepsilon} \sum_{i=1}^{N} {\bf D}(u^2; a_0 , a_i), \;\; \text{ with } {\bf D}(w;a,b):=\inf_{\Gamma: a \leadsto b}\int_{\Gamma} w\,d{\mathcal H}^1 \in [0,+\infty].
$$
Here the notation $\Gamma: a \leadsto b$ refers to a rectifiable curve $\Gamma\subset\Omega$ connecting $a$ and $b$.

Intuitively, and as previously, a phase field minimizer $u_{\varepsilon}$ of $\mathcal{F}_{\varepsilon}+\mathcal{G}_{\varepsilon}$ is expected to be of the form $u_{\varepsilon} = \varphi(\dist(x,\Gamma)/ \varepsilon)$ where the geodesic term $\mathcal{G}_{\varepsilon}(u_\e)$ forces $u_\varepsilon$ to vanish on a set $\Gamma$ connecting all $a_i$'s and whose length is approximated by the Ambrosio-Tortorelli term $\mathcal{F}_{\varepsilon}(u_{\varepsilon})$.\\

The $\Gamma$-convergence result proved in \cite{MR3337998} suggests that the minimization of this phase field model
provides an approximation of a Steiner solution, i.e. a Steiner tree.  However, despite the numerical experiments  provided in  \cite{MR3337998,MR4011534} show the ability of the model to approximate solutions of the Steiner problem in dimensions $2, 3$, they also show how difficult it is to minimize effectively  the geodesic term $\mathcal{G}_{\varepsilon}(u_\e)$ to guarantee the connectedness of $\Gamma$.
The conclusions are quite similar for the approaches developed in \cite{BonafiniOrlandiOudet,BonafiniOudet,ChambolleFerrariMerlet2019-1}
based on the measure-theoretic notion of current. In these latter approaches, the connectedness of $\Gamma$
is ensured by adding a divergence constraint of the form  $ \div(\tau) = \sum \alpha_i \delta_{a_i}$ for a suitable field $\tau$. 

In short, all these models provide reasonable approximations of solutions to the Steiner problem, but the algorithmic cost is rather high and the approximations lack of regularity. \smallskip

A motivation of this paper is therefore to propose a more straightforward approach that does not need to be coupled with sophisticated constraints, and which is directly related to the approximation of the mean curvature flow of non oriented interfaces.

\subsection{Unoriented phase field profile and neural networks}

As mentioned above, the key point that limits the effectiveness of phase field methods for approximating minimization problems involving the perimeter of a non oriented set $\Gamma$ is that they are based on Ambrosio-Tortorelli-type energies which need  to be coupled with other terms or constraints to ensure the stability of the optimal phase field profiles. In fact, unlike the Allen-Cahn equation, i.e. the gradient flow of the Cahn-Hilliard energy, which provides good approximations of the mean curvature flow of set boundaries, the gradient flow of the Ambrosio-Tortorelli model is not suitable for approximating the evolution by mean curvature of more general interfaces, in particular interfaces that are not boundaries.

\begin{figure}[!htb] \label{fig_phase_field representation}
\begin{center}
\includegraphics[width=0.30\textwidth,height=0.26\textwidth]{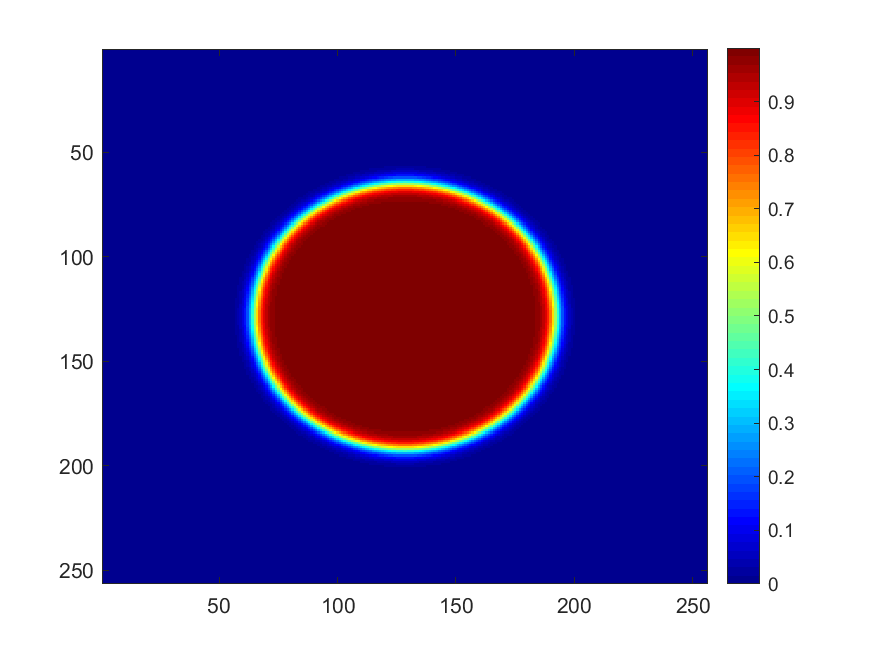}\hspace*{0.5cm}
\includegraphics[width=0.30\textwidth,height=0.26\textwidth]{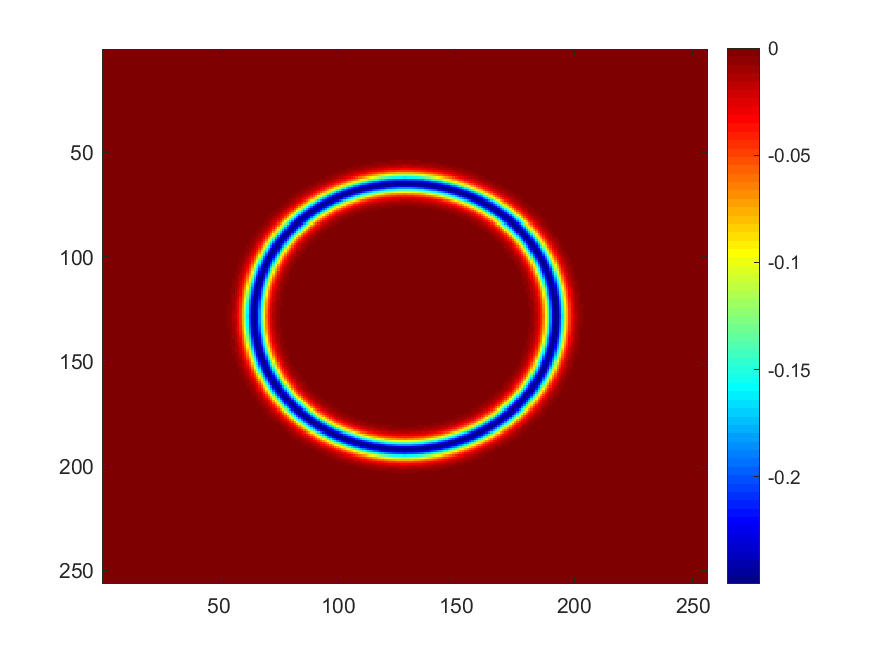}
\caption{Oriented ($q$, left) and non-oriented ($-q'$, right) phase field approximations of a disk / circle in dimension $2$.}
\end{center}
\end{figure}

Recently, it was proposed in \cite{MR4483340} to design and train suitable neural networks for approximating these latter evolutions.
More specifically, using the same notations as above and denoting as $\Gamma(t)$ the mean curvature flow of a smooth, general interface (not necessarily the boundary of a domain), the general idea in \cite{MR4483340} is to train a neural network $\S_{\theta}$ such that
$$ \S_{\theta}\left[-q'\left( \frac{\operatorname{d}(x,\Gamma(t))}{\e} \right) \right] \simeq  -q'\left( \frac{\operatorname{d}(x,\Gamma(t + \delta_t))}{\e} \right),$$
where $d$ is now the (unoriented) distance function. The main advantage of using the profile $q'$ is that it does not require any orientation of $\Gamma(t)$. Figure~\ref{fig_phase_field representation} illustrates the difference between oriented and non-oriented representations.

The training of $\S_{\theta}$, i.e. the calibration of its inner parameters, was performed in~\cite{MR4483340} with a training database made of circles (in dimension $2$) or spheres (in dimension $3$) whose flow by mean curvature is known exactly.

\begin{figure}[!htbp] \label{fig_phase_field neural_network_structure}
\begin{center}
\includegraphics[width=5cm]{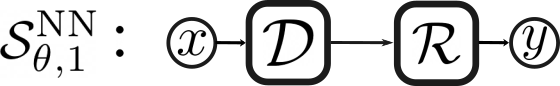}\hspace*{1cm}
\includegraphics[width=6cm]{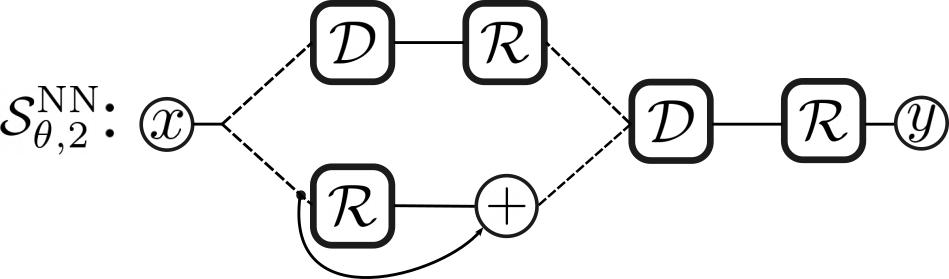}
\caption{Two network structures obtained from first order (left) and second order (right) discretisation schemes of the Allen-Cahn equation}
\end{center}
\end{figure}

As for the network structures proposed for $\S_{\theta}$  in~\cite{MR4483340}, they are very simple and inspired by discretisation schemes that alternate diffusion and reaction operations to approximate numerically the solutions to the Allen-Cahn equation. Two typical structures studied in~\cite{MR4483340} are shown in figure \ref{fig_phase_field neural_network_structure}, where each ${\mathcal D}$ denotes a convolution operation (two different ${\mathcal D}$ in the same network might have different parameters) and ${\mathcal R}$ is a pointwise reaction trainable function (in practice, a 3-layer perceptron).

\begin{figure}[!htbp] \label{fig_MMC_Neural_network}
\begin{center}
\includegraphics[width=.24\textwidth]{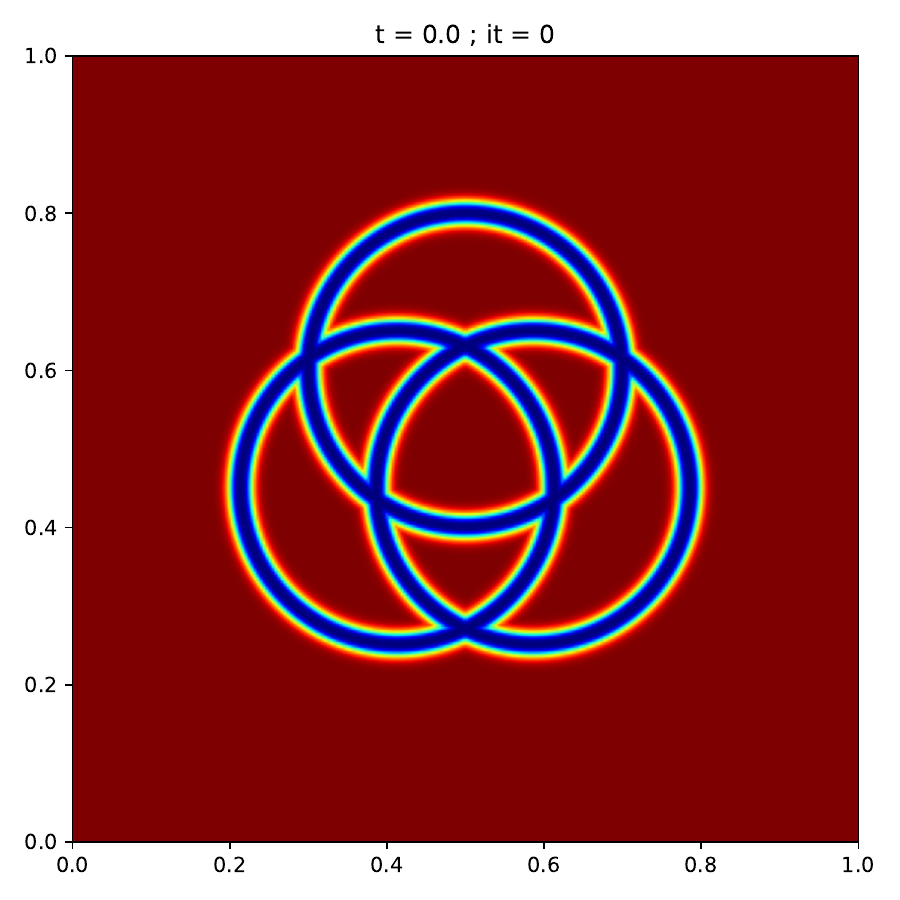}
\includegraphics[width=.24\textwidth]{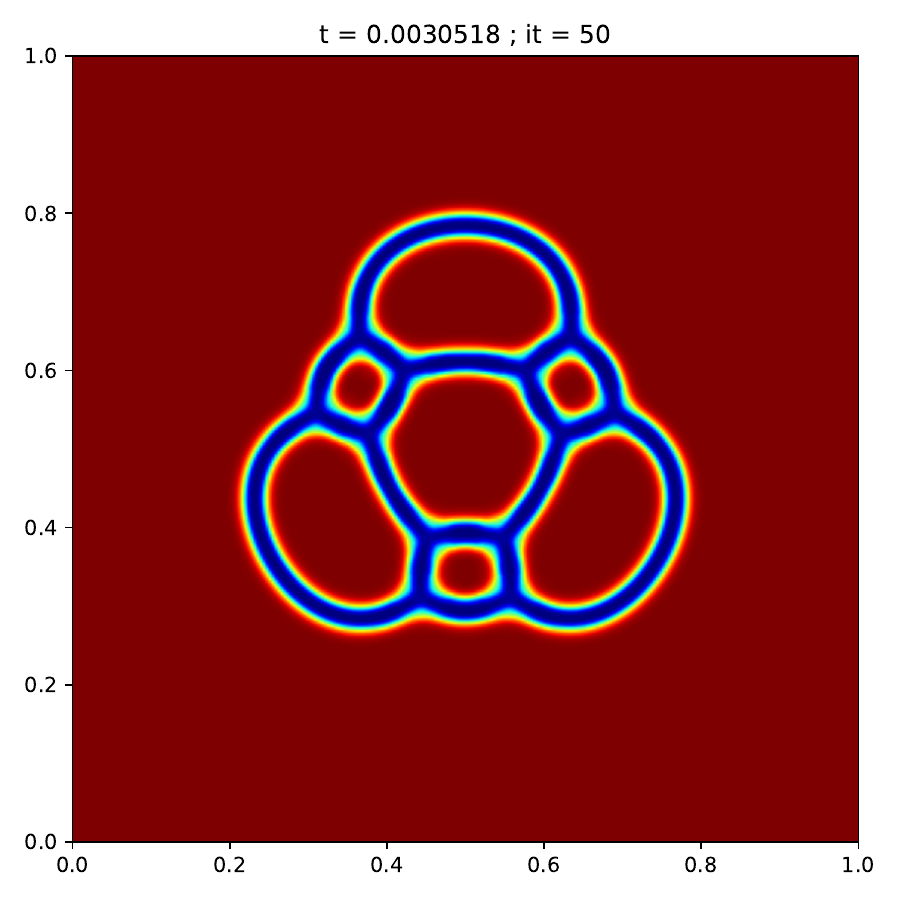}
\includegraphics[width=.24\textwidth]{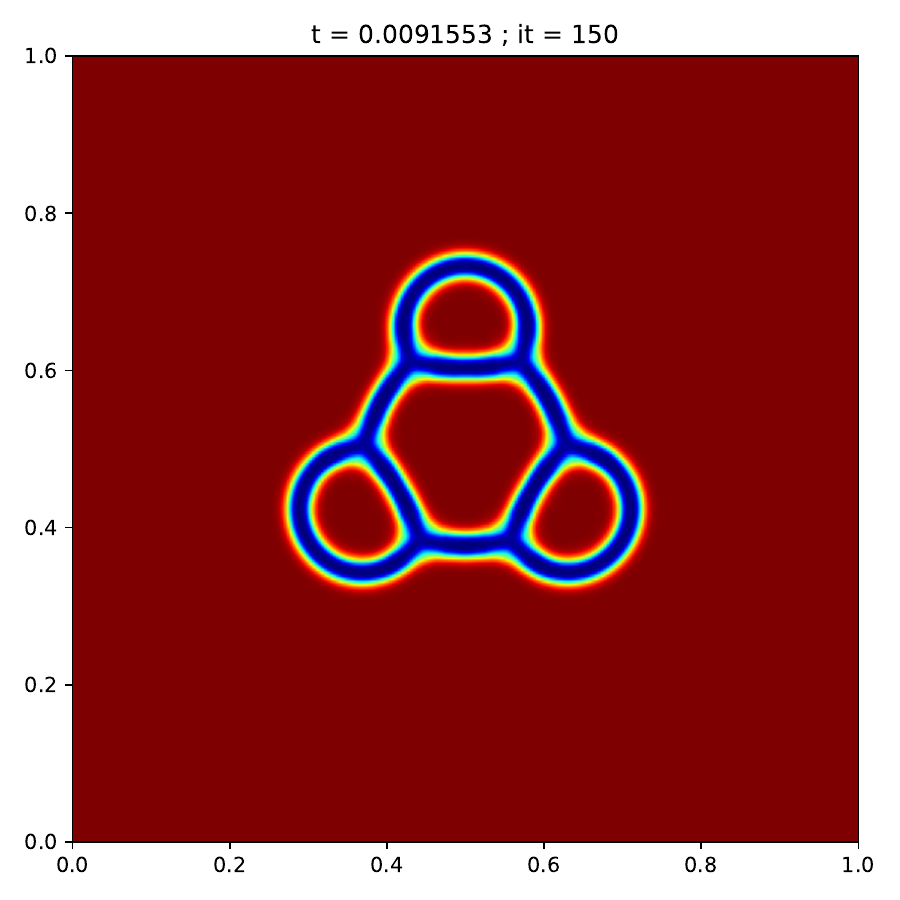}
\includegraphics[width=.24\textwidth]{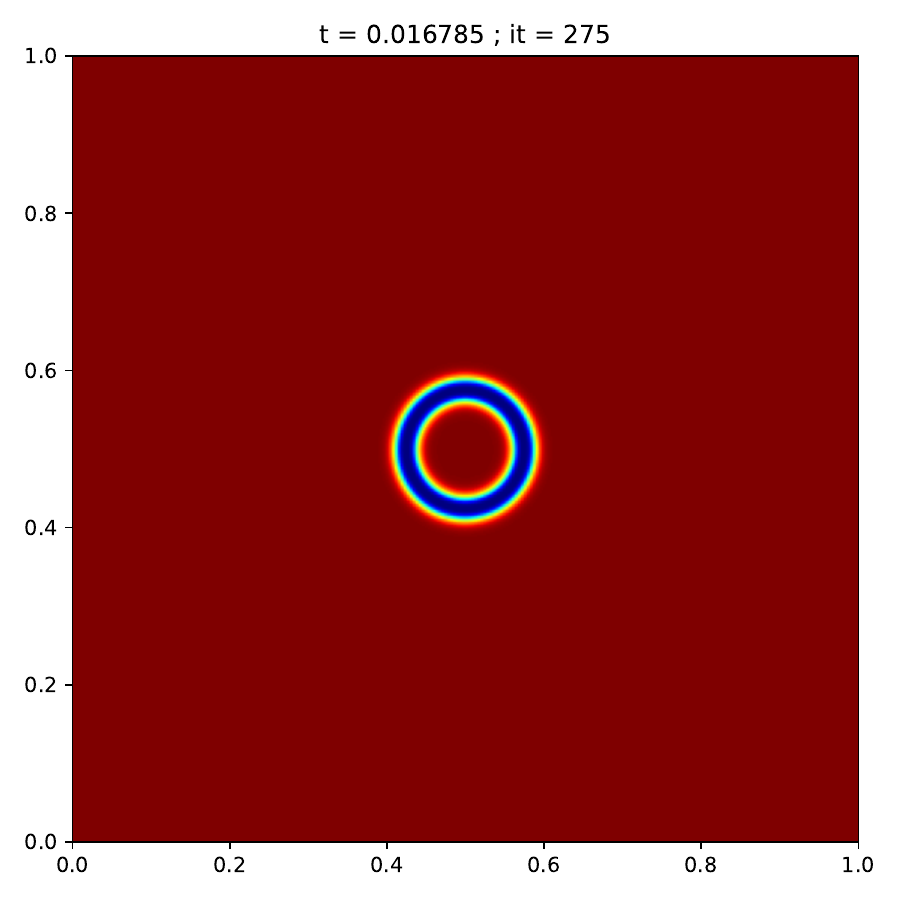}

\caption{Approximation with a second order neural network of the motion by mean curvature of non oriented interfaces.}
\end{center}
\end{figure}

The evaluation of the numerical performances of these networks leads to the following conclusion: oriented mean curvature flows can be learned with either first-order or second-order networks, but only second-order networks are able to learn non oriented mean curvature flows. More precisely, the first-order networks proposed in~\cite{MR4483340} extend time-discrete approximations of the Allen-Cahn equation
$$u_t=\Delta u-\frac 1{\e^2}W'(u),$$
hence the results of~\cite{MR4483340} suggest that any phase field model of the form 
\[ u_t =  D_1 [u] + R_1[u]\]
fails to approximate a non-oriented mean curvature flow.
 Instead, the good performances of second-order networks (as illustrated in figure~\ref{fig_MMC_Neural_network}) suggest that a phase field model such as
\begin{equation}\label{2ndorder}
\begin{cases}
    u_t &=   D_2 [v] + R_3(u) R_2 [v],   \\
    v &=  D_1 [u] + R_1[u],
   \end{cases}
\end{equation}
might be able to approximate non-oriented mean curvature flows.

\subsection{Derivation of a new phase field model}

Following these observations, let us try to identify a phase field model that would be suitable for the approximation of a nonoriented mean curvature flow. We start with the derivation of an ODE satisfied  by the profile $-q'$. Recall that $q$ satisfies
$$ q'(s) = - \sqrt{2 W(q(s))}, \quad   q''(s) = W'(q(s))\quad \text{ and } q(0) = \frac{1}{2},$$
where $W$ is the quadratic smooth potential  $W(s) = \frac{1}{2} s^2 (1-s)^2$. Exploiting the symmetry of the double-well potential $W$, it is not difficult to see that
\begin{eqnarray*}
 q^{(3)}(s) = W''(q(s)) q' =  ( 6 q(s)^2 - 6 q(s)  +   1 ) q'(s) = ( 1 + 6 q'(s)) q'(s).
\end{eqnarray*}
This shows  that the profile  $y = -q'$ is a solution to the equation
$$ y''(s) = F'(y(s)),$$
where the potential $F$ is defined as previously by
 $$ F(s)= \begin{cases}
          s^2(\frac{1}{2}-2s)  & \text{if }  s \leq \frac 1 4 \\
          + \infty  &\text{ otherwise }.
         \end{cases}
$$

A PDE naturally associated with the above ODE is 
the Allen-Cahn-type equation
$$u_t = \Delta u - \frac{1}{\e^2} F'(u),$$
i.e. the $L^2$-gradient flow of the modified Cahn-Hilliard energy obtained by replacing with $F$ the smooth potential $W$.
However, as shown by the numerical experiments performed with neural networks, such a model fails to preserve the profile $q'$ along the iterations.

A simple way to get stability consists in adding to the modified Cahn-Hilliard energy a second-order term that forces the equation $y'' = F'(y)$. More precisely, we consider the following phase field approximation of the perimeter plus Willmore energy:
\[
\E_\e(u) = \int_Q \Big(\frac{\e}{2} |\nabla u|^2 + \frac{1}{\e} F(u)\Big) dx +   \frac{\sigma_\e}{ 2 \e} \int_Q \left( \e \Delta u - \frac{1}{\e} F'(u) \right)^2dx
\]
where the parameter $\sigma_{\varepsilon}$ must be large enough to guarantee the profile's stability, but not so large that the Willmore term remains negligible
compared to the perimeter term. The $L^2$-gradient flow of $\E_\e$ writes
$$
\begin{cases}
 u_t &=  - \mu +   \sigma_\e  (\Delta v- \frac{1}{\e^2} F''(u)v), \\
 v &=  - \Delta u + \frac{1}{\e^2} F'(u),
\end{cases}
$$
which is exactly of the proposed form~\eqref{2ndorder}.

\subsection{Outline of the paper} The convergence analysis of $\E_\e$
as $\e\to 0$ is carried out in the next section. By restricting ourselves to the case of radial symmetry in dimension $d$, we will demonstrate a $\Gamma$-convergence result in the case of a relaxed version of $\E_{\e}$ which will make it possible to identify a limit set $E$ and for which $\E_{\e}(u_\e)$ gives an approximation of the sum of the perimeter and the Willmore energy of $E$.
\par Section~\ref{sec:numeric} focuses on the numerical discretization of the gradient flow of $\E_\e$, and in particular on its ability to approximate the mean curvature flow of non oriented interfaces. We illustrate this with various numerical simulations in dimensions $2$ and $3$.

\section{Convergence analysis}\label{sec:cva}

Let $Q$ be an open set of $\R^d$ and  consider the functional $\E_\e$ defined on $L^2(Q)$ as
\[
  \E_\e(u) =
  \begin{cases}\displaystyle
  \int_{Q}   {\bf m}_\e(u)   dx
  +  \frac{\sigma_\e }{2\e}  \int_{Q}  {\bf eul}_\e(u)^2
  dx & \text{ if } u\in \wpqO{2}{2}(Q),\;\;0\le u \le \frac{1}{4} \;\;\text{ a.e.,}
  \\
  +\infty & \text{ otherwise,}
\end{cases}
\]
with ${\bf m}_\e(u) = \frac{\e}{2}|\nabla u|^2+\frac{1}{\e}F(u)$ and ${\bf eul}_\e(u) = - \e \Delta u +  \frac{1}{\e}F'(u)$.

As explained previously, although this phase field model seems very close to approximation models for the area plus Willmore energy involving a quadratic potential $W(s)= \frac{1}{2}s^2(1-s)^2$, the analysis here is very different because of the particular choice of $F$ which admits an obstacle and not a double well in $s=1/4$. The first difficulty is that any minimizing sequence $(u_{\e})_\e$ of $\E_\e$ converges a.e. to $0$, hence does not allow us, in contrast with the classical Cahn-Hilliard approximation, to identify as the boundary of a limit domain a limit set $\Gamma$ where the area plus Willmore energy concentrates. Intuitively, this follows from  $F$ having a non vanishing derivative at $s = 1/4$, making the phase $\{ s = 1/4 \}$ unstable. An alternative to identify the limit set is to consider the support of the limit measure of diffuse perimeter measures.

Given a sequence  $(u_\e)_\e$ satisfying $\E_\e(u_\e)\le C<+\infty$, we define as in~\cite{MR1425577,moser,naga-tone,Roger_schatzle_2006} the absolutely continuous measures
\[
  \mu_\e=\left(\frac{\e }{2}|\nabla u_\e|^2+\frac{1}{\e}F(u_\e)\right){\mathcal L}^d.\]

There exists a Radon  measure $\mu$
such that, possibly passing to a subsequence,  
\begin{center}
$\mu_\e$  converges weakly-$*$ to $\mu$ as $\varepsilon \to 0$.
\end{center}

We study hereafter the properties of the support of $\mu$ and its relationship with the limits of $u_\e$ and  $\E_{\e}(u_{\e})$. We consider the problem only in the radially symmetric case in dimension $d \geq 1$. More precisely, in the case where $Q = B_{R_B} \setminus  B_{R_A}$ with $R_B > R_A > 0$, we prove  in Theorem \ref{thmdd} that the limit measure $\mu$ can be decomposed as
$$\mu= c_F  \sum_{i\in I}  m_i\haushd\res{\partial B_{r_i}},$$
where $c_F = 2\int_{0}^{1/4} F(s) ds$, $I$ is finite, and for all $i\in I$, $B_{r_i}$ is the ball of radius $r_i>0$ centered at $0$ and $m_i$ is a positive integer multiplicity. We will also prove that
\begin{equation}\label{decomp-energy}
\liminf_{\e \to 0} \E_{\e}(u_{\e}) \geq   c_F \sum_{i\in I} m_i \left[ \mathcal{P}(B_{r_i}) + \sigma_0 \mathcal{W}(B_{r_i}) \right]
\end{equation}
where $\sigma_0 = \lim_{\e \to 0} \sigma_{\e}$. This result can be interpreted as a $\Gamma$-liminf with a suitable relaxation. Consider indeed the following functional on $L^1(Q)\times\cM(Q)$ that depends on both the phase field function $u_{\e}$ and the associated measure $\mu_{\e}$:

$$
\mathcal{F}_\e(u_\e,\mu_\e) =  \begin{cases}
                   \E_\e(u_\e) & \text{ if } (u_\e,\mu_\e) \in \mathcal{D}_{\e} \\
                   + \infty & \text{ otherwise}
                  \end{cases}
$$
where  $\mathcal{D}_{\e}$ is the collection of all couples $(u_\e,\mu_\e)\in L^1(Q)\times\cM(Q)$ such that $u_\e\in \wpqO{2}{2}(Q)$ is radially symmetric and $\mu_\e={\bf m}_\e(u_\e){\mathcal L}^d=\left(\frac{\e }{2}|\nabla u_\e|^2+\frac{1}{\e}F(u_\e)\right){\mathcal L}^d$. Considering the $L^1$-topology for $u_{\e}$ and the topology of weak-* convergence for the Radon measure $\mu_{\e}$, inequality~\eqref{decomp-energy} above implies that the $\Gamma$-liminf of $\mathcal{F}_{\e}$ on $L^1(Q)\times\cM(Q)$ coincides with
$$
\mathcal{F}(u,\mu) =  \begin{cases}
               c_F \sum_{i\in I} m_i \left[P(B_{r_i}) + \sigma_0 \mathcal{W}(B_{r_i}) \right] & \text{ if } (u,\mu) = (0, c_F  \sum_{i\in I}  m_i\haushd\res{\partial B_{r_i}} ), \text{$I$ finite} \\
             + \infty &\text{ otherwise.}
            \end{cases}
$$
We can actually prove a stronger result:
\begin{theo}\label{theo-Gamma-conv} Let $Q = B_{R_B} \setminus  B_{R_A}$ with $R_B > R_A > 0$. Then $\mathcal{F}_\e$ $\Gamma$-converges to $\mathcal{F}$ on $L^1(Q)\times\cM(Q)$.

\end{theo}

\par The proof is organized as follows: in Section~\ref{sec:bluk_control} we show that, thanks to the Willmore term, the mass energy of the support of $\mu$ can be localized in sets where $u_{\e} > \eta$ for some $\eta > 0$. Then we develop in Sections~\ref{sec:Gammaliminf_1} and~\ref{sec:Gammaliminf_d}, respectively, the analysis of the $\Gamma$-liminf of $\mathcal{F}_{\e}$ in dimension $1$, and the analysis in the radially symmetric case in dimension $d \geq 2$. Although the radial $d$-case actually includes the result in dimension $d=1$, we think the exposition is clearer by separating the cases. Lastly, the $\Gamma$-limsup of $\mathcal{F}_{\e}$ is proved in Section~\ref{sec:Gammalimsup}.

\subsection{Localization of the support of $\mu$} \label{sec:bluk_control}
A first result, directly adapted from \cite{belmug04,MR1425577}, consists in showing, thanks to the contribution of the Willmore term, that in the neighborhood of the support of $\mu$, $u_\e$ does indeed involve a diffuse interface.
More precisely, we prove that for any $x_0$ in the support of $\mu$, there
exists, for some fixed $\eta>0$,  a subsequence $x_\e \to x_0$ such that $u_\e(x_\e)>\eta$.

  \begin{lemma} \label{lemme_xi}  Let $Q = \R^d$ and  $(u_\e)_\e$ be a sequence of functions in $\wpq{2}{2}(Q)$ such that $\E_\e(u_\e)\le C<+\infty$ and $ \mu_\e=\left(\frac{\e }{2}|\nabla u_\e|^2+\frac{1}{\e}F(u_\e)\right){\mathcal L}^d$  converges weakly-$*$ to a Radon measure $\mu$ as $\varepsilon \to 0^+$. Let $x_0 \in \spt(\mu)$, $r>0$
 and $\eta \in [0,\eta_0/2]$ where $\eta_0 = 1/12$. Then there exists $\overline{\e}$ such that
 $$ B_r(x_0) \cap \{ u_\e > \eta\} \neq \emptyset,   \quad \forall \e \in ]0,\overline{\e}[. $$
 \end{lemma}

 \begin{proof}
 The proof is a minor modification of the proof of Lemma $4.4$ in \cite{belmug04} adapted
 to our potential $F$.  We exploit some properties of $F$, in particular  that $F''(s)=1-12s \geq 0$ for all
 $s \leq  \eta_0 = 1/12$  as well as the existence of
 $c>0$ such that $F(s) \leq c F'(s)^2$ for all $s\leq \eta_0/2$. \\

 Let $A_1,A_2$ be two open subsets of $\R^d$ with $A_1 \subset \subset A_2 \subset \subset  \R^d$.
Let us prove  that there exists a positive constant $C_0$ such that
\begin{multline}\label{equ_control_bulk}
  \int_{A_1 \cap \{ u_\e \le \eta \}} [ {\bf m}_\e(u_\e) + \frac{1}{\e} F'(u_\e)^2] dx \\\leq   C_0 \eta
  \int_{A_2 \cap \{ u_\e > \eta \}} [\e |\nabla u_\e|^2 ]dx + C_0 \e \int_{\R^d} ({\bf eul}_\e(u_\e))^2 dx 
  +  C_0 \e^{1/2} \left[ \int_{\R^d}  {\bf m}_\e(u_\e) dx \right]^{1/2}
\end{multline}
for $\e$ sufficiently small. \medskip

Let $\phi \in C^{\infty}_{c}(A_2;[0,1])$ be such that $\phi = 1$ on $A_1$. Consider the function
$g$ defined by $g(s) = F'(s)$ for all $s \leq \eta$, $g(s) = 0$ for all $s>\frac{1}{12}$ and  $g$ affine
on $[\eta,\frac{1}{12}]$. In particular, $g'(s)<0$ for $\eta<s<\frac{1}{12}$, and
 $0 \leq g^2(s) \leq F'(s) g(s)$ for all $s \in \R$.
{}Using the identity
\[ \int_{\R^d} \phi ~ {\bf eul}_\e(u_\e)  g(u_\e)  dx = \int_{\R^d}  \left[ \e g'(u_\e) |\nabla u_\e|^2 \phi +  \frac{1}{\e} F'(u_\e) g(u_\e) \phi  +  \e g(u_\e)  \nabla u_\e \cdot \nabla \phi \right] dx,\]
and observing that
\begin{eqnarray*}
 \int_{\R^d} \phi ~ {\bf eul}_\e(u_\e)  g(u_\e) dx &\leq& \int_{\R^d} \frac{1}{2} \phi \left[ \e {\bf eul}_\e(u_\e)^2 +  \frac{1}{\e}g(u_\e)^2 \right]  dx \\
 &\leq& \int_\R^d \frac{1}{2} \phi \left[ \e {\bf eul}_\e(u_\e)^2 +  \frac{1}{\e}F'(u_\e)g(u_\e) \right]  dx,
\end{eqnarray*}
we obtain that
\[\int_{\R^d} \phi \left[ \e g'(u_\e) |\nabla u_\e|^2  \right] dx \leq  \int_{\R^d}  \phi \frac{1}{2} \left[ \e {\bf eul}_\e(u)^2  -  \frac{1}{\e}F'(u_\e)g(u_\e) \right] dx   - \int_{\R^d} \e g(u_\e)  \nabla u_\e \cdot \nabla \phi  dx.  \]
Then,
\begin{eqnarray*}
 \phantom{I}&\phantom{=}&  \int_{ \{ u_\e \leq \eta  \} \cap A_1} \left[ \e F''(u_\e) |\nabla u_\e|^2 + \frac{1}{2 \e} F'(u_\e)^2 \right] dx \\ & \leq &
 \int_{ \{ u_\e \leq \eta \}} \phi
 \left[ \e g'(u_\e) |\nabla u_\e|^2 + \frac{1}{2 \e} F'(u_\e)g(u_\e) \right] dx \\
 &\leq&
 \frac{\e}{2} \int_{\R^d} {\bf eul}_\e(u_\e)^2 dx - \int_{ \{ u_\e > \eta \}} \phi \left[ \e g'(u_\e) |\nabla u_\e|^2 \right]   - \int_{\R^d} \e g(u_\e)  \nabla u_\e \cdot \nabla \phi  dx,
 \\ &\leq&
 \frac{\e}{2} \int_{\R^d} {\bf eul}_\e(u_\e)^2 dx  - \int_{ \{ u_\e > \eta \} \cap A_2}  \left[ \e g'(u_\e) |\nabla u_\e|^2 \right] + \left| \int_{A_2} \e  g(u_\e)  \nabla u_\e \cdot \nabla \phi  dx \right|.
\end{eqnarray*}
Using Hölder's inequality, we obtain that there exists a constant C depending on $A_2, \phi$ and $F$ such that
\begin{eqnarray*}
  \left| \int_{A_2} \e  g(u_\e)  \nabla u_\e \cdot \nabla \phi dx \right| &\leq & \e^{1/2} \|\nabla \phi\|_{\infty} \sup_{[0,1/4]} |g| ~ |A_2 \setminus A_1|^{1/2}
  \left[ \int_{\R^d} \e | \nabla u_\e |^2 dx\right]^{1/2} \\
  &\leq& C\varepsilon^{1/2}   \left[ \int_{\R^d}  {\bf m}_\e(u_\e) dx \right]^{1/2}.
\end{eqnarray*}
Finally, as for all $s < \eta$ $ F''(s) \geq c_0$ and $ F(s) \leq c_1 F'(s)^2$, we deduce that there exists $C>0$ such that
  \[
  \int_{A_1 \cap \{ u \le \eta \}} [ {\bf m}_\e(u_\e) + \frac{1}{\e} F'(u_\e)^2] dx \leq  C \int_{ \{ u_\e \leq \eta  \} \cap A_1} \left[ \e F''(u_\e) |\nabla u_\e|^2 + \frac{1}{2 \e} F'(u_\e)^2 \right] dx.
  \]
Moreover, as there exists $c>0$ such that $|g'(s)| \leq c \eta $ for all $s \in [\eta,1/12]$, we deduce the existence of $C_0 > 0$ such that
\begin{eqnarray*}
  \int_{A_1 \cap \{ u_\e \le \eta \}} [ {\bf m}_\e(u_\e) + \frac{1}{\e} F'(u_\e)^2] dx &\leq&  C_0 \eta
  \int_{A_2 \cap \{ u_\e > \eta \}} [\e |\nabla u_\e|^2 ]dx + C_0 \e \int_{\R^d} ({\bf eul}_\e(u_\e))^2 dx \\
  &+&  C_0 \e^{1/2} \left[ \int_{\R^d}  {\bf m}_\e(u_\e) dx \right]^{1/2}.
 \end{eqnarray*}

We can now complete the proof of the lemma. Suppose by contradiction that there exists a sequence $\e\searrow 0^+$ such that $ B_r(x_0) \subseteq \{u_{\e} \leq \eta\}$. Using $A_1 = B_{r/2}(x_0)$ and
$A_2 = B_r(x_0)$, thus $ \{ u_\e >\eta \} \cap A_2 = \emptyset$, it follows that
\begin{eqnarray*}
 0 < \mu(B_{r/2}(x_0)) &\leq&  \liminf_{\e \to 0}  \int_{A_1 \cap  \{ u_\e \le \eta  \} }  {\bf m}_\e(u_\e) dx \\
 &\leq& C_0 \lim_{\e \to 0}  \Big(\e \int_{\R^d} ({\bf eul}_\e(u_\e))^2 dx + \e^{1/2} \left[ \int_{\R^d}  {\bf m}_\e(u_\e) dx \right]^{1/2}\Big) = 0,
\end{eqnarray*}
hence the contradiction.
 \end{proof}

\subsection{$\Gamma$ lim-inf: the one-dimensional case}  \label{sec:Gammaliminf_1}
In dimension one, as shown by Theorem~\ref{thm1d} below, the Willmore-type term does not contribute explicitly to the lower bound of the limit energy measure. However, it has a crucial role to guarantee the stability of the profile function $q'$. 

Let $Q = [a,b]$. For every $\e>0$, the one-dimensional energy $\E_\e$ of defined for all $u\in \wpqO{2}{2}([a,b])$ as:
\[
  \E_\e(u) =
  \begin{cases}\displaystyle
  \int_a^b \left(\frac{\e}{2}|u'|^2+\frac{1}{\e}F(u)\right) dx
  + \frac{\sigma_\e}{2\e}\int_a^b \left(\e u'' - \frac{1}{\e}F'(u)\right)^2
  dx & \text{ if } 0\le u \le \frac{1}{4} \text{ a.e.,}
  \\
  +\infty & \text{ otherwise.}
\end{cases}
\]
where $\sigma_\e$ is a positive function of $\e$ and $\sigma_\e \to \sigma_0 \geq 0$ as $\e \to 0$.

\begin{theo}\label{thm1d}
Assume that $\sigma_\e\to \sigma_0\geq 0$ and $\e^2/\sigma_\e\to 0$ as $\e\to 0^+$. Let $\{u_\e\}_\e\subset \wpqO{2}{2}([a,b])$ be a sequence such that $\E_\e(u_\e)\le C$ for all $\e>0$. Up to a subsequence, we can assume in addition that, as $\e\to 0$, $\mu_\e:=(\frac{\e }{2}|u'_\e|^2+\frac{1}{\e}F(u_\e)){\mathcal L}^1$ converges weakly-$*$ to a Radon measure $\mu$.
Then
 \begin{itemize}
\item  $(u_\e)$ converges a.e. to $u\equiv 0$.
\item The support of $\mu$ is a finite collection $\{x_i\}_{i\in I}$ of points and there exists a collection $\{m_i\}_{i\in I}$ of positive integers such that
$$ \mu = c_F \sum_{i \in I} m_i \delta_{x_i}.$$
\item  For $\e$ small enough, $u_\e\approx 1/4$ near every $x_i$.
\item  $ \liminf_{\varepsilon \to 0} \E_\e(u_\e) \geq  c_F \sum_{i \in I} m_i$.
\end{itemize}
\end{theo}
\begin{proof}
Since $\E_\e(u_\e)\le C<+\infty$, we have $\int_a^b F(u_\e)dx \le C\e$
and $\int_a^b \sqrt{2F(u_\e)}|u'_\e|dx  \le C$. This shows
that $w_\e(x) := \int_0^{u_\e(x)}\sqrt{2F(t)} dt$ is bounded in $BV(a,b)$. Remark that, for all $\e$, $u_\e\in {\mathrm C}^1$ by Sobolev injection.
Possibly passing to a subsequence, we may assume that $(w_\e)$ converges in $L^1(a,b)$
(or any $L^p(a,b)$, $p<+\infty$) and a.e.,
so that also $(u_\e)$ converges a.e. to some limit function $u$. Because $\int_a^b F(u_\e)dx\to 0$ we get that $F(u(x))=0$ a.e., so $u(x)\in \{0,\frac 14\}$ for a.e. $x$.\medskip


We first prove a uniform control of the pointwise discrepancy. We have, for all $x,y\in(a,b)$,
\[
  \frac{\e^2}{2} u_\e'^2(x) - F(u_\e(x)) = \frac{\e^2}{2} u_\e'^2(y) - F(u_\e(y))
  +\int_y^x \e^2 u_\e' u_\e''(t) - F'(u_\e(t))u'_\e(t) dt.
\]
Observe that
\begin{multline*}
  \int_y^x \e^2 u_\e' u_\e''(t) - F'(u_\e(t))u'_\e(t) dt
  \le \e\int_a^b |u'_\e(t)|\left|\e u_\e''(t)- \frac{1}{\e}F'(u_\e(t))\right|dt
  \\
  \le \e \|u'_\e\|_{L^2(a,b)} \left(\int_a^b\left(\e u''_\e - \frac{1}{\e}F'(u_\e)\right)^2dt\right)^{\frac{1}{2}}
    \le \frac{2\e}{\sqrt{\sigma_\e}} \E_\e(u_\e)
  \end{multline*}
  It follows that
  \[
    \left| \frac{\e^2}{2} u_\e'^2(x) - F(u_\e(x))\right| \le
    \e\left(\frac{\e}{2} u_\e'^2(y) + \frac{1}{\e}F(u_\e(y))\right)
    +
    \frac{2\e}{\sqrt{\sigma_\e}} \E_\e(u_\e)
  \]
  and averaging this over $y\in (a,b)$, we obtain
\begin{equation}\label{unif-discrep}
    \left| \e^2 u_\e'^2(x) - 2F(u_\e(x))\right| \le
    2\left(\frac{\e}{b-a}
    +
    \frac{2\e}{\sqrt{\sigma_\e}}\right) \E_\e(u_\e),
\end{equation}
so that if the energies are uniformly bounded, the quantity on the left-hand side goes uniformly to $0$.

\begin{figure}[!htbp]
\centering
	\includegraphics[width=.5\textwidth]{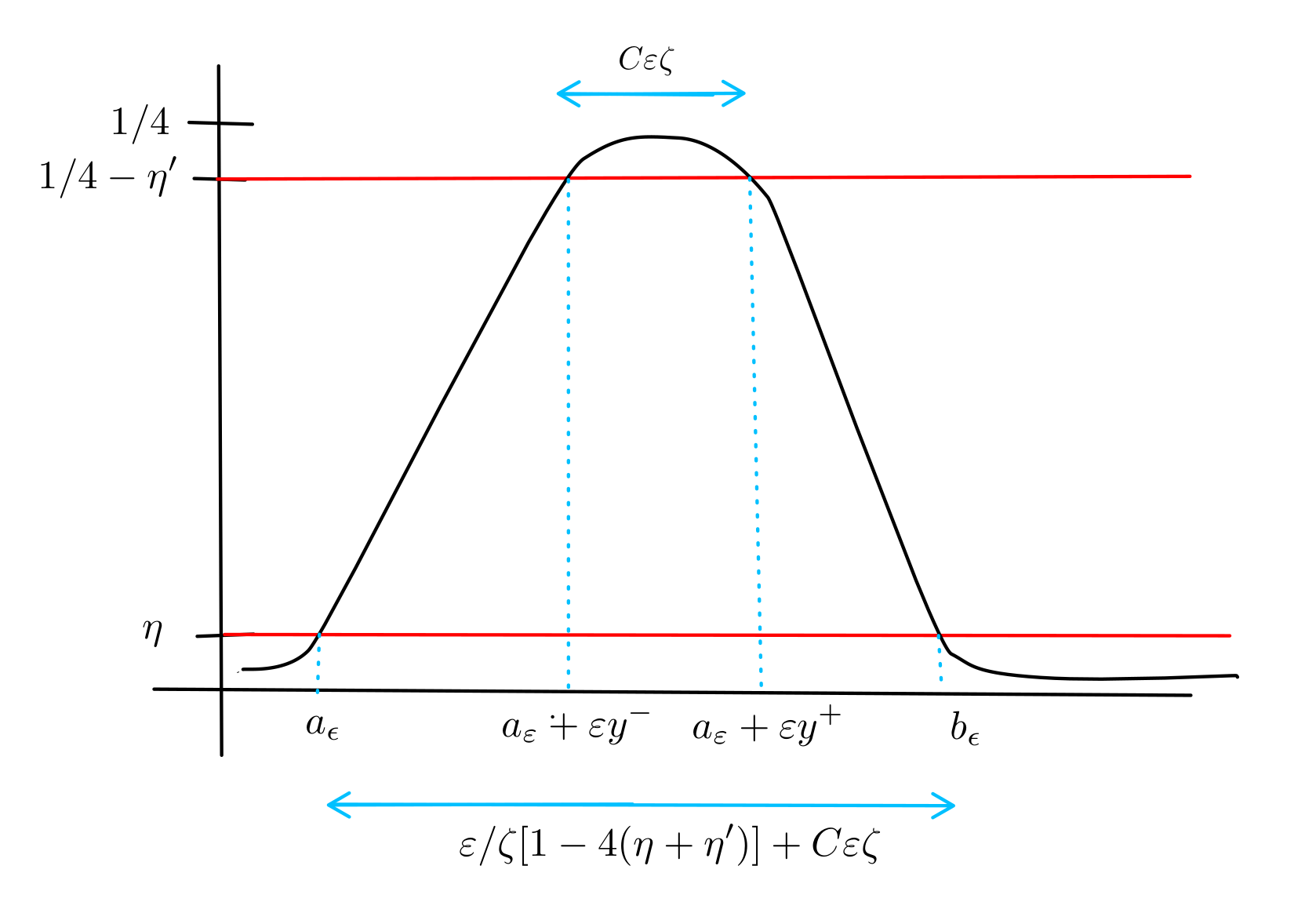}

\caption{Characteristic size of a diffuse interface}
\label{fig_discrepancy1}
\end{figure}

Let us now estimate the characteristic size of the transition zones between phases $s=0$ and $s=1/4$, as illustrated in Figure~\ref{fig_discrepancy1} .   Consider $\zeta>0$ small, and define $\eta, \eta'>0$ such that $\zeta=\sqrt{2F(\eta)}
= \sqrt{2F(1/4-\eta')}$. We assume that $\e>0$ is small enough, so that
  \[
     2\left(\frac{\e}{b-a}
    +
    \frac{2\e}{\sqrt{\sigma_\e}}\right) \E_\e(u_\e)\le \frac{\zeta^2}{4}.
\]
Observe that we can deduce from~\eqref{unif-discrep} that, for every $x\in(a,b)$, the two following inequalities hold:
\[
  \left| \e^2 u_\e'^2(x) - 2F(u_\e(x))\right| \le \frac{\zeta^2}{4}
 \qquad\text{ and }\qquad
  \left| \e |u_\e'(x)| - \sqrt{2F(u_\e(x))}\right| \le \frac{\zeta}{2}.
\]
(We use that if $|a^2-b^2|\le c^2$ for $a,b\ge 0$, then
either $a+b\ge c$ and $|a-b|\le c^2/c =c$, or $|a-b|\le a+b\le c$.)
Then in particular, if $\eta\le u_\e(x)\le 1/4-\eta'$, $\e^2u_\e'^2 (x)\ge 3\zeta^2/4$.
Let us now consider a connected component $(a_\e,b_\e)$ of $\{u_\e>\eta\}$.
By assumption, $u_\e(a)=u_\e(b)=0$, so $a_\e>a$ and $b_\e<b$. We consider, for $0\le y\le (b_\e-a_\e)/\e$, $v_\e(y) = u_\e(a_\e+ \e y)$. Since  $v_\e'(y)=\e u_\e'(a_\e+\e y)$, we have at $y=0$ that
$v'_\e(0)\ge \zeta/2 > 0$. Moreover, as long as $v_\e(y)\le 1/4-\eta'$, $v'_\e(y)\ge \zeta/2$, in particular $v_\e$ is  increasing. Denote $y^-$ the first point where it reaches the value $1/4-\eta'$,
then, by the Mean Value Theorem, $1/4-\eta' - \eta \ge \zeta y^-/2$.
Next, denote $y^+>y^-$ the first point where $v_\e$ reaches again the value
$1/4-\eta'$ ($y^+>y^-$ does exist because $v'_\e(y^-)>0$ and $b_\e<b$). Since $v_\e'(y^+)<0$ and $|v'_\e(y)|>0$ if $v_\e(y)\in[\eta,1/4-\eta']$, $v_\e$ has to decrease, until it reaches
the value $\eta$ again (which by definition is at the point $y=(b_\e-a_\e)/\e$).
Again, one has $1/4-\eta'-\eta \ge \zeta[(b_\e-a_\e)/\e - y^+]/2$.

 We now estimate
$y^+-y^-$. For all $y\in[y^-,y^+]$, one has $v_\e(y)\ge 1/4-\eta'$ and, if $1/4-\eta'>1/12$ (which is true if $\zeta$ is small enough), $0> F'(1/4-\eta')\geq F'(v_\e(y))$.
Observe that
\[
  \E_\e(u_\e)\ge \frac{\sigma_\e}{2\e^2} \int_{y^-}^{y^+} \left( v''_\e - F'(v_\e)\right)^2 dy
\]
and
\begin{multline*}
  \int_{y^-}^{y^+} \left( v''_\e - F'(v_\e)\right)^2 dy
  \ge \int_{y^-}^{y^+} -2F'(v_\e)v''_\e dy + (y^+-y^-)[F'(1/4-\eta')]^2
 \\ = 2F'(1/4-\eta')(v'_\e(y^-)-v'_\e(y^+)) +2\int_{y^-}^{y^+} v_\e'^2 F''(v_\e)
   +(y^+-y^-)[F'(1/4-\eta')]^2
 \end{multline*}
 Note that for $y\in[y^-,y^+]$, $v_\e'^2(y) \le 2F(v_\e(y)) + \zeta^2/4\le 5\zeta^2/4$
 and $|F''(v_\e(y))|\le 2$, therefore
 \[
   (y^+-y^-) ( [F'(1/4-\eta')]^2 - 5\zeta^2) \le \frac{2\e^2}{\sigma_\e}\E_\e(u_\e)
   + 2\sqrt{5}\zeta F'(1/4-\eta')
 \]
Since $F'(1/4)=-\frac 1 8$, there exist $\tilde\zeta,\delta>0$ such that  $[F'(1/4-\eta')]^2 - 5\zeta^2\ge \delta$ for all $\zeta\in (0,\tilde\zeta)$.

 Hence there is a constant $C$ depending only on $F$
 such that if $\e>0$ is small enough (depending on $\sup_\e\E(u_\e)$),
 $$y^+-y^-\le C\zeta.$$ 
 
 All in all, we end up finding that
 \begin{equation}\label{control-b-a}
   b_\e-a_\e \le \frac{\e}{\zeta}\left(1 - 4(\eta+\eta')\right) + C\e\zeta.
 \end{equation}

We now estimate the energy inside the interval $[a_\e,b_\e]$. We note
 that on $[0,y^-]$,
 \[
   \int_0^{y^-} \big(\frac{v_\e'^2}{2} + F(v_\e)\big) dy
   =    \int_0^{y^-} \Big(\sqrt{2F(v_\e)}|v'_\e| + \frac{1}{2}\left(|v'_\e| - \sqrt{2F(v_\e)}\right)^2\Big) dy
 \]
 so that
 \[
\left|   \int_0^{y^-} \big(\frac{v_\e'^2}{2} + F(v_\e)\big) dy  -\int_{\eta}^{1/4-\eta'} \sqrt{2F(t)}dt\right|\le  y^-\frac{\zeta^2}{8}
\le \frac{\zeta}{8}\left(\frac{1}{2}-2(\eta'+\eta)\right).
\]
The same estimate holds on $[y^+,(b_\e-a_\e)/\e]$. In $[y^-,y^+]$,
one has $v_\e'^2/2 + F(v_\e)\le 2F(v_\e)+\zeta^2/4 \le 5\zeta^2/4$ so
that
\[
  \int_{y^-}^{y^+} \big(\frac{v_\e'^2}{2} + F(v_\e)\big) dy
  \le \frac{5C}{4}\zeta^3.
\]
We deduce that (still with $c_F = 2\int_0^{1/4}\sqrt{2F(t)}dt$)
\begin{equation}\label{control-energy-a-b}
  \int_{a_\e}^{b_\e} \big(\frac{\e}{2}u_\e'^2 + \frac{1}{\e}F(u_\e)\big) dx
  = c_F + o_\zeta(1).
\end{equation}

As $ \E_{\e}(u_\e) \leq C$, we deduce that, for $\eta$ given, there exists a finite number $N$ of connected components, the same $N$ for all $\varepsilon$ sufficiently small, such that
$ \{ u_{\varepsilon} > \eta \} = \bigcup_{i=1}^N   (a^{i}_{\varepsilon},b^{i}_{\varepsilon})$ with
$$\cdots<a_\e^{i} < b_\e^{i} < a_\e^{i+1} <b_\e^{i+1} <\cdots$$
In addition, \eqref{control-b-a} and \eqref{control-energy-a-b} imply that
$ \mu_{\e}([a^{i}_{\varepsilon},b^{i}_{\varepsilon}]) = c_F +  o_{\zeta}(1)$ and $\lim_{\e \to 0} a^{i}_{\varepsilon} = \lim_{\e \to 0} b^{i}_{\varepsilon} =: x_i \in [a,b]$. Since, by Lemma \ref{lemme_xi}, for all $x \in  \text{supp}(\mu) $ there exists $x_\e \in \{ u_\e > \eta\} \to x$, we deduce that the support of $\mu$ coincides with $\{x_i\}_{i=1,\dots,N}$.
The finite number of connected components of $\{ u_\e > \eta \}$, together with the fact that $b_\e-a_\e \le \frac{\e}{\zeta}\left(1 - 4(\eta+\eta')\right) + C\e\zeta$ and $y^{-}_i<y^{+}_i$, imply that $(u_\e)$ converges a.e. to $u\equiv 0$ and, for $\e$ small enough, $u_\e\approx 1/4$ near every $x_i$.

Moreover, as
$u_{\e} \in \wpqO{2}{2}([a,b])$, an extension argument and inequality \eqref{equ_control_bulk} with $A_1 = [a,b]$ and $A_2 = [a-1,b+1]$ show that
$$  \lim_{\e \to 0} \int_{\{ u_{\e} \leq \eta\}} (\frac{\e}{2}u_\e'^2 + \frac{1}{\e}F(u_\e))dx = O(\eta) = o_{\zeta}(1).$$

By combining with \eqref{control-energy-a-b} and passing to the limit $\zeta \to 0$, we obtain that
$$ \lim_{\e} \mu_{\e}([a,b]) =  N c_F  \qquad\text{ and }\qquad  \mu = \sum_{i=1}^{N} c_F \delta_{x_i}.$$
 Two limit points $x_i$, $x_j$ need not be distinct, so we consider the collection $\{x_i, \,i\in I\}$ of reindexed distinct limit points, together with their positive integral multiplicities $\{m_i, i\in I\}$, and we finally conclude that
 $$ \liminf_{\e \to 0} \E_\e(u_\e) \geq   \liminf_{\e \to 0} \mu_{\e}([a,b])  =  c_F \sum_{i \in I} m_i,$$
and
  $$\mu =  c_F  \sum_{i \in I}  m_i \delta_{x_i}.$$

\end{proof}

\subsection{$\Gamma$-liminf: the radially symmetric case in dimension $\mathbf{d \geq 2}$} \label{sec:Gammaliminf_d}

\par We now consider the problem in $\R^d$, $d\geq 2$, looking only at radially symmetric functions.
The analysis is very similar to the 1D case, it shows in addition that, in dimension $d\ge 2$, the Willmore approximation term not only stabilizes the
profile $q'$ but also contributes explicitly in the limit to the energy of the limit support.

For $u\in \wpqO{2}{2}(Q)$, we consider
\[
  \E_\e(u) =
  \begin{cases}\displaystyle
  \int_Q {\bf m}_\e(u) dx
  + \frac{\sigma_\e}{2 \e}  \int_Q  \left[{\bf eul}_\e(u)\right]^{2}
  dx & \text{ if } 0\le u \le \frac{1}{4} \text{ a.e.,}
  \\
  +\infty & \text{otherwise,}
\end{cases}
\]
 where $\sigma_\e>0$, $\sigma_\e\to \sigma_0 \geq 0$ as $\e \to 0$, and 
 \begin{eqnarray*}
 {\bf m}_\e(u) &=&\ds   \frac{\e}{2}|\nabla u|^2+\frac{1}{\e}F(u)\\
 {\bf eul}_\e(u) &=& \e \Delta u - \frac{1}{\e}F'(u).
 \end{eqnarray*}

\begin{theo}\label{thmdd}
Let $R_B>R_A>0$ and $Q = B_{R_B} \setminus B_{R_A}$. Let $(u_\e)_\e$ be a sequence of radially symmetric functions in $\wpqO{2}{2}(Q)$ such that $\E_\e(u_\e)\le C<+\infty$. Assume that $\lim_{\e \to 0^+} \sigma_\e=\sigma_0\geq 0$ and $\lim_{\e \to 0^+}\e^2/\sigma_\e\to 0$. Possibly taking a subsequence, we assume that there exists a Radon measure $\mu$
such that, as $\e\to 0^+$, $\mu_\e:={\bf m}_\e(u_\e){\mathcal L}^d$  converges weakly-$*$ to $\mu$. 
 Then, possibly passing to a subsequence,
 \begin{itemize}
\item  $(u_\e)$ converges a.e. to $u\equiv 0$.
\item The support of $\mu$ is concentrated on a finite collection $\{\partial B_{r_i},\, i\in I\}$ of $(d-1)$-spheres centered at $0$ and there exists a collection of positive integer multiplicities $\{m_i,\,i\in I\}$ such that
$$ \mu =    c_F  \sum_{i\in I}  m_i \haushd\res{\partial B_{r_i}}.$$
\item For $\e$ small enough, $u_\e\approx 1/4$ near $\partial B_{r_i}$ for all $i\in I$.
\item $\displaystyle \liminf_{\e \to 0}   \E_\e(u_{\e}) \geq  c_F \sum_{i\in I} m_i \left( P(B_{r_i}) + \sigma_0 \mathcal{W}(B_{r_i}) \right)$.
\end{itemize}
\end{theo}

\begin{proof}
We first prove that, possibly passing to a subsequence, for a.e. $x$, $u_\e(x)\to u(x)\in \{0,\frac 14\}$.  Indeed, as in the one-dimensional case,
 we have $\int_{Q} F(u_\e)dx \le C\e$ and $\int_{Q} \sqrt{2F(u_\e)}|\nabla u_\e|dx  \le C$. This shows
that $w_\e(x) := \ds\int_0^{u_\e(x)}\sqrt{2F(t)} dt$ is bounded in $\BV(Q)$.
Possibly passing to a subsequence, we may assume that $(w_\e)$ converges in $L^1(Q)$
(or any $L^p(B_R)$, $p<+\infty$) and almost~everywhere,
so that also $(u_\e)$ converges a.e. to some limit function $u$. Because $\int_{Q} F(u_\e)dx\to 0$ we get that $F(u)=0$ a.e., therefore $u(x)\in \{0,\frac 14\}$ a.e. \medskip


Using the radial symmetry of $u_{\e}$, let us introduce $\overline{u}_\e : [R_A,R_B] \to \R$ such that,
for all $x \in Q$,  $u_\e(x) = \overline{u}_\e(|x|)$.
Notice that
\begin{equation}\label{controlemasse}
\int_{Q} {\bf m}_\e(u_\e) dx = |S^{d-1}|  \int_{R_A}^{R_B} \left[ \e \frac{|\overline{u}_\e'(t)|^2}{2} + \frac{F(\overline{u}_\e)}\e \right] t^{d-1}  dt \leq C,
\end{equation}
and
$$ \frac{\sigma_\e}{\e} \int_{Q} {\bf eul}_\e(u_\e)^2 dx = |S^{d-1}|  \frac{\sigma_\e}{\e} \int_{R_A}^{R_B} \left[ -\e \overline{u}_\e''(t) + \frac{F'(\overline{u}_\e'(t))}\e  - \frac{\e (d-1)}{t} \overline{u}_\e'(t) \right]^2  t^{d-1} dt \leq C,$$
where $S^{d-1}= \{ x \in \R^d,\, \|x \| = 1 \}$ is the unit sphere.

The first additional difficulty compared to the one-dimensional case comes from the presence of
the  term $w_\e = -\e \frac{(d-1)}{t} \overline{u}_\e'$
which does not allow us to have a direct control with $\E_{\e}(u_\e)$ of the term $e_{\e} = -\e \overline{u}_\e'' + F'(\overline{u}_\e')/\e$. However, for all $(x,y)\in [R_A,R_B]^2$, we have

\[
  \frac{\e^2}{2} \overline{u}_\e'^2(x) - F(\overline{u}_\e(x)) = \frac{\e^2}{2} \overline{u}_\e'^2(y) - F(\overline{u}_\e(y))
  +\int_y^x [\e^2 \overline{u}_\e' \overline{u}_\e''(t) - F'(\overline{u}_\e(t))\overline{u}'_\e(t)] dt,
\]
therefore

\begin{multline*}
\left|  \frac{\e^2}{2} \overline{u}_\e'^2(x) - F(\overline{u}_\e(x))\right|
  \le\e \frac{1}{R_A^{d-1}}
  \left( \frac{\e}{2} u_\e'^2(y) + \frac{1}{\e}F(u_\e(y))\right)|y|^{d-1}
  \\+ \frac{1}{R_A^{d-1}}\int_{[x,y]} \e |u_\e'(t)| \left|\e u_\e''(t) -\tfrac{1}{\e} F'(u_\e(t))\right||t|^{d-1} dt.
\end{multline*}

Then  we calculate that
\begin{multline*}
  \int_{[x,y]} \e |\overline{u}_\e'(t)| \left|\e \overline{u}_\e''(t) -\tfrac{1}{\e} F'(\overline{u}_\e(t))\right||t|^{d-1} dt
  \\
  \le   \int_{[x,y]} \e |\overline{u}_\e'(t)| \left|\e (\overline{u}_\e''(t) +\tfrac{d-1}{t}\overline{u}_\e'(t)) -\tfrac{1}{\e} F'(\overline{u}_\e(t))\right||t|^{d-1} dt
  + \int_{[x,y]} \e |u_\e'(t)| \e\tfrac{d-1}{t}|\overline{u}_\e'(t)||t|^{d-1}dt
  \\
  \le \e^2\left(\tfrac{1}{2\sqrt{\sigma_\e}}+\tfrac{d-1}{x}\right)\int_{[x,y]}|\overline{u}_\e'(t)|^2 |t|^{d-1}dt
  + \frac{\sqrt{\sigma_\e}}{2}\int_{[x,y]} \left|\e (\overline{u}_\e''(t) +\tfrac{d-1}{t}\overline{u}_\e'(t))
    -\tfrac{1}{\e} F'(u_\e(t))\right|^2|t|^{d-1} dt
  \\
  \le C\left(\frac{\e}{R_A} + \frac{\e}{\sqrt{\sigma_\e}}\right)\E_\e(u_\e).
\end{multline*}

Averaging over $y\in (R_A,R_B)$, we obtain as in 1D a uniform  control of the pointwise discrepancy
\[
  \left|\e^2 \overline{u}_\e'^2(x) - 2F(\overline{u}_\e(x))\right|
  \le  C\left( \frac{\e}{R_A^{d}} + \frac{1}{R_A^{d-1}} \frac{\e}{\sqrt{\sigma_\e}}  \right) \E_\e(u_\e).
\]
We assumed that $R_A>0$, but the case $R_A=0$ can be handled similarly by integrating $y$ on $[R_B/2,R_B]$ to obtain a control that is only locally uniform around $0$.

As in the one-dimensional case, we can now estimate the length of a diffuse interface.  Let us consider $\zeta>0$ small, and define $\eta, \eta'>0$ such that $\zeta=\sqrt{2F(\eta)} = \sqrt{2F(1/4-\eta')}$. We also assume that $\e>0$ is small enough, so that
  \[
    \left|\e^2 \overline{u}_\e'^2(x) - 2F(\overline{u}_\e(x))\right|\le \frac{\zeta^2}{4}.
\]
Let us consider a connected component $(a_\e,b_\e) \subset [R_A,R_B]$ of $\{\overline{u}_\e>\eta\}$.
As in the 1D case,  as $\overline{u}_\e(R_A)=\overline{u}_\e(R_B)=0$,
we have $a_\e>R_A$ and $b_\e<R_B$. We consider
$v_\e(y) = \overline{u}_\e(a_\e+ \e y)$, for $0\le y\le (b_\e-a_\e)/\e$.  By construction,  $v_\e$ is  increasing and  reaches the value $1/4-\eta'$ at $y^-$.  Next, as in 1D, denote $y^+>y^-$ the first point where $v_\e$ reaches again the value
$1/4-\eta'$. Then $v_\e$ keeps decreasing, until it reaches
the value $\eta$ again (which by definition happens at $y=(b_\e-a_\e)/\e$).
So, using the same arguments as in the 1D case,
as
$$
 \int_{y^-}^{y^+} \left( v''_\e - F'(v_\e)\right)^2 dy \leq  C\left( \frac{\e}{R_A^{d}} + \frac{1}{R_A^{d-1}} \frac{\e}{\sqrt{\sigma_\e}}  \right) \E_\e(u_\e).
$$
we can show that $y^{+} - y^{-} \leq C \zeta$ and $b_\e-a_\e \le  \e C_{\zeta}.$

We now estimate the mass energy inside the interval $[a_\e,b_\e]$. We note  that on $[0,y^-]$, there exists $c_{\e,\zeta} \in [a_\e^{d-1},(a_\e + \e y^{-})^{d-1}]$ such that
$$
 \int_0^{y^-} \big(\frac{v_\e'^2}{2} + F(v_\e)\big) |a_{\e} + \e  y |^{d-1} dy = c_{\e,\zeta} \int_0^{y^-} \big(\frac{v_\e'^2}{2} + F(v_\e)\big)  dy.
$$
We can deduce as in the 1D case that
$$
\int_{B_{b_\e} \setminus B_{a_\e}} {\bf m}_\e(u_\e) dx =
|S^{s-1}|\int_{a_\e}^{b_\e} [\frac{\e}{2}\overline{u}_\e'^2 + \frac{1}{\e}F(\overline{u}_\e)] x^{d-1} dx = |S^{d-1}| c_{\e,\zeta} c_F +  o_\zeta(1).$$

As $ \E_{\e}(u_\e) \leq C$, it follows that, for $\zeta$ given, there exists a finite number $N$ of connected components, the same for all $\varepsilon$ sufficiently small, such that
$ \{ \overline{u}_{\varepsilon} > \eta \} = \bigcup_{i=1:N}   [a^{i}_{\varepsilon},b^{i}_{\varepsilon}]$ with
$$ \mu_{\e}(B_{b^{i}_{\varepsilon}} \setminus B_{a^{i}_{\varepsilon}}) = |S^{s-1}| c_{\e,\zeta}^i c_F +  o_\zeta(1) $$
Moreover, defining
$\lim_{\e \to 0} a^{i}_{\varepsilon} = \lim_{\e \to 0} b^{i}_{\varepsilon} =: r_i \in [R_a,R_b]$,
Lemma \ref{lemme_xi} shows that the support of $\mu$ coincides with the union of the spheres $\partial B_{r_i}$ centered at $0$, with radius $r_i$, and as

Moreover, as
$u_{\e} \in \wpqO{2}{2}(Q)$, an extension argument and \eqref{equ_control_bulk} with $A_1 = B_{R_B}\setminus B_{R_A}$ and $A_2 = B_{R_B+1}\setminus B_{0.9 R_A}$ show that
$$  \lim_{\e \to 0} \int_{\{ u_{\e} \leq \eta\}} \frac{\e}{2}|\nabla u_\e|^2 + \frac{1}{\e}F(u_\e) dx = O(\eta) = o_{\zeta}(1).$$

Passing to the limit $\zeta \to 0$, it follows that
$$ \lim_{\e \to 0} \mu_{\e}(Q) =  c_F |S^{d-1}| \sum_{i=1}^{N} r_i^{d-1}\qquad\text{and}\qquad\mu =     c_F  \sum_{i=1}^{N} \haushd\res{\partial B_{r_i}}.$$
 As two radii $r_i, r_j$ are not necessarily distinct, these decompositions can be rewritten using a collection of distincts radii $r_i$, $i\in I$ (using the same notation for simplicity) associated with positive integral multiplicities $m_i \in \N^{*}$ such that
 $$ \mu =     c_F  \sum_{i\in I} m_i \haushd\res{\partial B_{r_i}},$$
 and
 $$ \liminf_{\e \to 0}  \int_Q {\bf m}_\e(u) dx \geq   c_F |S^{d-1}| \sum_{i \in I} m_i r_i^{d-1} = c_F \sum_{i \in I} m_i P(B_{r_i}).$$

Let us now evaluate the contribution of the Willmore term around a given radius $r_i$.
Up to a subsequence, $ \mu_{\e} =  {\bf m}_\e(u_\e)dx  \wto \sum_i m_i c_F \haushd\res{\partial B_{r_i}}$,
where, for every $i$, the multiplicity $m_i$ corresponds to the number of connected components
$(a^{i}_\e,b^{i}_\e)$ that concentrate at $r_i$. Assuming we have first extracted
a subsequence $(\e_k)$ such that $\displaystyle\liminf_{\e} \int_Q {\bf eul}(u_\e)dx = \lim_k \int_Q {\bf eul}(u_{\e_k})dx$
 we can extract another subsequence
such that ${\bf m}_{\e_k}(u_{\e_k})(x)\to 0$ for a.e.~$x$. For simplicity, we drop the index $k$ in the sequel. \\

As illustrated on Figure \ref{fig_willmore1}, we select a particular limit radius $r_i$ and $\delta>0$ such that
${\bf m}_{\e}(u_\e)(r_i\pm \delta) \to 0$, and $[r_i-\delta,r_i+\delta]$ does
not contain any other radius $r_j$. Let $\alpha:=r_i-\delta$, $\beta:=r_i+\delta$.

\begin{figure}[!htbp]
\centering
	\includegraphics[width=.5\textwidth]{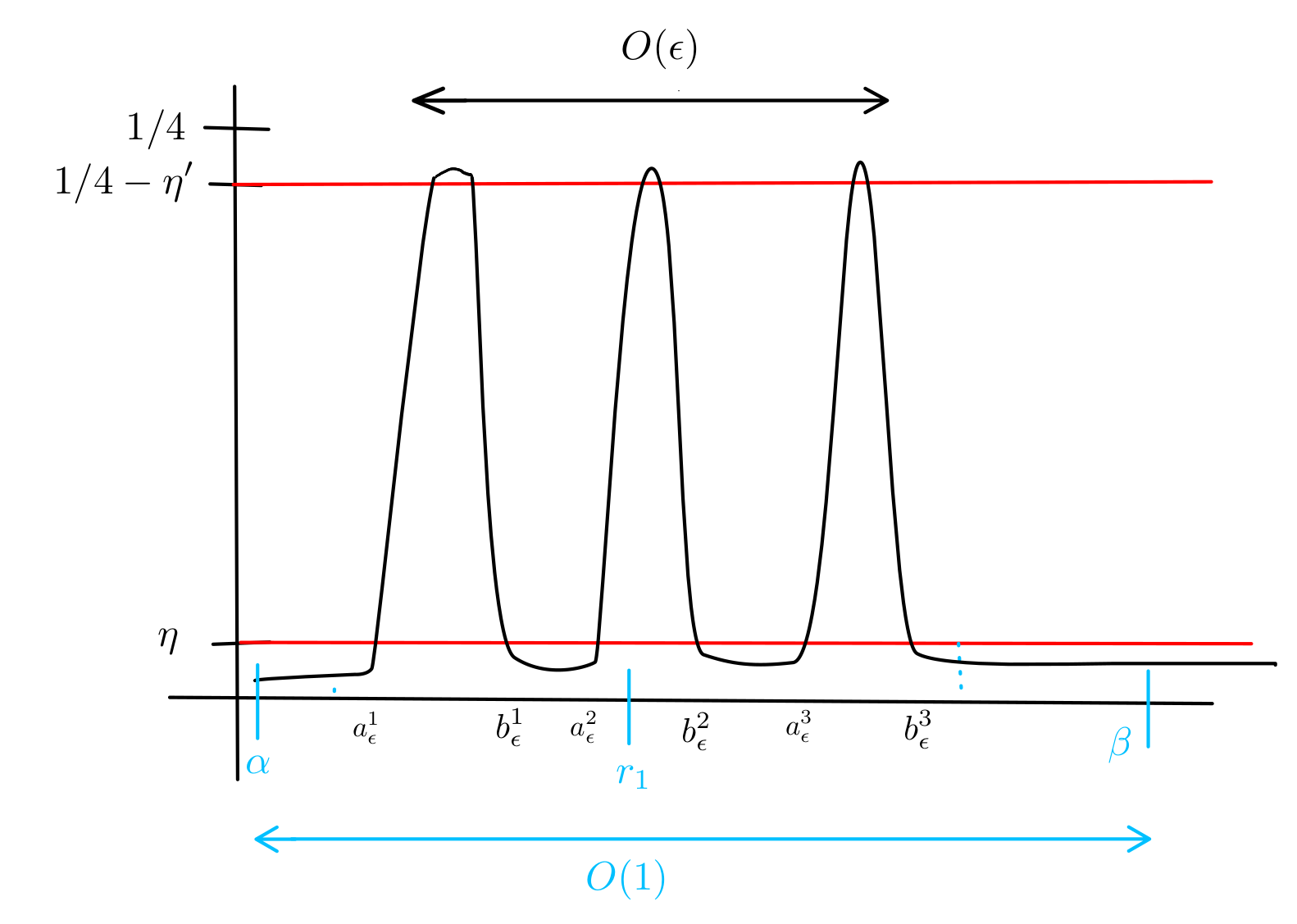}

\caption{Estimating the Willmore term on $[r_i-\delta,r_i+\delta]$}
\label{fig_willmore1}
\end{figure}

We focus on the term
\[
 \frac{\sigma_\e}{2\e} \int_{B_{\beta} \setminus B_{\alpha}  } {\bf eul }_\e^2(u_{\e}) dx = \frac{\sigma_\e}{2\e}|S^{d-1}|\int_{\alpha}^{\beta} \left(
    \e u_\e''(t) + \e\frac{d-1}{t}u'_\e(t) - \frac{1}{\e}F'(u_\e(t))\right)^2 t^{d-1}dt
\]
First, we write that
\begin{align*}
  \frac{\sigma_\e}{2\e}\int_{\alpha}^{\beta}& \left(
    \e u_\e''(t) + \e\frac{d-1}{t}u'_\e(t) - \frac{1}{\e}F'(u_\e(t))\right)^2 t^{d-1}dt
  \\ & = \frac{\sigma_\e}{2\e}\int_{\alpha}^{\beta}\frac{d-1}{t} \left(
        \e\sqrt{\frac{d-1}{t}}u'_\e(t)+
    \sqrt{\frac{t}{d-1}}\Big(\e u_\e''(t)  - \frac{1}{\e}F'(u_\e(t))\Big)\right)^2 t^{d-1}dt
  \\ & \ge (d-1)\frac{\sigma_\e \alpha^{d-2}}{2\e}\int_{\alpha}^{\beta} \left(
        \e\sqrt{\frac{d-1}{t}}u'_\e(t)+
    \sqrt{\frac{t}{d-1}}\Big(\e u_\e''(t)  - \frac{1}{\e}F'(u_\e(t))\Big)\right)^2 dt
  \\ & \ge (d-1)\frac{\sigma_\e \alpha^{d-2}}{2\e}\int_{\alpha}^{\beta}
       \Big[\e^2\frac{d-1}{t} u'_\e(t)^2 dt + 2         \e u'_\e(t)
       \Big(\e u_\e''(t)  - \frac{1}{\e}F'(u_\e(t))\Big)\Big] dt
  \\ & \ge (d-1)^2 \sigma_\e \frac{\alpha^{d-2}}{\beta}\int_{\alpha}^{\beta}
       \frac{\e}{2} u'_\e(t)^2 dt + (d-1)\frac{\sigma_\e \alpha^{d-2}}{\e}\int_{\alpha}^{\beta} u'_\e(t)
       \Big(\e^2 u_\e''(t)  - F'(u_\e(t))\Big) dt.
\end{align*}
The last term satisfies
\begin{multline*}
  (d-1)\frac{\sigma_\e \alpha^{d-2}}{\e} \Big(\e^2 \frac{u'_\e(\beta)^2-u'_\e(\alpha)^2}{2} - F(u_\e(\beta))+F(u_\e(\alpha))\Big)
  \\=
  (d-1)\sigma_\e \alpha^{d-2} \Big(\frac{\e }{2}u'_\e(\beta)^2- \frac{1}{\e}  F(u_\e(\beta)) -\frac{\e}{2}u'_\e(\alpha)^2+ \frac{1}{\e}F(u_\e(\alpha))\Big)
  \stackrel{\e\to 0}{\to} 0
  \end{multline*}
  by our choice of $\alpha,\beta$.

As for the first term, given $\zeta>0$ as before and $\e$ small enough, we
  write:
  \[
     \sigma_\e (d-1)^2  \frac{\alpha^{d-2}}{\beta}\int_{\alpha}^{\beta}
    \frac{\e}{2} u'_\e(t)^2 dt
    \ge
    \sigma_\e (d-1)^2  \frac{\alpha^{d-2}}{\beta}
    \sum_{l=1}^{m_i} \frac{\e}{2} \int_{a_\e^l}^{b_\e^l}u'_\e(t)^2 dt
  \]
  where $m_i$ is the integral multiplicity of the limit concentrated mass measure at $r_i$ and
  $(a_\e^\ell,b_\e^\ell)$, $\ell=1,\dots,m_i$, are the connected components of $\{u_\e>\eta\}$
  that collapse at $r_i$. Using the limit mass and the control of the discrepancy, we have 
  \[
    \liminf_\e \e \int_{a_\e^\ell}^{b_\e^\ell}  u'_\e(t)^2dt \ge (c_F - c(\eta)),
  \]
  with $c(\eta)\to 0$ as $\eta\to 0$. Therefore, letting  $\delta \to 0$ and $\zeta \to 0$,  we can deduce that

\begin{eqnarray*}
 \liminf_{\e\to 0 }\frac{\sigma_\e}{2\e} \int_{B_{\beta} \setminus B_{\alpha}  } {\bf eul }_\e^2(u_{\e}) dx
    &\geq & \frac{1}{2} \sigma_0 m_i  (c_F - c(\zeta))  | S^{d-1}| (d-1)^2  \frac{(r_i - \delta)^{d-2}}{(r_i + \delta)} \\
    &\geq&  \frac{1}{2}  \sigma_0  c_F     m_i |S^{d-1}|  (d-1)^2 r_i^{d-3}, \\
    & \geq&  \sigma_0 m_i  \mathcal{W}(B_{r_i}).
\end{eqnarray*}
which shows that

\begin{eqnarray*}
 \liminf_{\e \to 0} \E_\e(u_\e) &\geq&  c_F |S^{d-1}| \sum_{i \in I} m_i \left( r_i^{d-1} +   \frac{\sigma_0}{2} \frac{(d-1)^2}{r_i^{3-d}} \right) \\
&\geq&  c_F \sum_{i \in I} m_i \left[ P(B_{r_i}) + \sigma_0 \mathcal{W}(B_{r_i}) \right].
\end{eqnarray*}
\end{proof}

The previous theorem can be extended to the domain $Q = B_R$. However, the statement is slightly less precise because various concentration phenomena can occur at $r=0$.
\begin{theo}\label{thmdd2}
Let $R>0$ and $Q = B_{R}$. Let $(u_\e)_\e$ be a sequence of radially symmetric functions in $\wpqO{2}{2}(Q)$ such that $\E_\e(u_\e)\le C<+\infty$. Assume that $\lim_{\e \to 0^+} \sigma_\e=\sigma_0\geq 0$ and $\lim_{\e \to 0^+}\e^2/\sigma_\e\to 0$. Possibly taking a subsequence, we assume that there exists a Radon measure $\mu$
such that, as $\e\to 0^+$, $\mu_\e:={\bf m}_\e(u_\e){\mathcal L}^d$  converges weakly-$*$ to $\mu$. 
 Then, possibly passing to a subsequence,
 \begin{itemize}
\item  $(u_\e)$ converges a.e. to $u\equiv 0$.
\item The support of $\mu$ is concentrated on a (at most) countable collection $\{\partial B_{r_i},\, i\in I\}$ of $(d-1)$-spheres centered at $0$ with radius $r_i>0$, and it possibly contains $0$. Moreover,
there exists a collection of multiplicities $\{m_i \in \mathcal{N}^*,\,i\in I\}$ and a real $\alpha \geq 0$ such that
$$ \mu =    c_F  \sum_{i\in I}  m_i \haushd\res{\partial B_{r_i}} + \alpha \delta_0$$
\item For $\e$ small enough, $u_\e\approx 1/4$ near $\partial B_{r_i}$ for all $i\in I$.
\item $\displaystyle \liminf_{\e \to 0}   \E_\e(u_{\e}) \geq  c_F \sum_{i\in I} m_i \left( P(B_{r_i}) + \sigma_0 \mathcal{W}(B_{r_i}) \right)$.
\end{itemize}

\end{theo}

\begin{remk}
The main difference with Theorem~\ref{thmdd} lies in the lack of controllability of the measure $\mu$  at $r=0$, which prevents us from ensuring that $\alpha$
is an integer and thus from taking its contribution to the limit energy of $\E_{\varepsilon}$ into account. Indeed, even if we can establish a control of the discrepancy around $r=0$, this control is no longer uniform and therefore no longer allows us to apply the previous reasoning.
\end{remk}

\begin{proof} The result follows from applying Theorem~\ref{thmdd} to the set $Q_r = B_R \setminus B_r$, for all $r$ in a sequence of positive infinitesimals, combined with a diagonal argument.
The novelty comes from the fact that  $u_{\varepsilon}$ may not vanish on the boundary $\partial B_r$.
In particular, if $\overline{u}_{\varepsilon}(r) > \delta$ (using the same notation and arguments as in the previous proof),
the first connected component has the form $(a_{\varepsilon},b_{\varepsilon}) = (r,b_{\varepsilon})$ and the contribution
$ \mu_{\varepsilon}( B_{b_{\varepsilon}} \setminus B_r)$
cannot be evaluated as $c_F + o_{\zeta}(1)$ like for the other components. Instead, combining with Theorem~\ref{thmdd},  we have that
$$ \mu\res{Q_r} =    c_F  \sum_{i\in I_{r}}  m_i \haushd\res{\partial B_{r_i}} + \alpha_r \haushd\res{\partial B_{r}},$$
where $\alpha_r \in \R^{+}$ .  Notice also that
\begin{eqnarray*}
 \liminf_{\e \to 0} \E^{r}_\e(u_\e) &=& \int_{Q_r} {\bf m}_\e(u) dx
  + \frac{\sigma_\e}{2 \e}  \int_{Q_r}  \left[{\bf eul}_\e(u)\right]^{2}
  dx \\
 &\geq& c_F \sum_{i \in I_{r}} m_i \left[ P(B_{r_i}) + \sigma_0 \mathcal{W}(B_{r_i}) \right].
\end{eqnarray*}
The final result follows from a diagonal argument with $r\to 0$.
 \end{proof}

\subsection{$\Gamma$-limsup} \label{sec:Gammalimsup}

Consider a bounded subset $E \subset Q$ with  $C^{2}$ boundary. Following the same arguments
as in~\cite{bellettini_1996}, one may construct a sequence $(u_{\e})_\e$
such that
$$ \limsup_{\e \to 0} \E_\e(u_\e)  \leq  c_F \left( \mathcal{P}(E) + \sigma_0 \mathcal{W}(E) \right).$$
A suitable sequence is for instance  
$$ u_\e(x) = \gamma_\e \left( \frac{\operatorname{dist}(x,E)}{\e}  \right),$$
where $\gamma_\e$ is the $C^{1,1}$ truncated profile   defined by
$$\gamma_\e(s) =  \begin{cases}
                 -q'(s)  &\text{ if } |s| \leq x_\e \\
                 p_\e(s)& \text{ if }   x_\e < |s| \leq 2 x_\e \\
                 0 & \text{ if }  |s| > 2x_\e
                \end{cases}
$$
where $x_\e = |\log(\e)| $ and  $p_\e$ is a degree 3 polynomial satisfying
$$ p_\e(x_\e) = -q'(x_\e), \;\;  p'_\e(x_\e) = -q''(x_\e), \;\; p_\e(x_\e) = -q'(x_\e) \;\text{ and }\;   p'_\e(2x_\e) = 0.$$
Moreover, it is not difficult to see  that $u_{\e} \to 0$ and that
$$ \mu_\e = \frac{\epsilon^2}{2} |\nabla u_\e|^2 + \frac{1}{\e} F(u_\e) \rightharpoonup  c_F\haushd\res{\partial E}.$$
In the case that concerns us in this work, we are trying to approximate measures of the form
$$
\mu = \sum_i m_i c_F \haushd\res{\partial B_{r_i}}.$$
By a density argument, it is possible to assume that $m_i=1$  and that the distance between each radius is greater than $\delta$. In this case, for $\e$ sufficiently small ( $\e |\log(\e)| < \delta/2$), it is possible to show, using arguments as in~\cite{bellettini_1996}, that the sequence
$$ u_\e(x) = \gamma_\e \left( \frac{\operatorname{dist}(x,\sum_{i \in I} \partial B_{r_i} )}{\e}  \right),$$
satisfies
\begin{itemize}
 \item $u_\e \in  \wpqO{2}{2}(B_{R_B} \setminus B_{R_A})$ with $R_A = \min_{i \in I} r_i - \delta$ and $R_B = \max_{i \in I} r_i +\delta$.
 \item $u_\e \to 0$ in $L^1(Q)$
 \item $\mu_\e(u_\e) = \frac{\epsilon^2}{2} |\nabla u_\e|^2 + \frac{1}{\e} F(u_\e) \rightharpoonup  c_F \sum_{i \in I} \haushd\res{\partial B_{r_i}}$.
 \item $ \limsup_{\e \to 0} \E_\e(u_\e) \leq  c_F  \sum_{i\in I} \left( P(B_{r_i}) +  \sigma_0 \mathcal{W}(B_{r_i})\right)$.
\end{itemize}
from which the $\Gamma$-limsup property needed for Theorem~\ref{theo-Gamma-conv} ensues.

\section{Numerical experiments}\label{sec:numeric}

We propose in this section a simple numerical scheme to approximate the $L^2$-gradient flow of $\E_\e$, whose definition for $u \in \wpq{2}{2}(Q)$ satisfying $0\leq u \leq 1/4$ is recalled:
\[
  \E_\e(u,Q) = \int_Q \left(\frac{\e}2|\nabla u|^2+\frac{1}{\e}F(u)\right) dx
  + \frac{\sigma_\e}{2\e}  \int_Q  \left( \e \Delta u - \frac{1}{\e}F'(u)\right)^{2}
  dx,
\]
Its $L^2$-gradient flow is described, up to a time rescaling, by the system
\begin{equation}\label{eq-flow}
 \left\{\begin{array}{lll}
 u_t &=& \mu +   \sigma_{\e} \left(-\Delta \mu + \frac{1}{\epsilon^2} F''(u) \mu    \right) \\
 \mu &=& \Delta u - \frac{1}{\varepsilon^2} F'(u).
\end{array}
\right.
\end{equation}
The numerical approximation of the solutions will be performed in a square calculation box $Q$ with periodic boundary conditions. \medskip

Our primary objective in this section is to demonstrate that the numerical solutions to the phase field flow starting from the initial condition $u(x,0)= -q'(\dist(\Gamma,x)/\e)$ are accurate approximations, at least in the regular case, of the non-oriented mean curvature flow  $t \mapsto \Gamma(t)$ starting from $\Gamma(0)=\Gamma$. Here, $\dist$ denotes the classical distance function.

Given a positive time discretization parameter $\delta_t$, we now review several classical  schemes used to construct sequences
$(u^n)_{n \geq 0}$ that approximate exact solutions $u(\cdot, n\delta_t)$ of \eqref{eq-flow}.

Due to the homogeneity of the differential operators involved, a fast semi-implicit Fourier spectral method as in \cite{Chen_fourier,BrasselBretin,MR3738845,bretin_droplet,bretin_largephases,BRETIN2018324} can be employed, see also \cite{DU2020425}  for a recent review of numerical methods for phase fields approximations of various geometric flows.

Recall that the Fourier $\boldsymbol K$-approximation of a function $u$ defined in
$Q $ is given by
$$u^{\boldsymbol K}(x) = \sum_{{\boldsymbol k}\in K_N  } c_{\boldsymbol k} e^{2i\pi{\boldsymbol \xi}_k\cdot x},$$
where  $K_N =  [ -\frac{N_1}{2},\frac{N_1}{2}-1 ]\times [ -\frac{N_2}{2},\frac{N_2}{2}-1] \cdots \times   [ -\frac{N_d}{2},\frac{N_d}{2}-1] $,   ${\boldsymbol k} = (k_1,\dots,k_d)$ and ${\boldsymbol \xi_k} = (k_1/\ell_1,\dots,k_d/\ell_d)$. In this expression, the coefficients $c_{\boldsymbol k}$'s denote the first discrete Fourier coefficients of $u$.
The inverse discrete Fourier transform (IFFT) gives the values of $u^{\boldsymbol K}$ at points
$x_{\boldsymbol k} = (k_1 h_1, \cdots, k_d h_d)$, $h_{\alpha} = \ell_{\alpha}/N_{\alpha}$ for $\alpha\in\{1,\cdots,d\}$, more precisely $u^{\boldsymbol K}_{\boldsymbol k} =   \textrm{IFFT}[c_{\boldsymbol k}]$. Conversely,
$c_{\boldsymbol k}$ can be computed as the discrete Fourier transform of $u^{\boldsymbol K}_{\boldsymbol k},$ i.e. $c_{\boldsymbol k} =
\textrm{FFT}[u^{\boldsymbol K}_{\boldsymbol k}]$.\medskip

Semi-implicit approaches typically treat explicitly the nonlinear terms of phase field models, which requires
 stabilization techniques to ensure computational efficiency.

In this context, a well-known approach proposed by Eyre in \cite{MR1676409} uses a convex-concave splitting of the Cahn-Hilliard energy. This method provides a simple, effective, and stable scheme for approximating various evolution problems with gradient flow structures, as illustrated in \cite{MR2418360, MR2519603, MR2799512, MR3100769, MR3564350, MR3682074, Schonlieb2011}. More recent approaches stabilize these schemes by introducing relaxation techniques based on the \textit{Scalar Auxiliary Variable} (SAV) method. These methods have been used in a wide variety of studies \cite{HOU2019307, SHEN2018407, doi:10.113717M1150153, doi:10.113719M1264412, Huang2020} to address a broad class of gradient flow systems, particularly the Cahn-Hilliard equation \cite{wang2022application, huang2023structurepreserving, yang, MR4696121}, dissipative systems \cite{bouchriti2020remarks, zhang2022generalized}, and their variants \cite{ESAV, ExpSAV}.

In the remainder of this paper, we will focus exclusively on the convex-concave splitting approach, which is simpler to implement. \medskip

 The code used for the numerical experiments described in the next paragraphs is  available on a GitHub\footnote{\url{https://github.com/eliebretin/UMCF.git}}
 repository.

\subsection{A convex-concave Fourier spectral numerical method}

 Consider the gradient flow of an energy  $\E_\e$
$$ \e \partial_t u  = - \nabla_u \E_\e(u),$$
and decompose  $\E_\e$ as the sum of a convex energy $\E_{\e,c}$ and a concave energy $\E_{\e,e}$ (a classical explicit decomposition will be proposed later):
$$\E_\e= \E_{\e,c} + \E_{\e,e}.$$
An effective numerical scheme associated with this decomposition can be designed by integrating implicitly the gradient of the convex part, and explicitly the gradient of the concave part:
$$ \e \frac{u^{n+1} - u^{n}}{\delta_t} = - \left( \nabla_u  \E_{\e,c}(u^{n+1}) +  \nabla_u  \E_{\e,e}(u^n)\right).$$
The energy stability can be easily proven by interpreting this scheme as one step, starting from
$u^n$, of the implicit discretization of the semi-linearized PDE   $\partial_t u =  -\nabla_u \overline{\E}_{\e,u^{n}}(u)$
 where $\overline{\E}_{\e,u^{n}}$ is  defined as
 $$
\overline{\E}_{\e,u^{n}}(u)  =\E_{\e,c}(u) + \E_{\e,e}(u^n) +
\langle \nabla_u  \E_{\e,e}(u^n)  , (u - u^n) \rangle.
$$
It is easily seen that $\overline{\E}_{\e,u^{n}}(u^{n+1}) \leq \overline{\E}_{\e,u^{n}}(u^{n})$. Using the concavity of $\E_{\e,e}$ it can be proved that the energy $\E_\e$ decreases along the iterations:
$$\E_\e(u^{n+1}) \leq \E_\e(u^{n}).$$

For our phase field model
\[
  \E_\e(u,Q) = \int_Q \left(\frac{\e}{2}|\nabla u|^2+\frac{1}{\e}F(u)\right) dx
  + \frac{\sigma_\e}{2 \e}  \int_Q  \left( \e \Delta u - \frac{1}{\e}F'(u)\right)^{2}
  dx,
\]
we define, using $\alpha,\beta>0$, the convex energy
$$  \E_{\e,c}(u,Q) = \frac\e 2 \int_Q \left(|\nabla u|^2 +  \sigma_\e (\Delta u)^2   + {\alpha} u^2 +   {\beta}|\nabla u|^2\right) dx.$$
For sufficiently large coefficients  $\alpha,\beta$,  the energy
$$ \E_{\e,e}(u,Q) = \frac \e 2\int_Q \left( \frac{2}{\e^2}F(u) +  {\sigma_{\e}}\left( ( \Delta u - \frac{1}{\e^2}F'(u))^{2} - (\Delta u)^2 \right) - {\alpha} u^2 -   {\beta}|\nabla u|^2\right)  dx,
$$
is concave,  at least on a subset of $\wpq{2}{2}(Q)$  where the energy $\E_{\e}(\cdot, Q)$ is bounded. The numerical
scheme associated with the decomposition $\E_\e= \E_{\e,c} + \E_{\e,e}$ is given by 
$$u^{n+1} = L[g(u^n)],$$
 where $g$ is a nonlinear operator defined by
$$ g(u) =  u + \delta_t \left( - \frac{F'(u)}{\e^2}  +  \sigma_\e \left(- \Delta  [F'(u)/\e^2] + \frac{F''(u)}{\e^2} \left[ \Delta u - F'(u)/\e^2 \right] \right)  + \alpha u - \beta \Delta u\right), $$
and $L$ is the homogeneous linear operator $\left( I_d +  \delta_t (- \Delta + \sigma_\e \Delta^2 + \alpha I_d - \beta \Delta  ) \right)^{-1}$, whose Fourier symbol is
 $$\Hat{L}(\xi) = 1/(1 + \delta_t ( 4 \pi^2 |\xi|^2 + \sigma_\e 16 \pi^4 |\xi|^4 + \alpha + \beta 4 \pi^2 |\xi|^2 ).$$

\begin{remk}
In the scheme above, we have not taken into account the singularity of $F$ for values greater than $1/4$. A straightforward way to address this issue and incorporate it into our scheme is by introducing a projection step:
  $$u^{n+1} = \min(L[g(u^n)],1/4).$$
\end{remk}

\begin{remk}
Using the parameters $\alpha, \beta$ in the convex-concave splitting ensures high stability of the numerical scheme. However, it generally results in reduced accuracy. As the main purpose of the numerical results presented in the next sections is to compare precisely our flow with the mean curvature flow, we need high accuracy. So we opt for the settings $\alpha=\beta=0$ to avoid accuracy losses and a time step $\delta_t=0.01\e^2$ significantly below the natural stability condition of the PDE. Readers who are not so demanding in terms of precision, but who are interested in stable simulations obtained with larger time steps might opt for positive values of $\alpha,\beta$, typically $\alpha=\frac1{\e^2}$, $\beta=1$, and a time step $\delta_t=\e^2$.
\end{remk}

\subsection{$2D$ numerical simulations}

We first focus on evolving interfaces confined in $Q=[0,1]^2$ discretized with $N$ nodes in both directions. 

\subsubsection{Evolution of a circle and comparison with mean curvature flow }

A circle of initial radius $R_0$ flowing by mean curvature remains a circle with radius given by
 $$R(t) = \sqrt{R_0^2 - 2t}$$
until the extinction time $t_{end} = \frac{1}{2} R_0^2$.  In our first numerical experiment, we compare this evolution with the flow provided by our phase field model for two different values of $\e$.

We set  $N=2^8$, $\e = 1.5/N$ or $\e = 3/N$, $\sigma_\e = 4 \e^2$,  $\alpha = \beta = 0$ and $\delta_t = 0.01 \epsilon^2$.   This time step setting, well below the natural stability condition, is more suitable to illustrate numerically the convergence order of the phase field model with respect to $\e$, while drastically reducing the numerical discretization errors.

\begin{figure}[!htbp]
\centering
	\includegraphics[width=.245\textwidth]{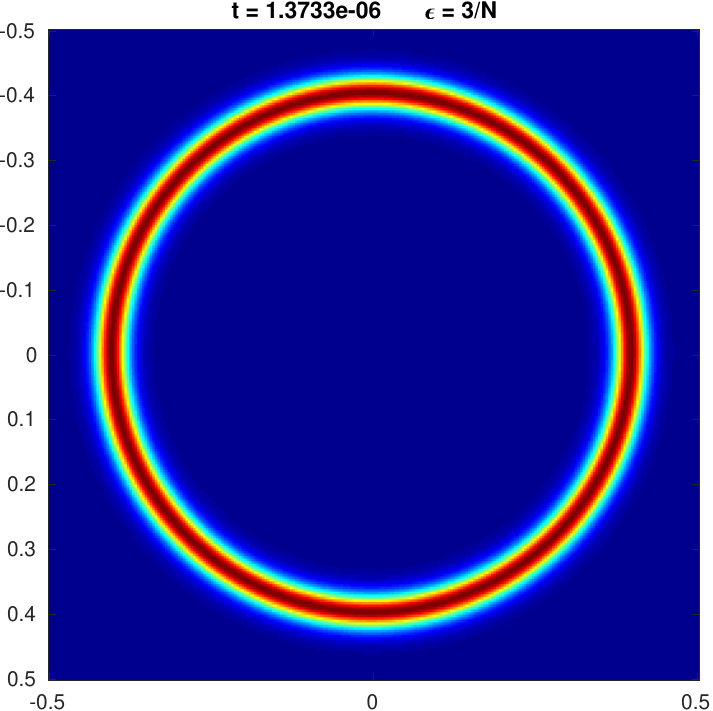}
    \includegraphics[width=.245\textwidth]{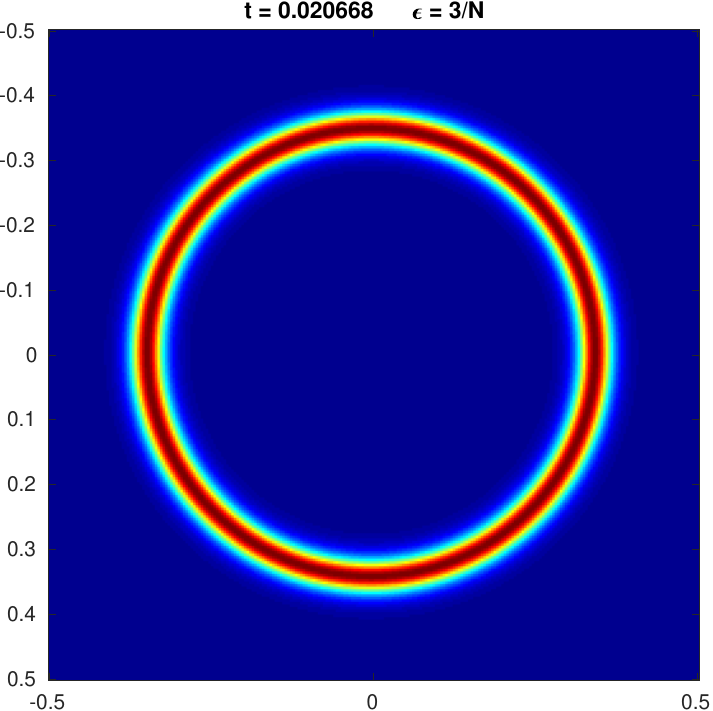}
    \includegraphics[width=.24\textwidth]{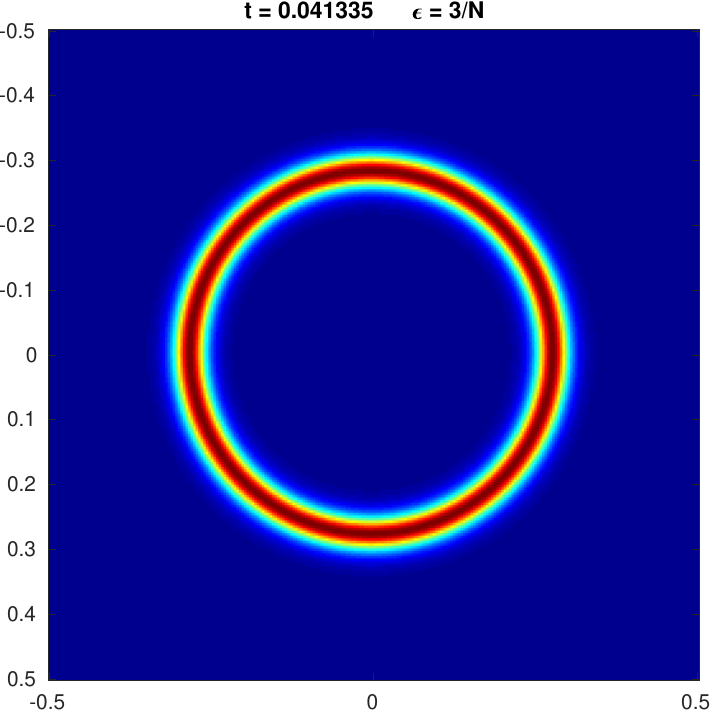}
    \includegraphics[width=.24\textwidth]{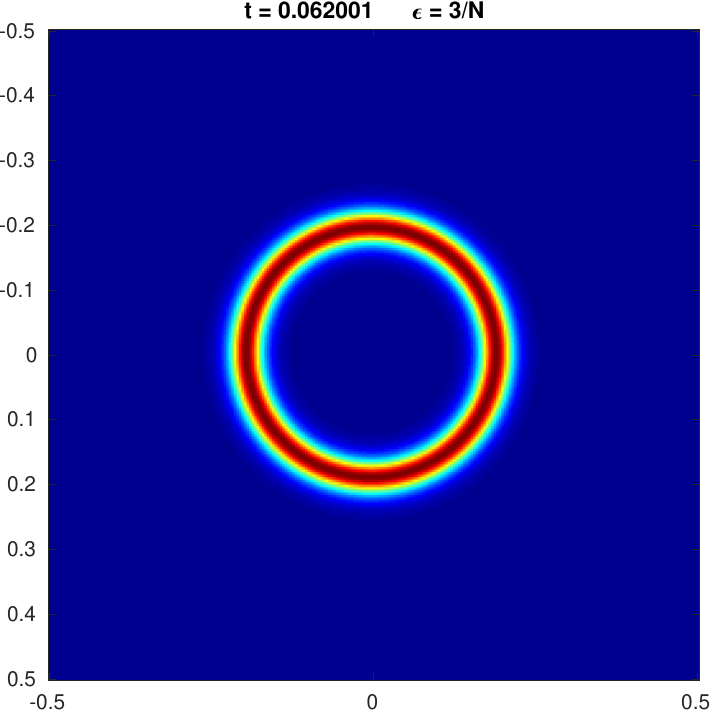}\\
	\includegraphics[width=.24\textwidth]{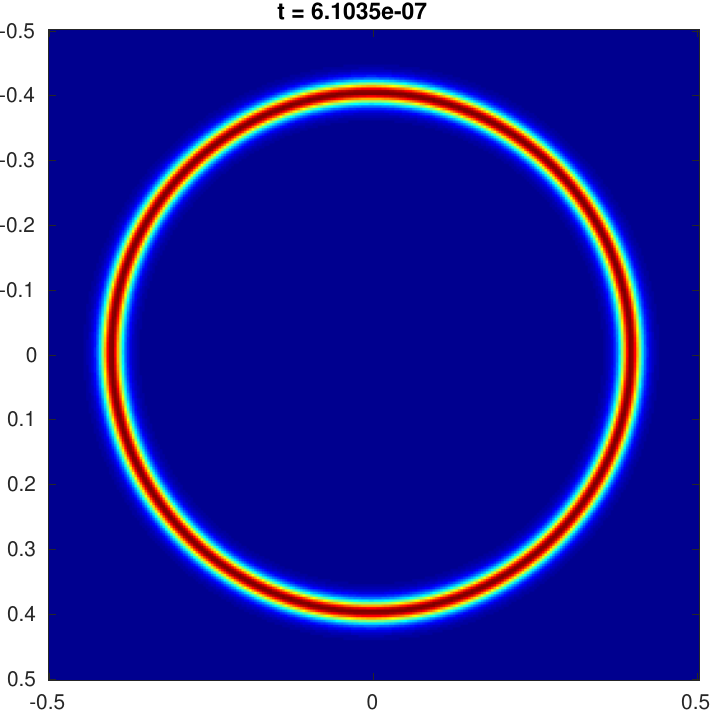}
    \includegraphics[width=.24\textwidth]{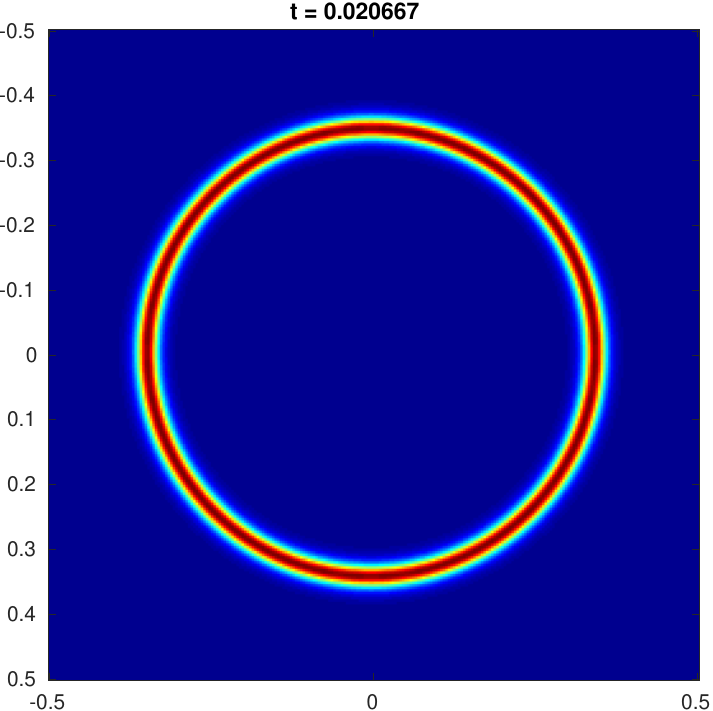}
    \includegraphics[width=.24\textwidth]{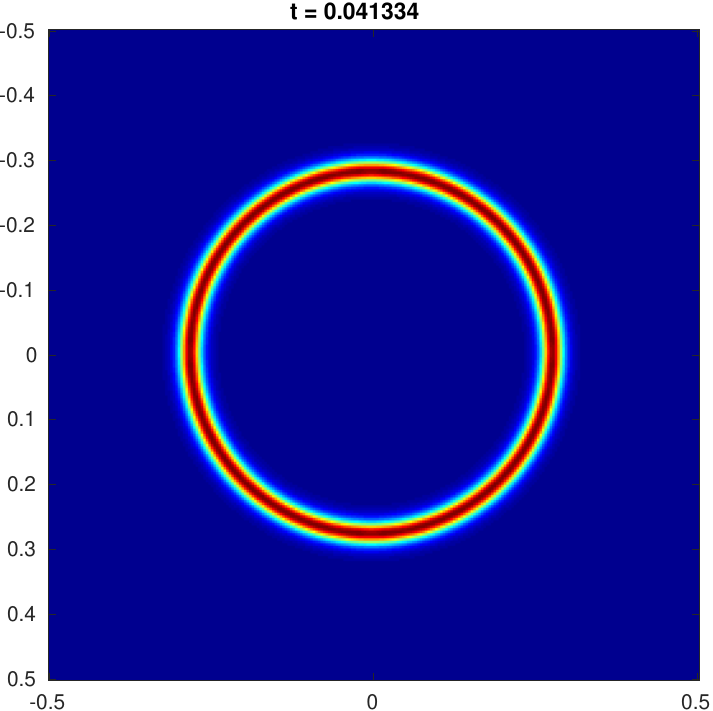}
    \includegraphics[width=.24\textwidth]{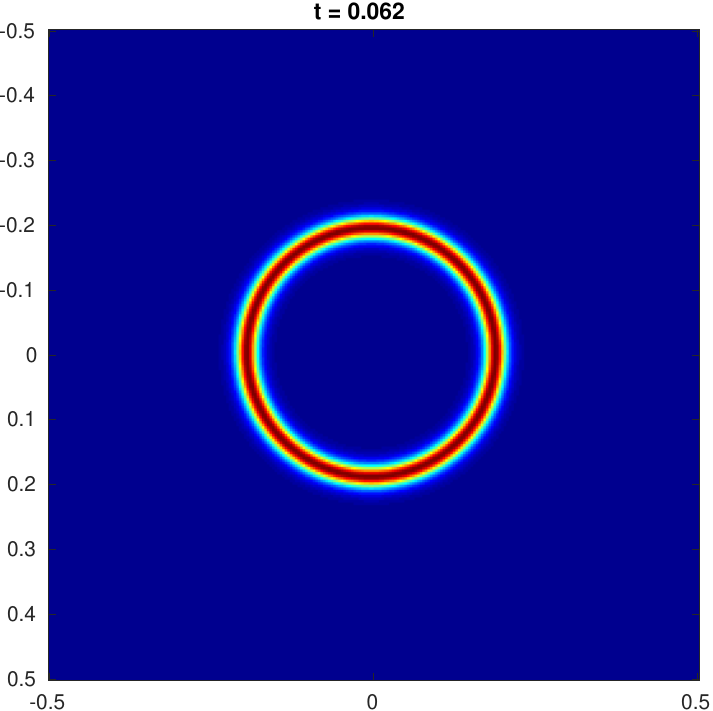}
\caption{Approximate evolution of a circle across iterations using two values of $\e$: first line with $\e = 3/N$, second line with $\e = 1.5/N$.}
\label{fig_test_circle1}
\end{figure}

\begin{figure}[!htbp]
\centering

	\includegraphics[width=7cm]{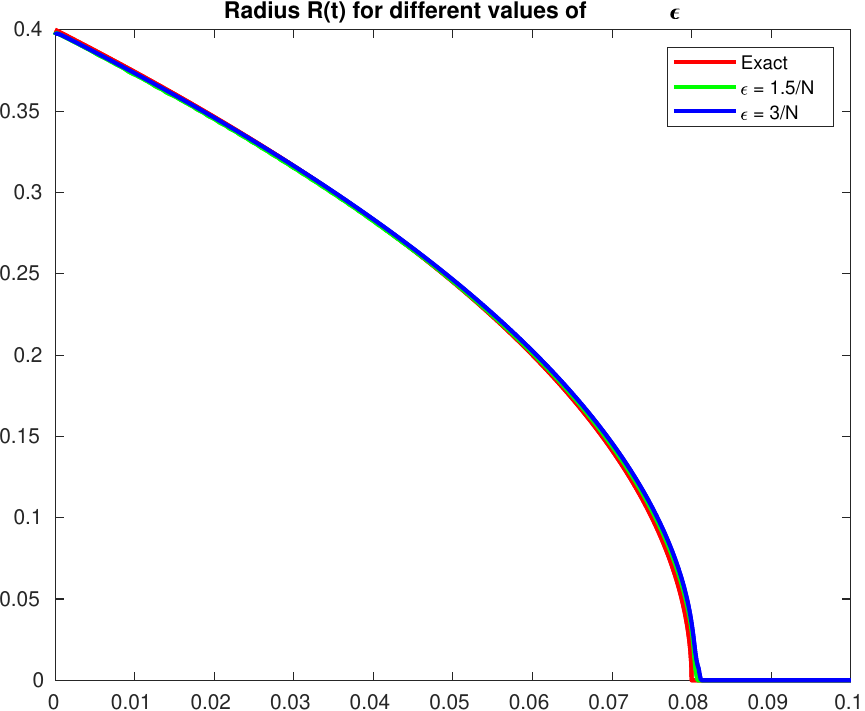}
    \includegraphics[width=7cm]{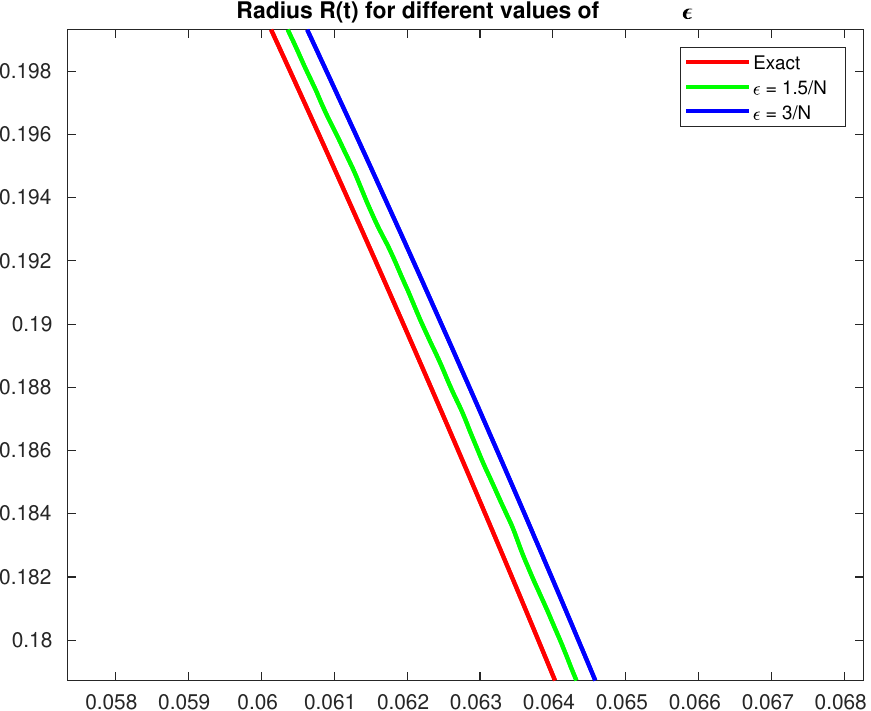}
\caption{Comparison of the numerical radius of the circle evolved with our model vs the exact radius $t\mapsto R(t)$ provided by mean curvature flow. Left: on the time interval $[0,0.1]$. Right: zoom on the time interval $[0.058,0.068]$.}
\label{fig_test_circle2}
\end{figure}

Figure  \ref{fig_test_circle1} shows several iterates $u^n$ for the setting $\e = 1.5/N$ and $\e = 3/N$.   The influence of the parameter $\e$  on the size of the diffuse interface is clearly observable.
 Furthermore, the initial circle appears to evolve in a consistent manner, with its radius decreasing over the iterations. To compare the radius law with the theoretical law of motion by mean curvature, we plot in Figure  \ref{fig_test_circle2}
 a comparison between the different evolutions $t \mapsto R(t)$. Specifically, the approximate radius
 $R_{\e}(t)$ is estimated from the phase field function $u_\e(t)$ using the formula
$$ R_\e =  \frac{1}{2 \pi \e} \int_Q u_\e dx.$$

Notably, the left side of Figure   \ref{fig_test_circle2} shows that the two laws obtained with
$\e = 1.5/N$ and $\e = 3/N$ are very close to the theoretical one. A zoom on the interval $[0.058,0.068]$ shown
on the right side of Figure
\ref{fig_test_circle2} indicates that the observed error is approximately twice as significant when using  $\e = 3/N$ compared to
$\e = 1.5/N$.  This suggests that our phase field model should converge to the mean curvature flow with an error of order
  $O(\e)$. In contrast, this error is typically of order $O(\e^2)$ for the Allen-Cahn equation. \medskip

  This numerical result goes well beyond our initial theoretical expectations: we expected the phase field profile to be stable for sufficiently large choices of $\sigma_\e$, but we did not expect the approximate flow to approximate well the mean curvature flow. To theoretically justify such a result, we would need to adapt the asymptotic developments of Allen-Cahn solutions to the context of non-oriented phase field functions.

\subsubsection{Influence on stability of parameter $\sigma_\e$}

The second numerical experiment aims to demonstrate the influence of parameter
$\sigma_\e$  on the stability of interfaces evolving by our phase field flow. We have already justified in the analysis of the model the importance of
$\sigma_\e$  as  $\e$ approaches $0$. Here, we aim to illustrate numerically that the choice of this parameter is not trivial and must be made with care. Specifically, a parameter that is too small may fail to ensure the stability of the interface, while a parameter that is excessively large can increase the computational cost and amplify the contribution of the Willmore term in the flow.

\begin{figure}[!htbp]
\centering
    \includegraphics[width=.24\textwidth]{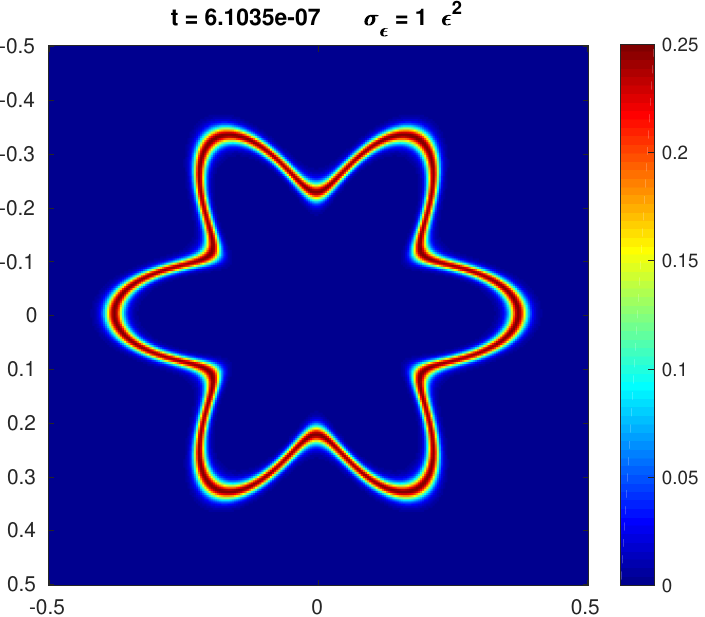}
    \includegraphics[width=.24\textwidth]{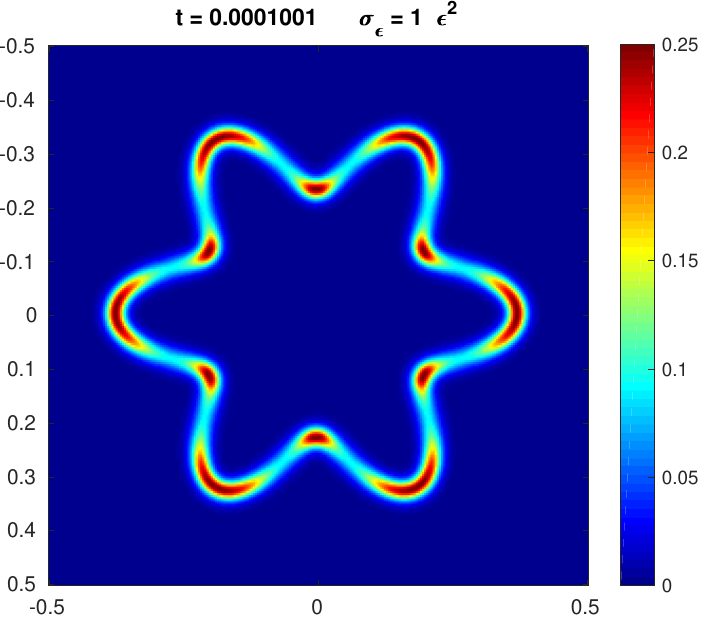}
    \includegraphics[width=.24\textwidth]{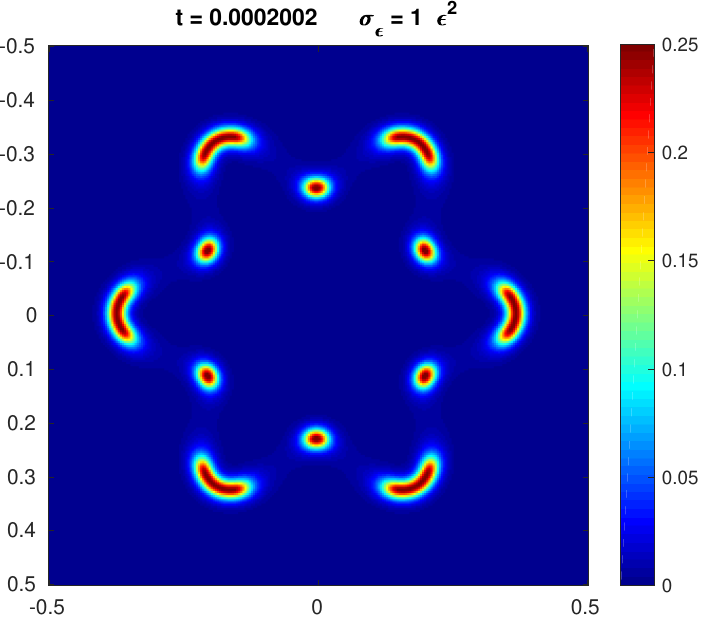}
    \includegraphics[width=.24\textwidth]{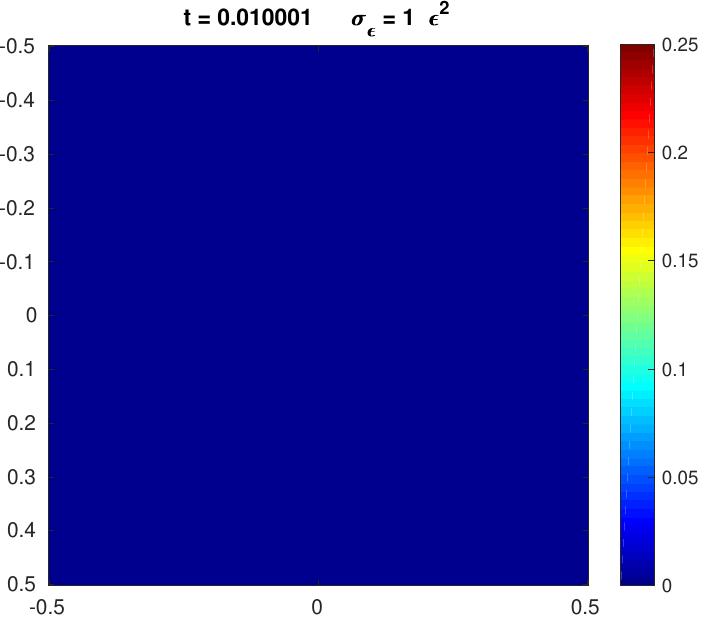}\\
    \includegraphics[width=.24\textwidth]{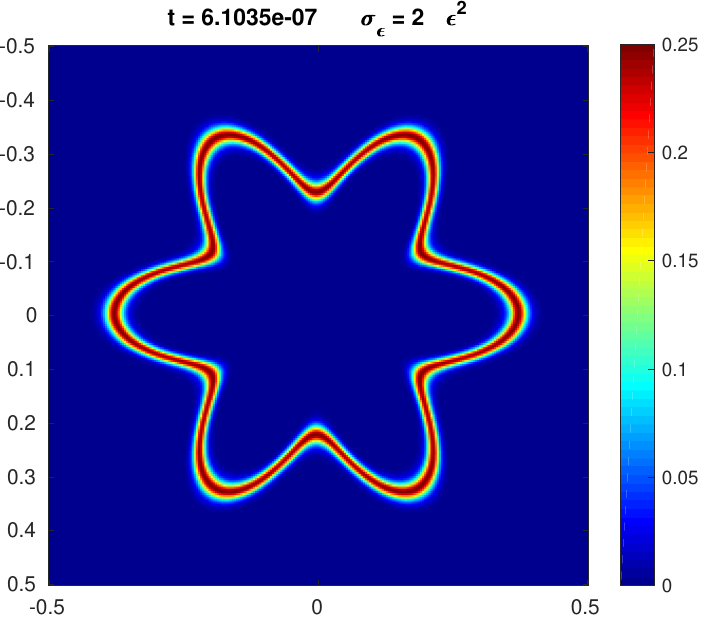}
    \includegraphics[width=.24\textwidth]{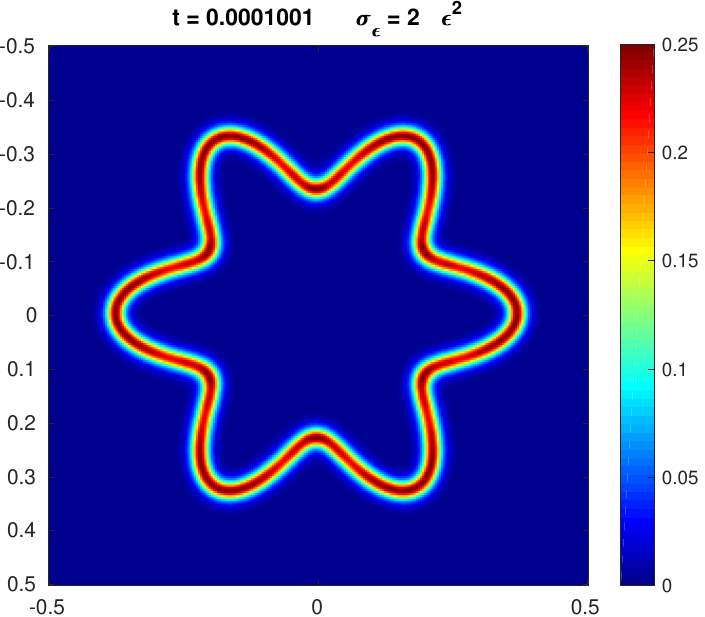}
    \includegraphics[width=.24\textwidth]{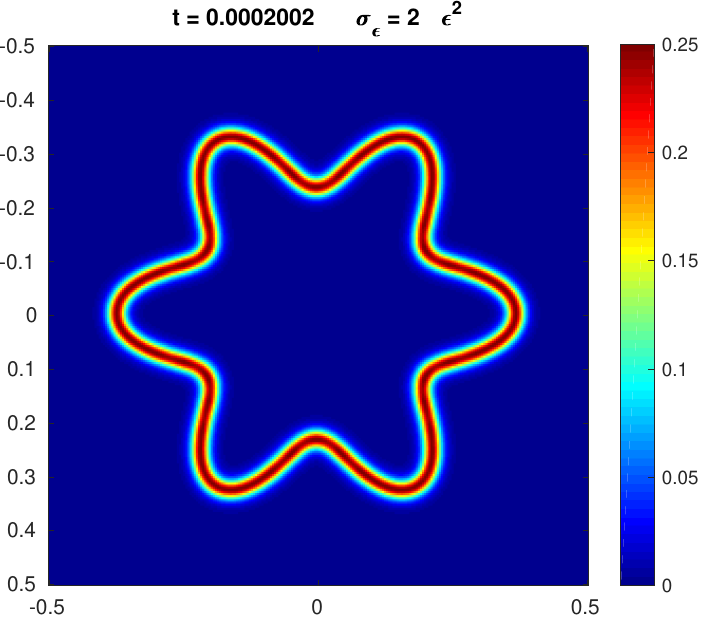}
    \includegraphics[width=.24\textwidth]{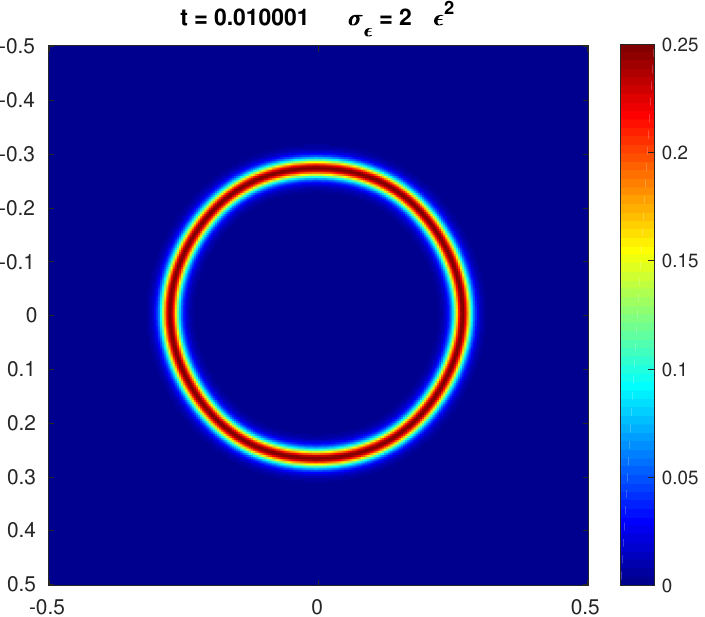} \\
    \includegraphics[width=.24\textwidth]{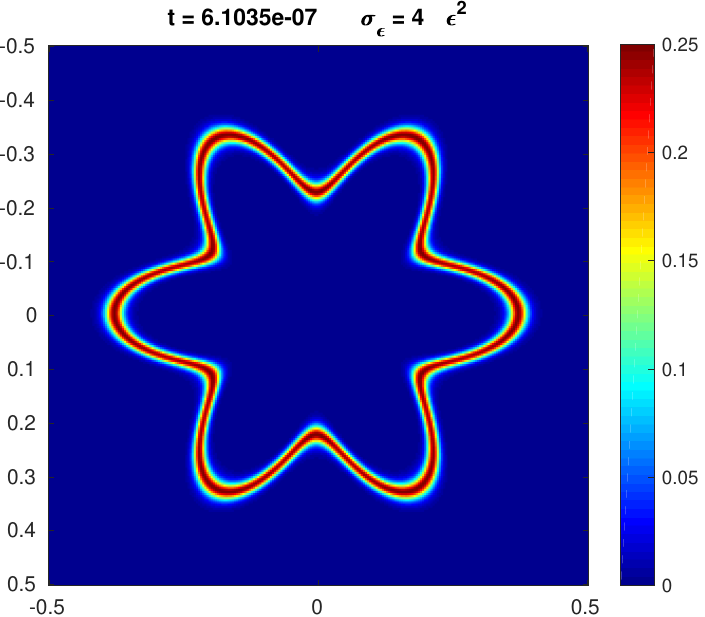}
    \includegraphics[width=.24\textwidth]{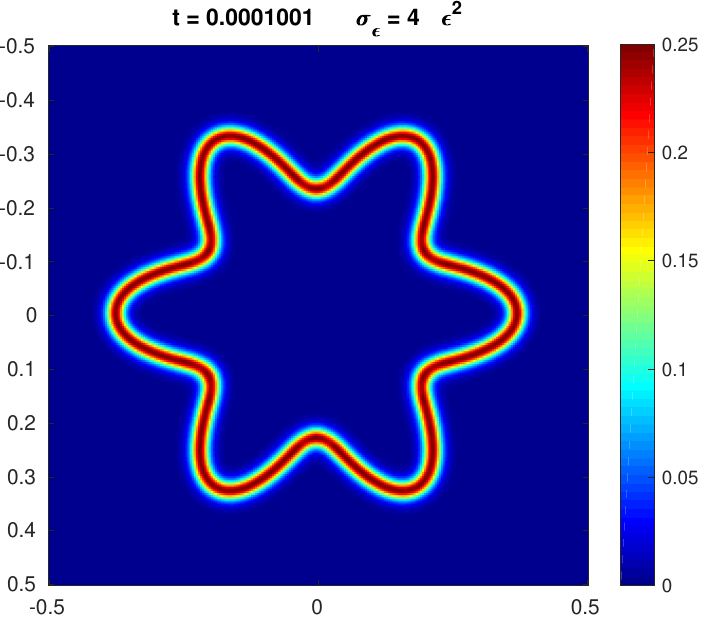}
    \includegraphics[width=.24\textwidth]{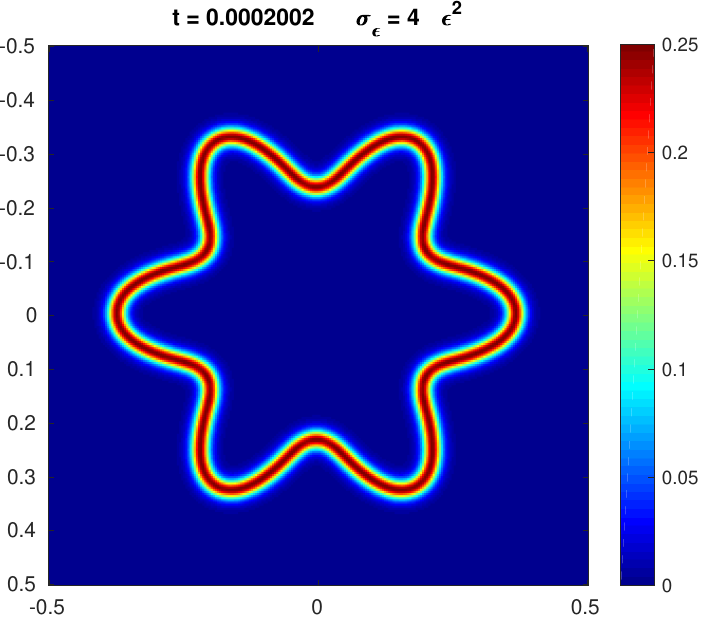}
    \includegraphics[width=.24\textwidth]{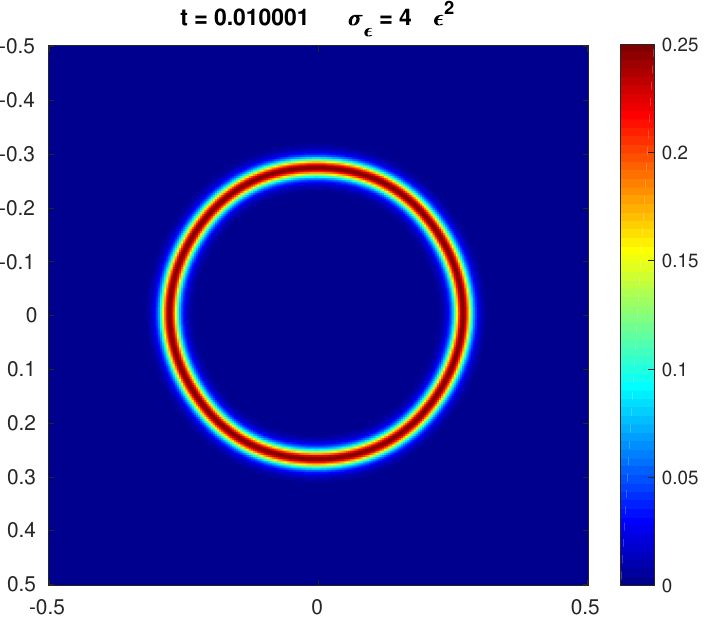}

\caption{Evolution of a shape for three values of $\sigma_\e$. First line: $\sigma_\e = \e^2$, second line: $\sigma_\e = 2 \e^2$, third line: $\sigma_\e = 4 \e^2$.}
\label{fig_test_sigma_e}
\end{figure}

Figure \ref{fig_test_sigma_e} illustrates three different evolutions
 $t \to u_\e(t)$ obtained using  $\sigma_\e = \e^2$, $\sigma_\e = 2 \e^2$ and   $\sigma_\e = 4 \e^2$.
 The other parameters remain constant:  $N = 2^8$, $\e = 2/N$, $\delta_t = 0.01 \e^2$ and $\alpha = \beta = 0$.

We first observe that the choice
 $\sigma_\e = \e^2$
  does not ensure the stability of the interface. Surprisingly, the interfaces with low curvature appear to be the most unstable.
  
The second row of Figure \ref{fig_test_sigma_e}, corresponding to $\sigma_\e = 2 \e^2$, shows an evolution that looks better. However, upon examining the second figure more closely, we note that the profile does not fully reach the value of $1/4$ when the curvature is zero, indicating that the evolution is not completely stable. In contrast, the choice $\sigma_\e = 4 \e^2$  seems to ensure, at least in this example, a satisfactory diffuse profile across the entire interface long the iterations.

\subsubsection{Non smooth interfaces and evolution of triple junctions}

\begin{figure}[!htbp]
\centering
    \includegraphics[width=.24\textwidth]{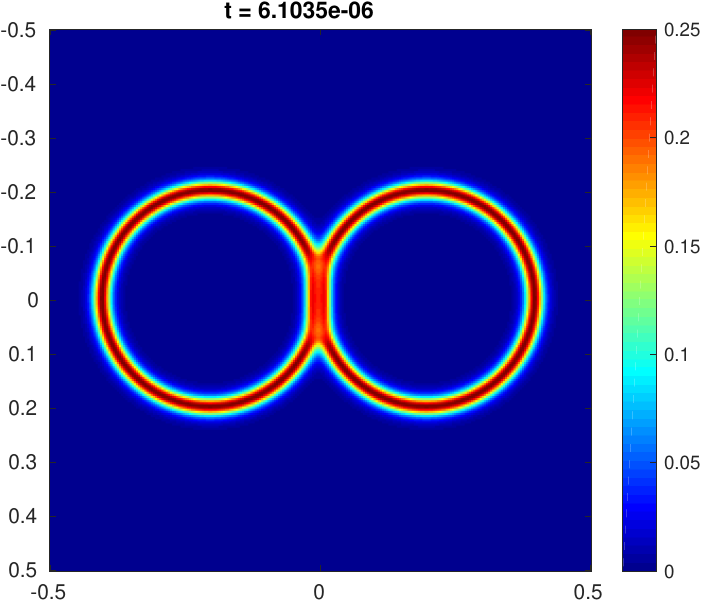}
     \includegraphics[width=.24\textwidth]{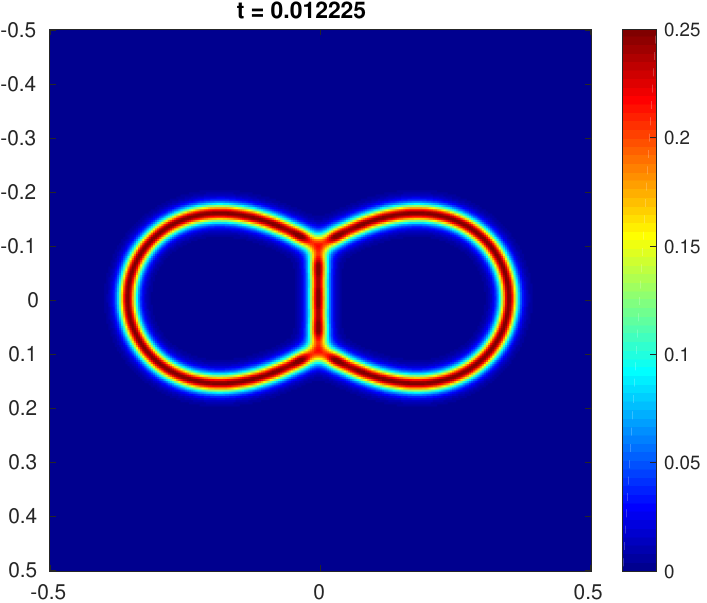}
     \includegraphics[width=.24\textwidth]{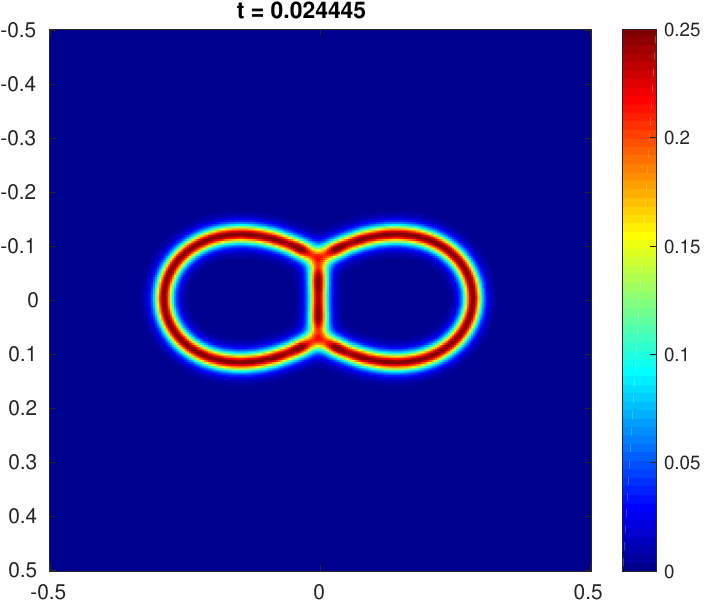}
     \includegraphics[width=.24\textwidth]{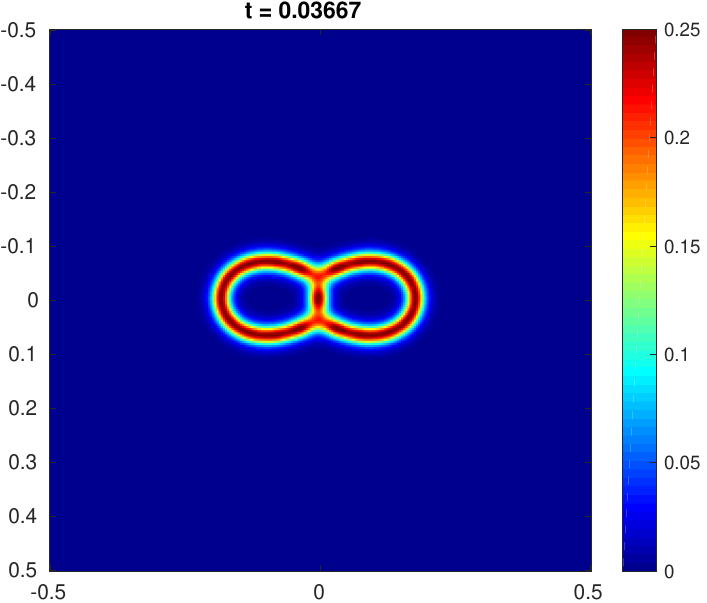} \\

    \includegraphics[width=.24\textwidth]{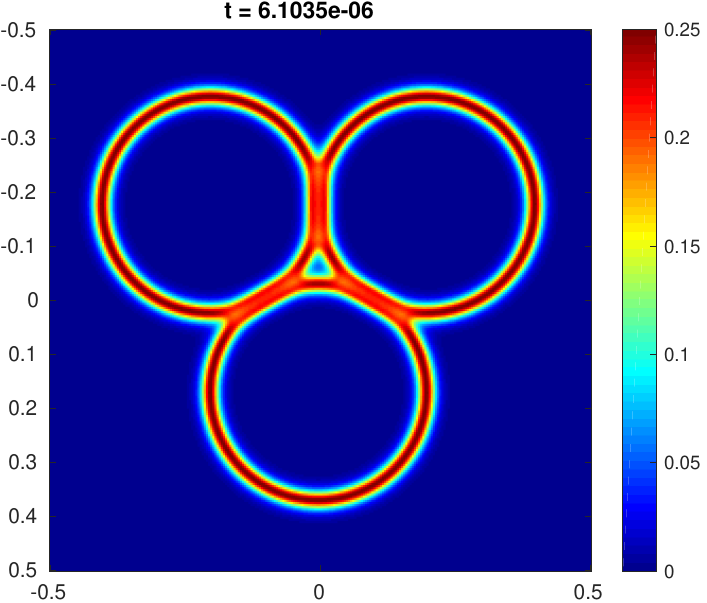}
     \includegraphics[width=.24\textwidth]{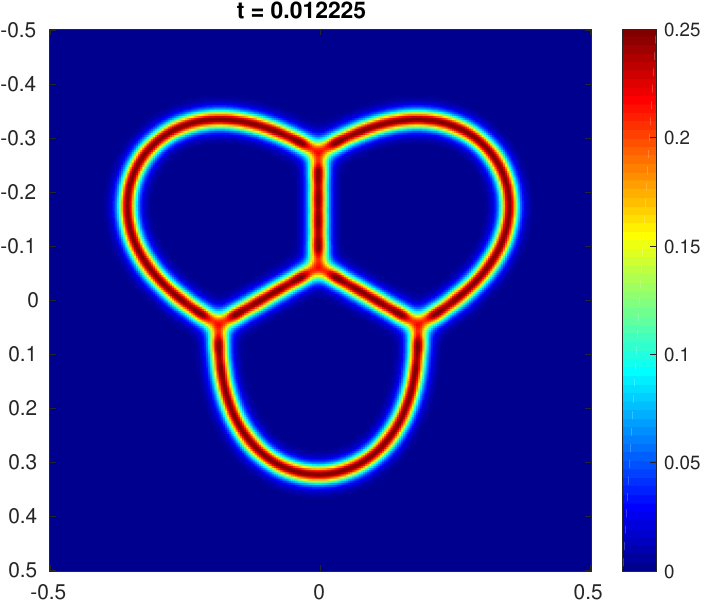}
     \includegraphics[width=.24\textwidth]{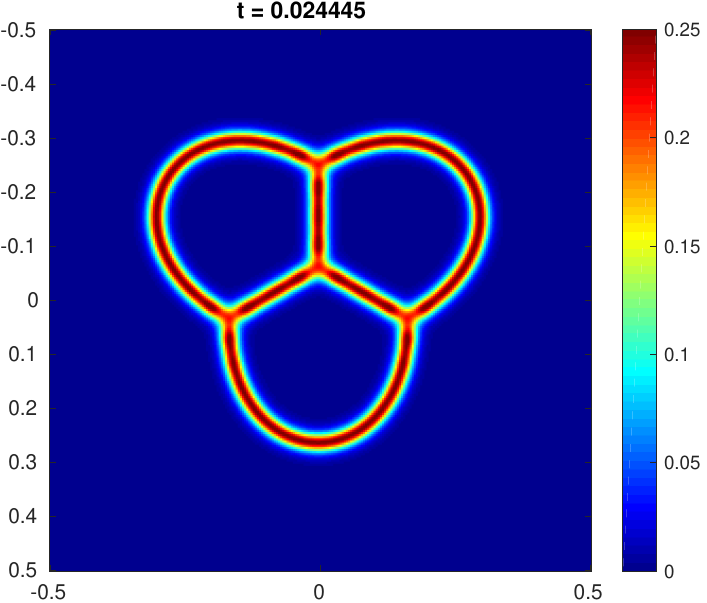}
     \includegraphics[width=.24\textwidth]{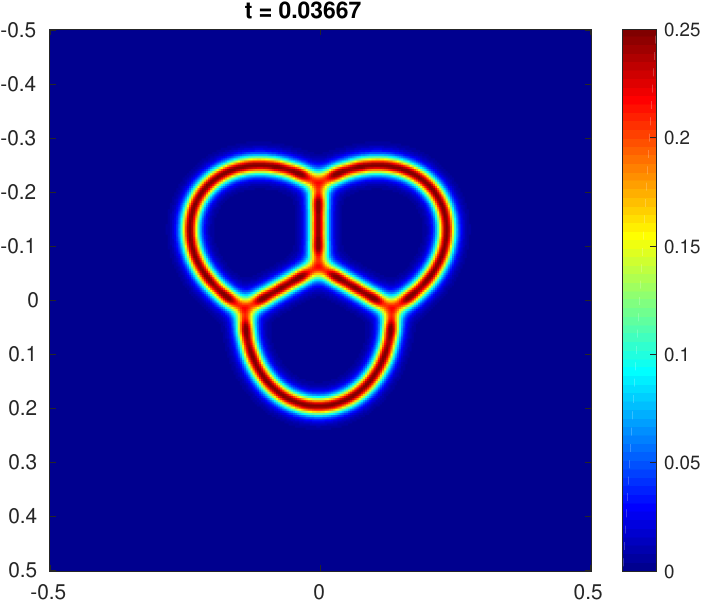} 
\caption{Non-smooth evolving shapes. Each row corresponds to a different initial condition. The figures show the phase field across iterations.
}
\label{fig_test_non_smooth}
\end{figure}

The aim of this next numerical experiment is to demonstrate that the proposed phase field model can effectively handle the evolution of singular interfaces with triple junctions. We consider two examples where the initial condition involves either two or three glued circles. The mean curvature flow of these shapes cannot be approximated with the classical Allen-Cahn equation, it is necessary to use a multiphase Allen-Cahn model where the interior of each circle represents a distinct phase. Across iterations produced by this multiphase model, the angles between phases at triple points evolve toward Herring's equiangular equilibrium configuration. The numerical simulations obtained with our approach are presented in Figure~\ref{fig_test_non_smooth}, revealing stable triple points that exhibit a quasi-symmetry consistent with Herring's law. The numerical results closely resemble the outcomes that a multiphase model would yield. For this experiment, we set  $N=2^8$, $\e = 2/N$,
$\sigma_\e = 4 \e^2$, $\delta_t = 0.1 \e^2$ and $\alpha = \beta = 0$.

\subsection{$3D$ numerical simulations}
As in the two-dimensional case, we first examine the evolution of a sphere. The other experiments illustrate the influence of the Willmore term in our model. It not only ensures the stability of the interface, it also plays a significant role when the interface undergoes topological changes, as in the case of a dumbbell flowing by mean curvature. An interesting aspect of the Willmore term is its ability to stabilize fine structures, such as tubular sets, which is not achievable with conventional mean curvature approximation models. We will illustrate this point by showing how our phase field model can approximate well the mean curvature flow of thin structures of codimension two, such as circles or, more generally, tubular sets.

In all subsequent numerical experiments, we display the $\frac 16$-level set
 of each approximate solution $u_\e$ to represent the evolving interface associated to the phase field function $u_\e$ (recall that the potential $F$ satisfies $F'(\frac 1 6)=0)$.

\subsubsection{Evolution of a sphere and comparison with mean curvature flow}

\begin{figure}[!htbp]
\begin{center}
\hspace*{-0.4cm}
	\includegraphics[width=.26\textwidth]{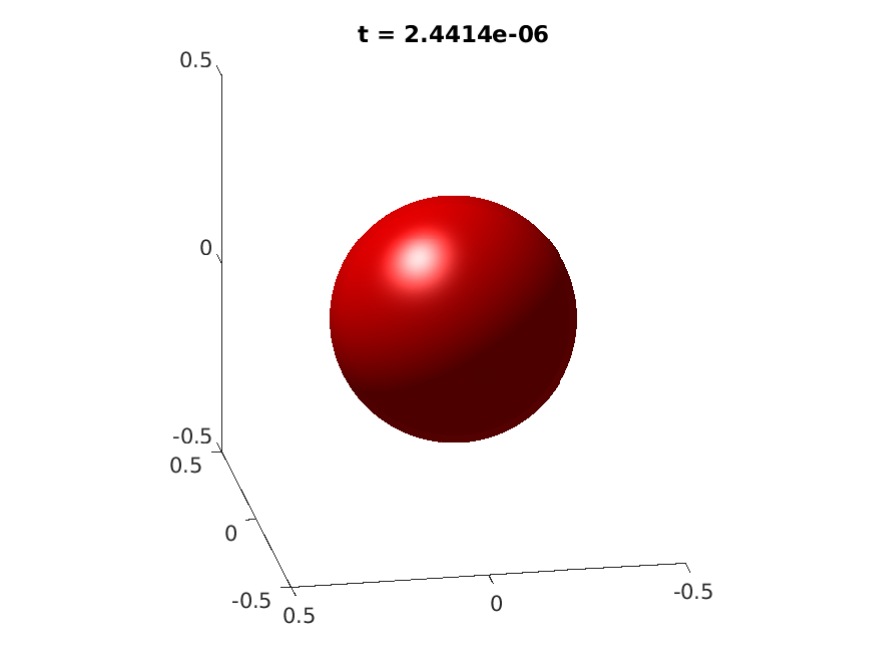}\hspace*{-0.4cm}
	\includegraphics[width=.26\textwidth]{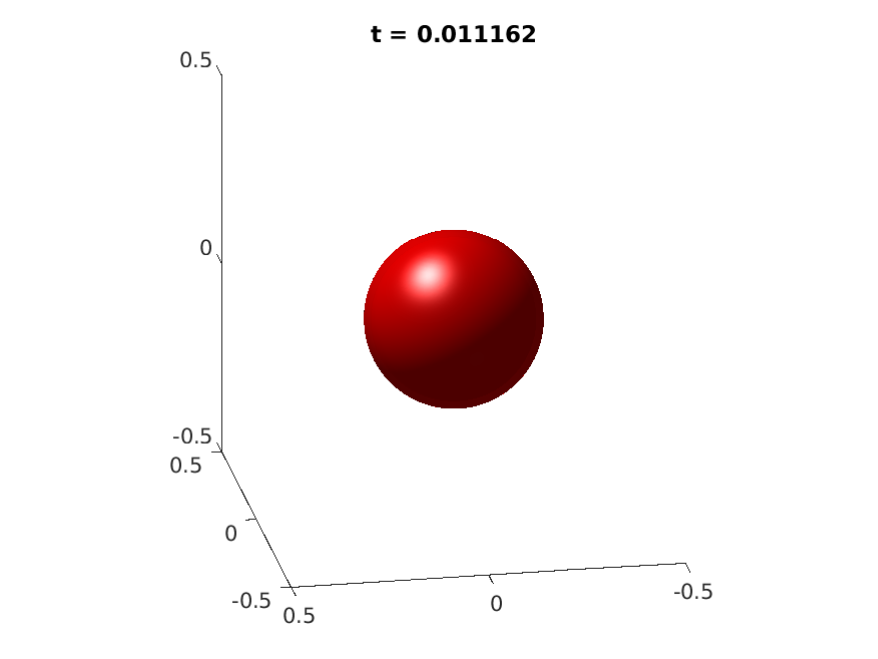}\hspace*{-0.4cm}
	\includegraphics[width=.26\textwidth]{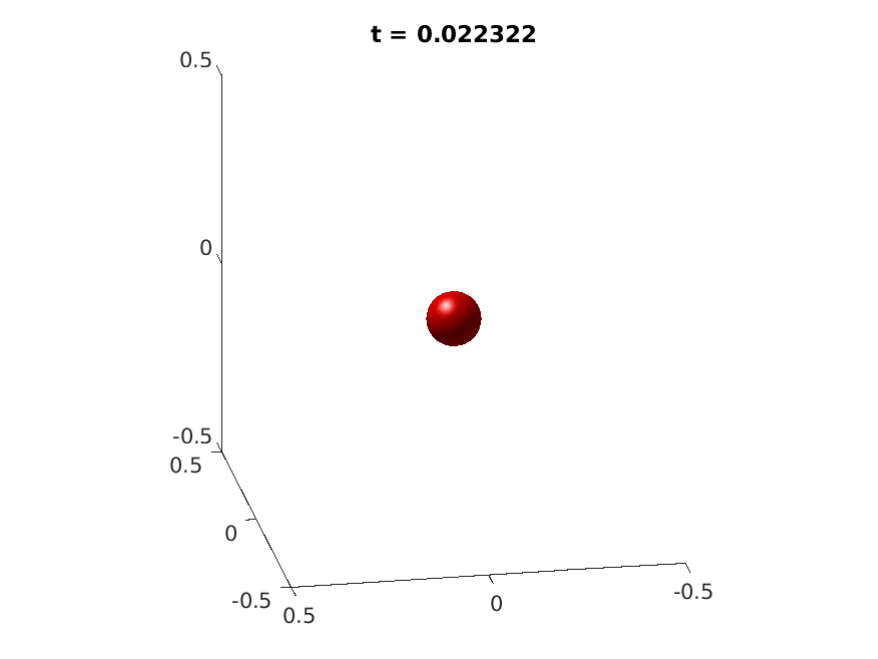}\hspace*{-0.4cm}
	\includegraphics[width=.24\textwidth]{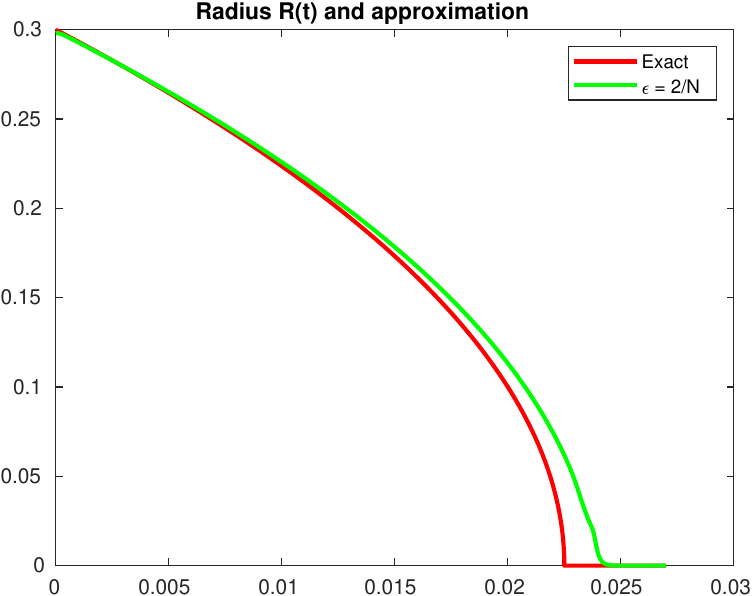}

\caption{Evolution obtained with our model starting from an initial sphere, and comparison of the numerical radius with the radius of a sphere evolving exactly by mean curvature.}
\label{fig_test_sphere}
\end{center}
\end{figure}

An initial sphere of radius $R_0$ flowing by mean curvature remains a sphere whose radius satisfies the law $R(t) = \sqrt{R_0 - 4 t}$. In the first three images of Figure \ref{fig_test_sphere}, we display 
 the $\frac 16$-level set of the phase field function $u_\e$ provided by our numerical model at three different iterations
The computation was made using the numerical parameters $N = 2^7$, $\e = 2/N$, $\sigma_\e = 2 \e^2$, $\delta_t = 0.01 \e^2$ and $\alpha = \beta = 0$. The evolving shape closely resembles a sphere with radius decreasing over time. 
The last image of Figure \ref{fig_test_sphere} plots the numerical radius  of the shapes computed across iterations using
the formula $ R_\e = \sqrt{\frac{1}{4 \pi \e} \int u_\e dx}$.  The evolution of the numerical radius closely aligns with the
theoretical law, although the error appears to be larger than in the two-dimensional case. This could be due the parameter settings, in particular the resolution is coarser than the resolution used in $2D$.

\subsubsection{Evolution of a dumbbell and influence of the Willmore term}

\begin{figure}[!htbp]
\centering
	\includegraphics[width=.24\textwidth]{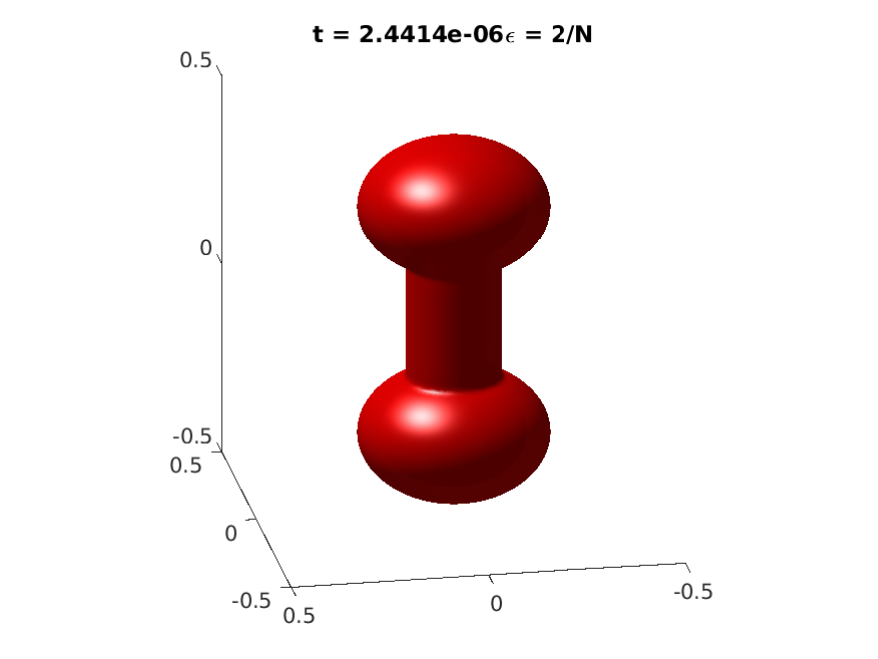}
    \includegraphics[width=.24\textwidth]{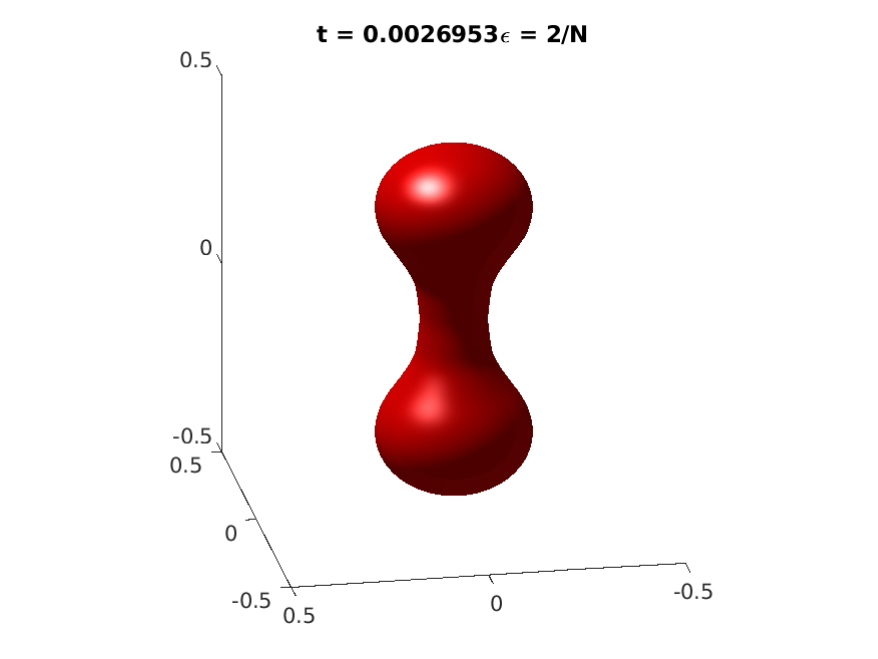}
    \includegraphics[width=.24\textwidth]{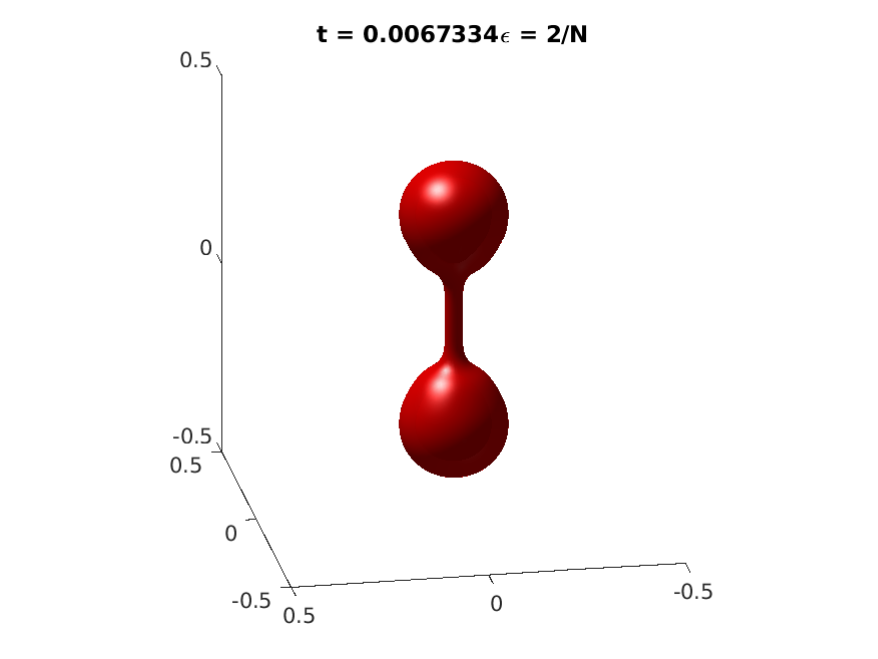}
    \includegraphics[width=.24\textwidth]{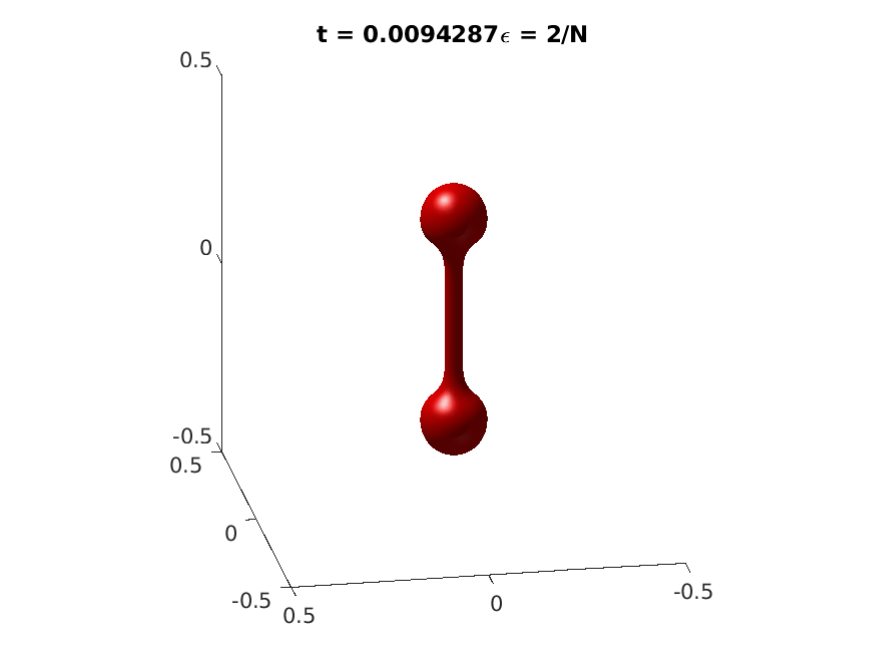} \\
	\includegraphics[width=.24\textwidth]{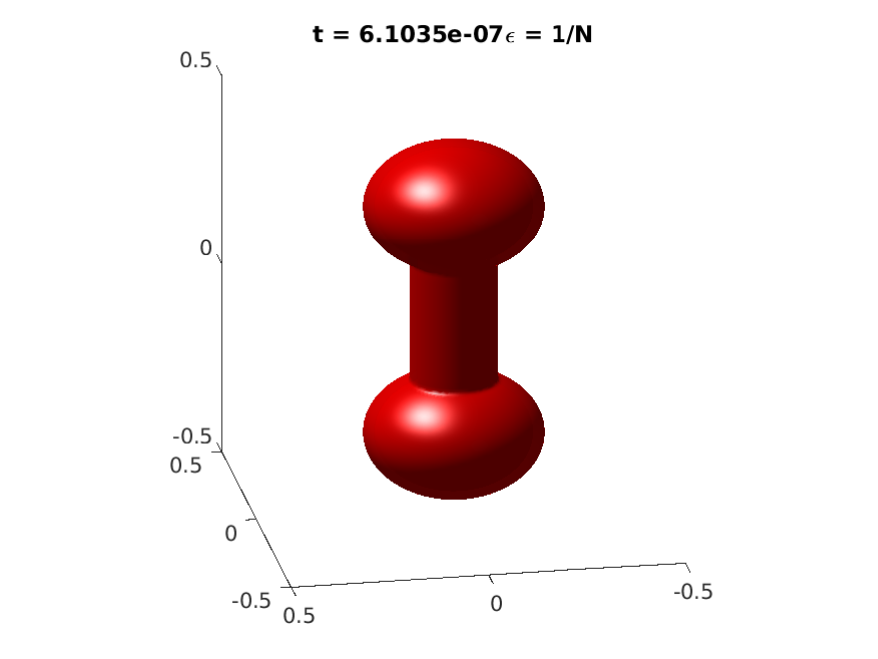}
    \includegraphics[width=.24\textwidth]{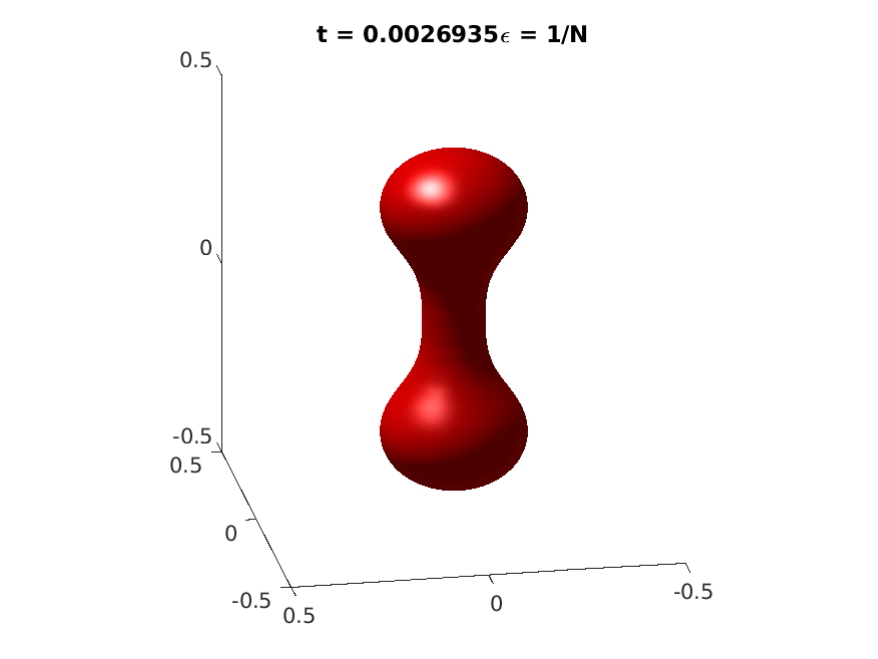}
    \includegraphics[width=.24\textwidth]{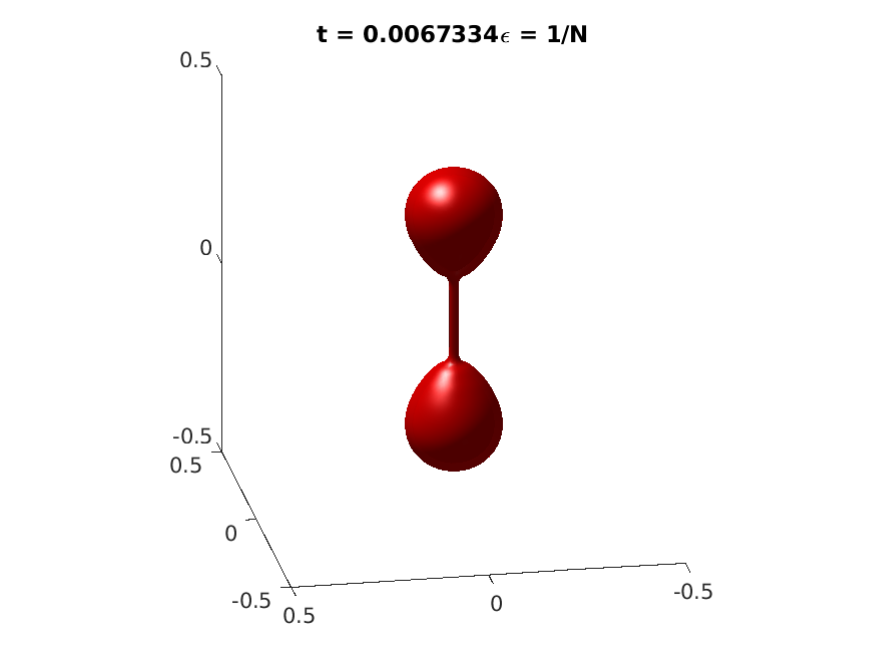}
    \includegraphics[width=.24\textwidth]{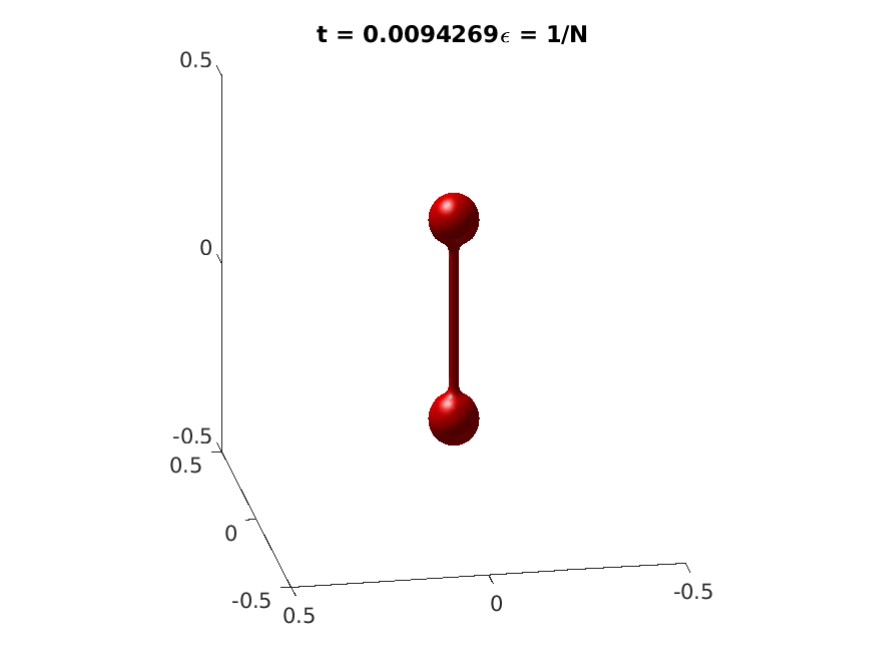} 
\caption{A dumbbell flowed by our phase field model, for two different values of $\e$.}
\label{fig_test_dumbbell}
\end{figure}

A classical example of mean curvature flow where singularities appear in finite time is the evolution of a dumbbell. The classical Allen-Cahn model provides an approximate mean curvature flow beyond singularities which is compatible with Brakke's notion of weak mean curvature flow. When singularities appear, the dumbbell moved by the Allen-Cahn model separates into two sphere-type shapes that evolve until they vanish. In the proposed numerical experiment, we apply our model to a dumbbell, noting that the Willmore term is expected to prevent singular interfaces, thus topological changes.
This is indeed what we observe in the two numerical experiments shown in
Figure ~\ref{fig_test_dumbbell}, where each row corresponds to different values of
$\e$. The experiments were carried out with the following parameters: $N = 2^7$, $\e = 2/N$, $\sigma_\e = 2 \e^2$, $\delta_t = 0.01 \e^2$ and $\alpha = \beta = 0$. In both cases the two sphere-type shapes remain connected by a cylinder, the thickness of the cylinder being directly related to the value of $\e$. Notably, the cylinder appears to have minimal influence on the dynamics, probably because its lowest principal curvature is zero.

A preliminary conclusion is that our model effectively approximates the smooth mean curvature flow before the  singularities begin to appear, a point at which the Willmore term begins to have a substantial impact on the flow.

\subsubsection{Interfaces with triple junctions}
\begin{figure}[!htbp]
\centering
    \includegraphics[width=.24\textwidth]{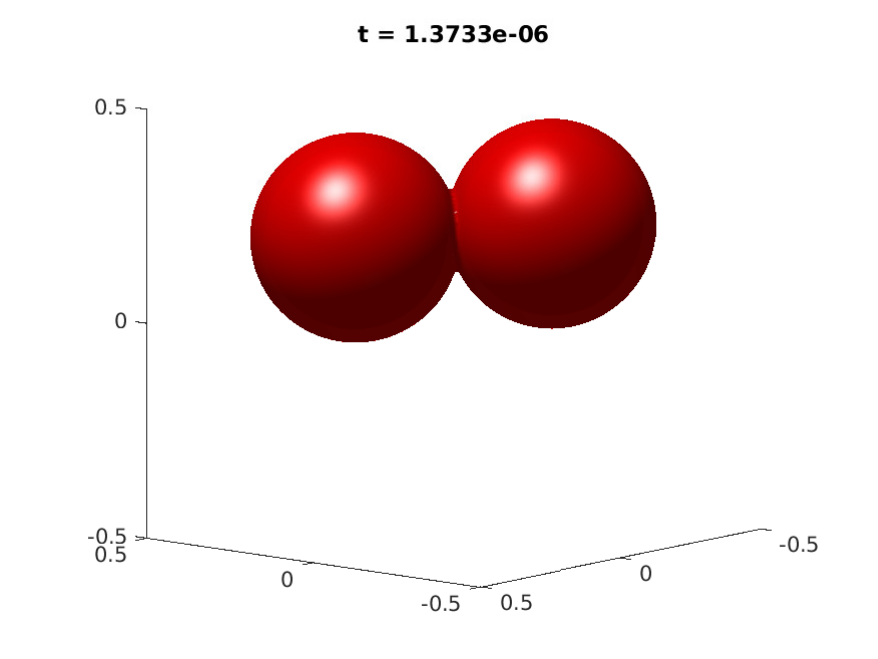}
     \includegraphics[width=.24\textwidth]{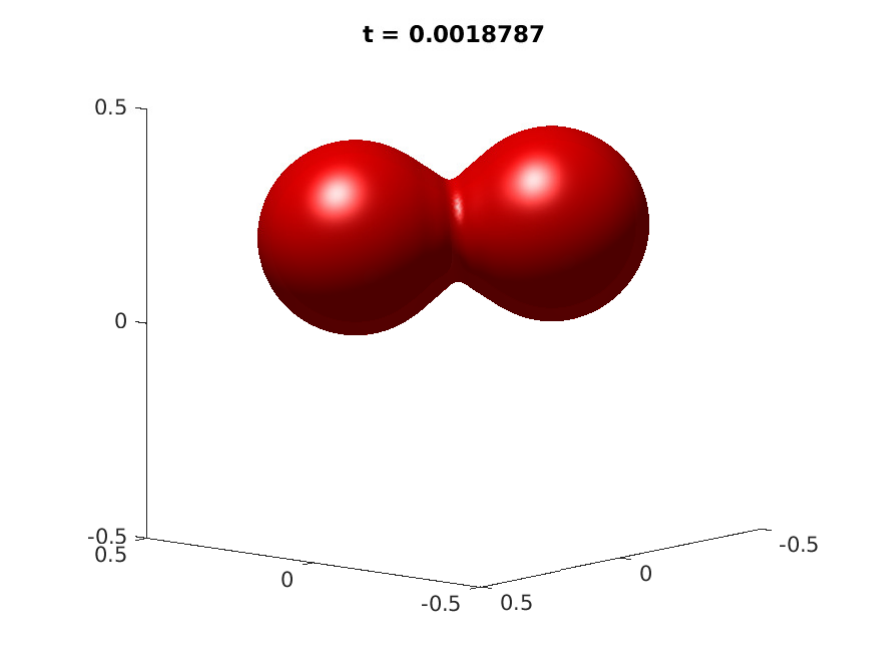}
     \includegraphics[width=.24\textwidth]{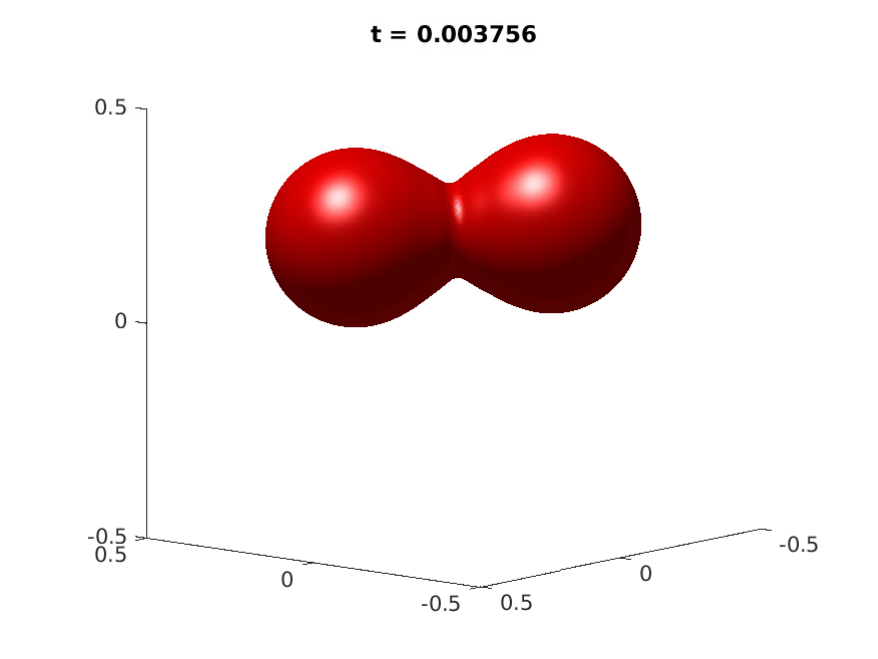}
     \includegraphics[width=.24\textwidth]{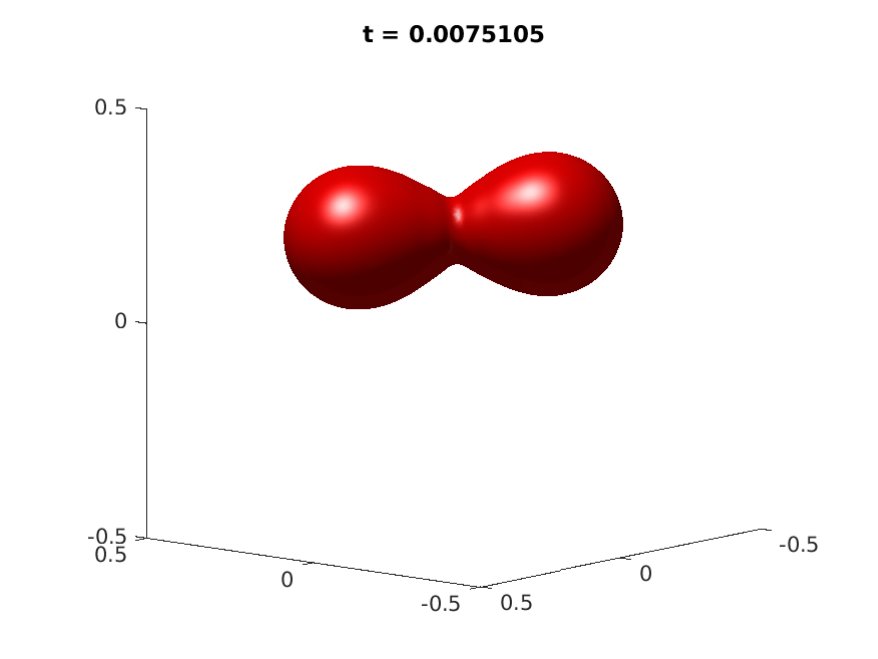} \\
    \includegraphics[width=.24\textwidth]{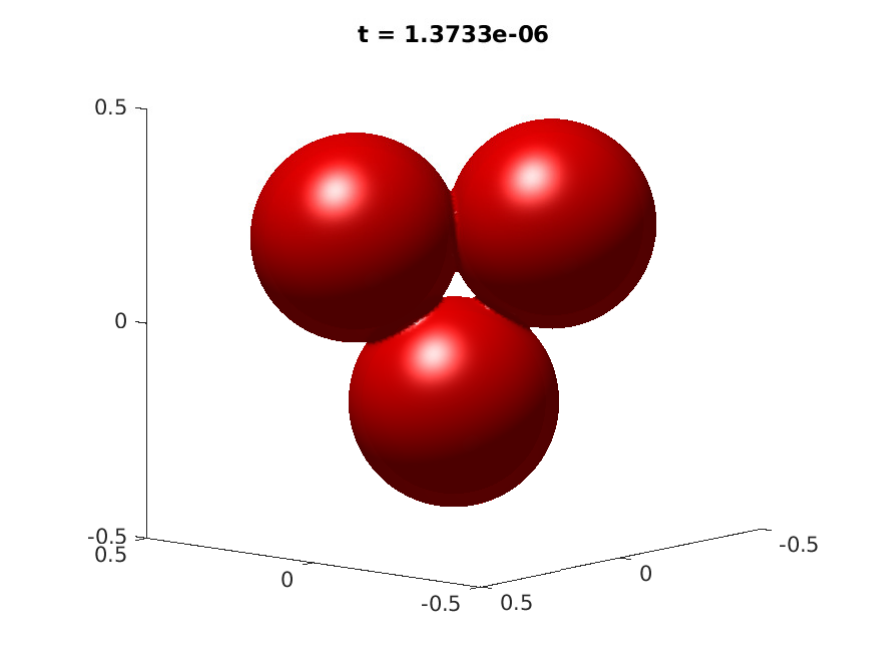}
     \includegraphics[width=.24\textwidth]{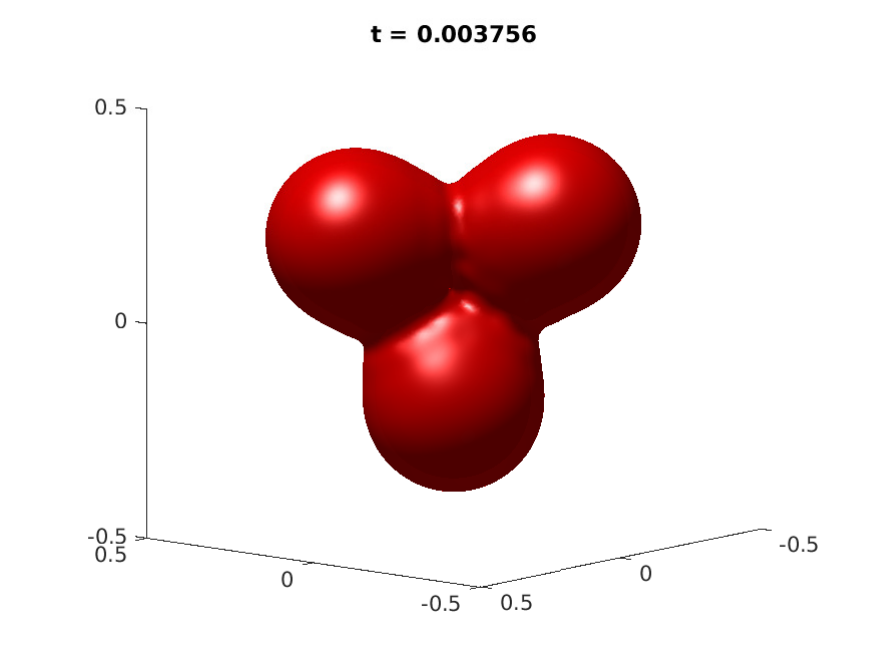}
     \includegraphics[width=.24\textwidth]{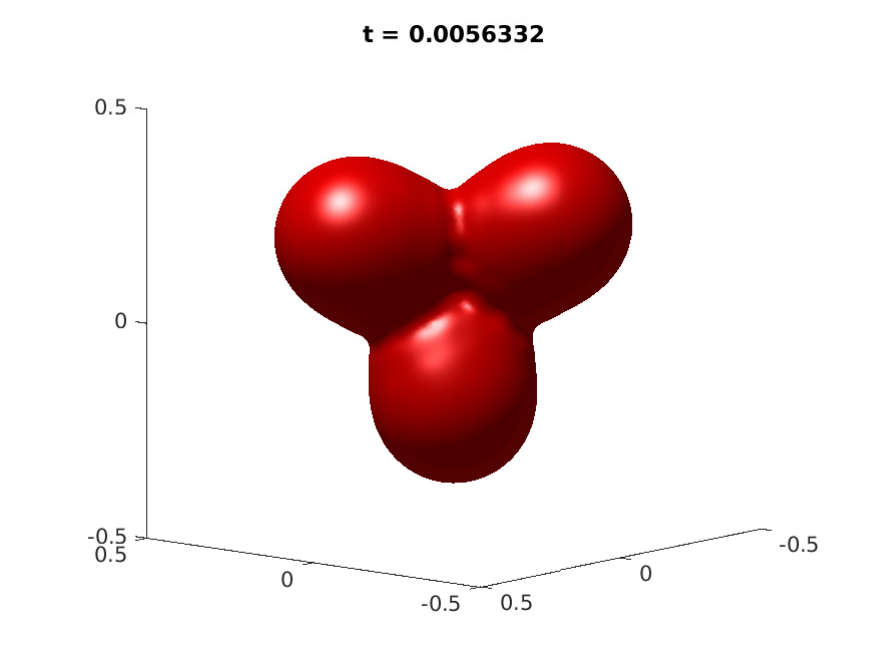}
     \includegraphics[width=.24\textwidth]{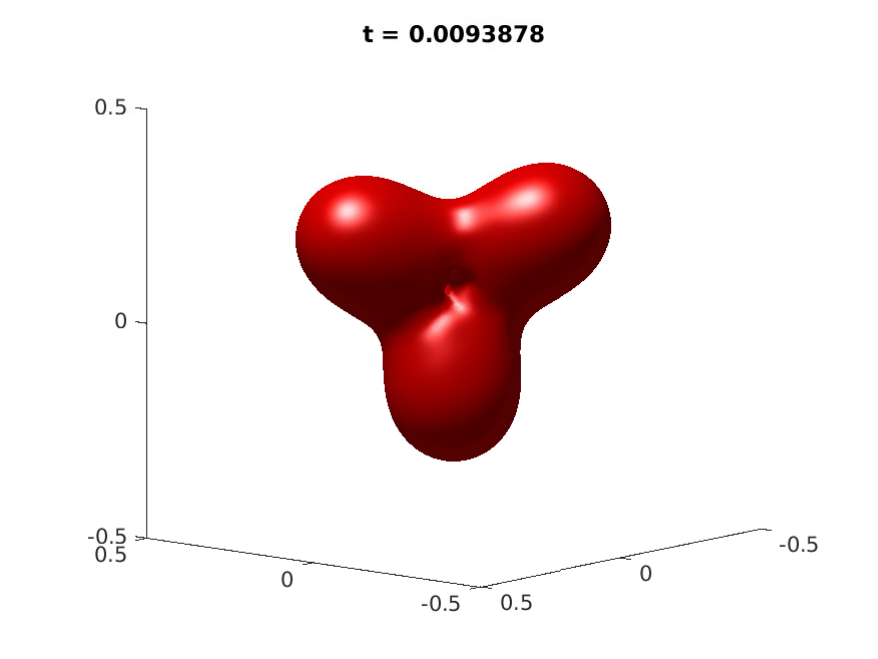} \\

\caption{Two nonsmooth shapes evolved by our model.}
\label{fig_test_non_smooth_3D}
\end{figure}

We illustrate in Figure~ \ref{fig_test_non_smooth_3D} the evolution by our phase field model of two-dimensional shapes in 3D made of two or three spheres glued together. Like in the two-dimensional case, the evolution is consistent with the result that would provide a multiphase approach, exhibiting triple lines that evolve to satisfy Herring's equilibrium conditions.

\subsubsection{Evolution of tubular sets with triple junctions}

\begin{figure}[!htbp]
\centering
	\includegraphics[width=.24\textwidth]{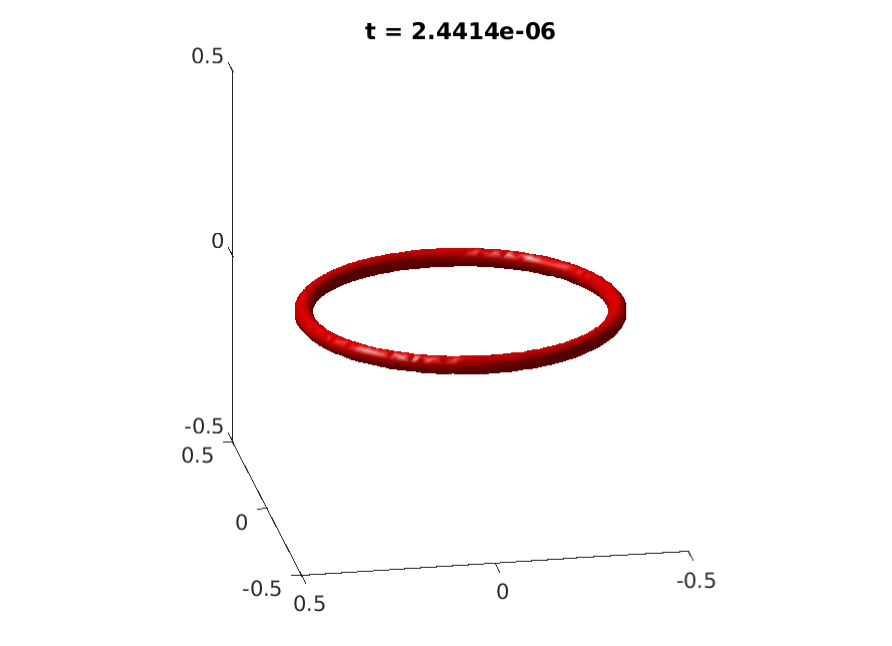}
	\includegraphics[width=.24\textwidth]{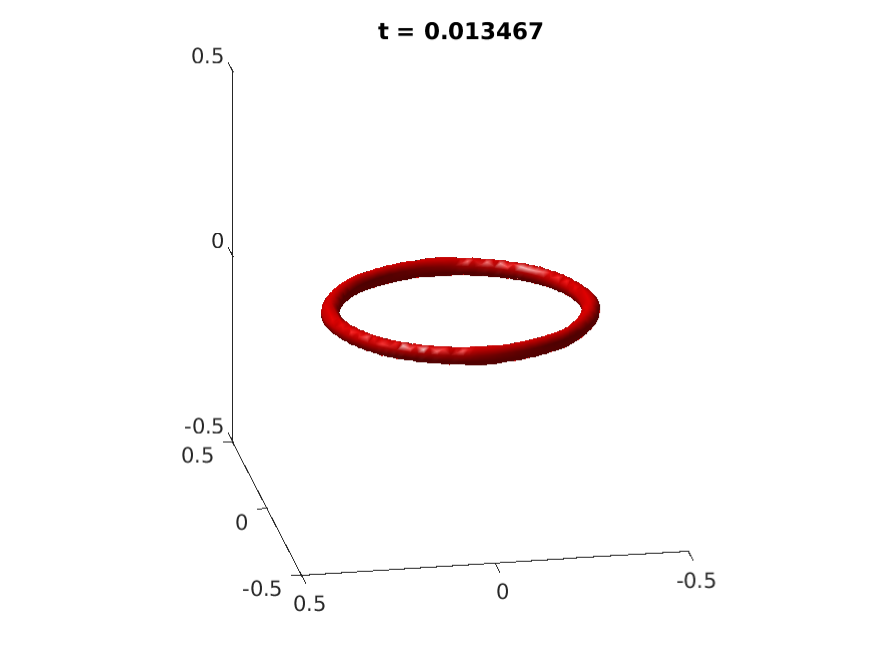}
	\includegraphics[width=.24\textwidth]{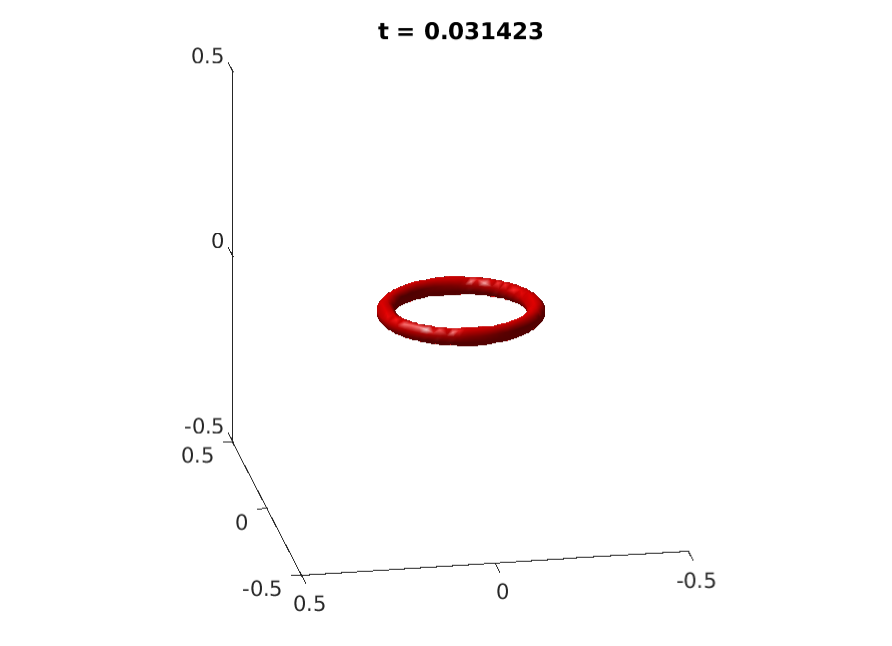}
	\includegraphics[width=.2\textwidth]{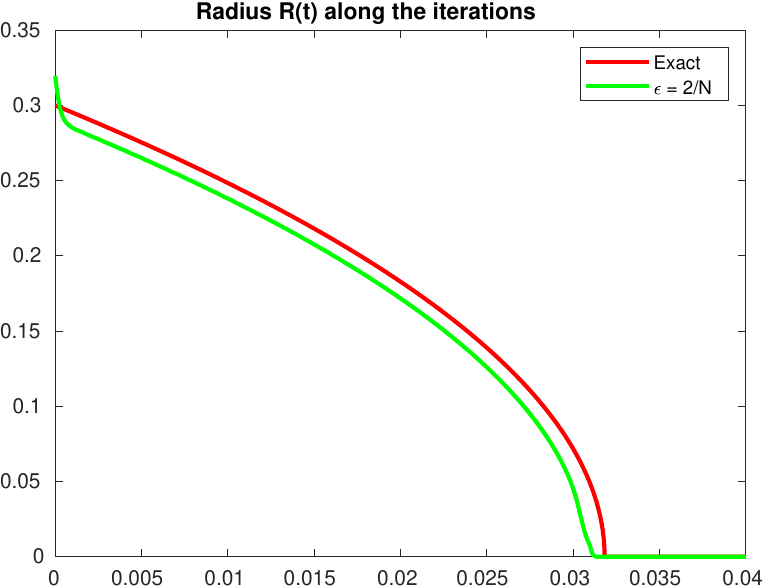}
\caption{Evolution with our model of an approximate circle in 3D, and comparison of its radius with the exact radius law of mean curvature flow.}
\label{fig_test_circle3D}
\end{figure}

As observed in a previous experiment, the Willmore term ensures the stability of the phase field profile, and its geometric properties guarantee the stability of thin structures. We illustrate now that our approach can effectively approximate the mean curvature flow of sets of codimension two represented approximately by tubular structures. The first example involves
the evolution of a circle in dimensions three. In Figure \ref{fig_test_circle3D}, we display
the level set $1/6$ across iterations of a solution $u_{\e}$ approximating the evolution of a 1D circle in 3D. The parameters are $N = 2^7$, $\e = 2/N$, $\sigma_\e = 2 \e^2$, $\delta_t = 0.01 \e^2$ and $\alpha = \beta = 0$.
As anticipated, the shapes approximate circles with radius decreasing over time. A closer look at the radius law reveals that it is consistent with a motion by mean curvature, differing only by a multiplicative factor.

Finally, we present in Figure \ref{fig_test_circle3D} the last numerical experiments performed on two filaments with periodic boundary conditions, one with triple points, the other without. Both filaments evolve toward minimal sets.

\begin{figure}[!htbp]
\centering
	\includegraphics[width=.24\textwidth]{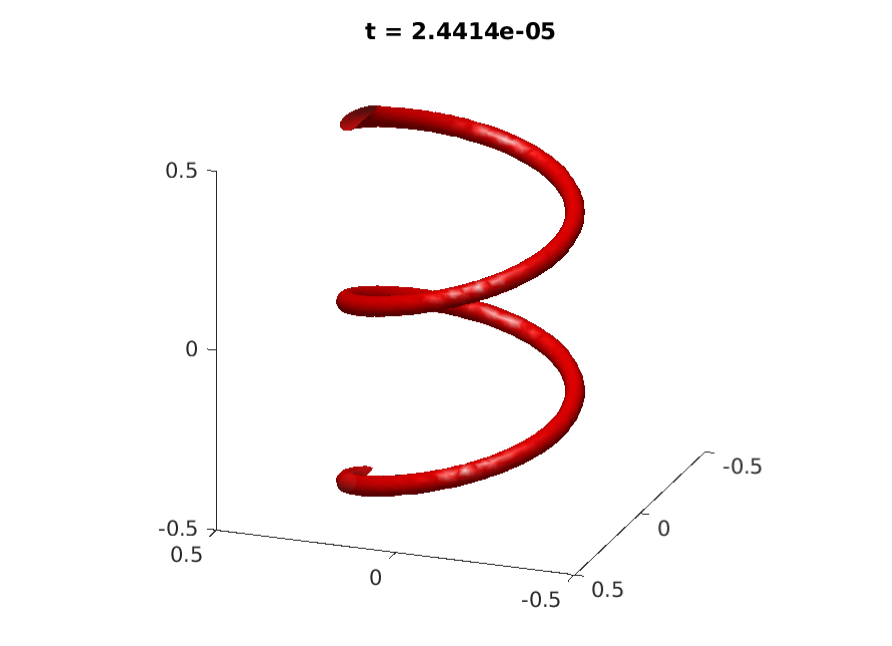}
	\includegraphics[width=.24\textwidth]{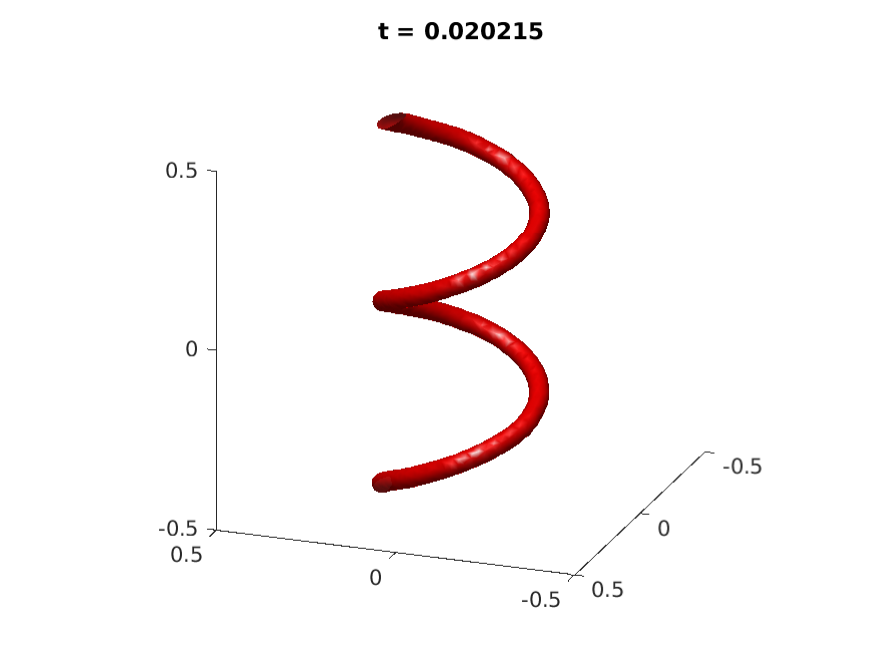}
	\includegraphics[width=.24\textwidth]{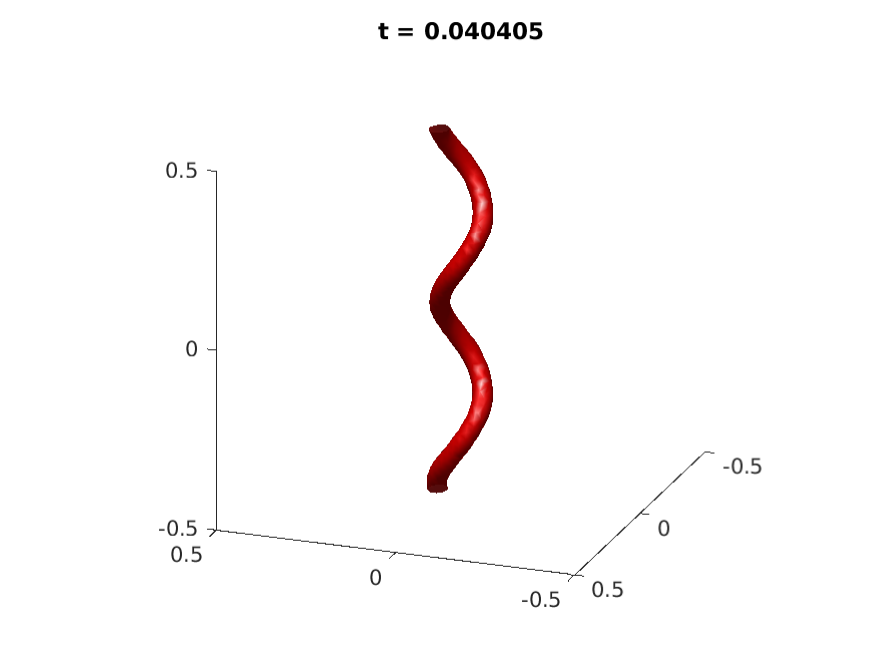}
	\includegraphics[width=.24\textwidth]{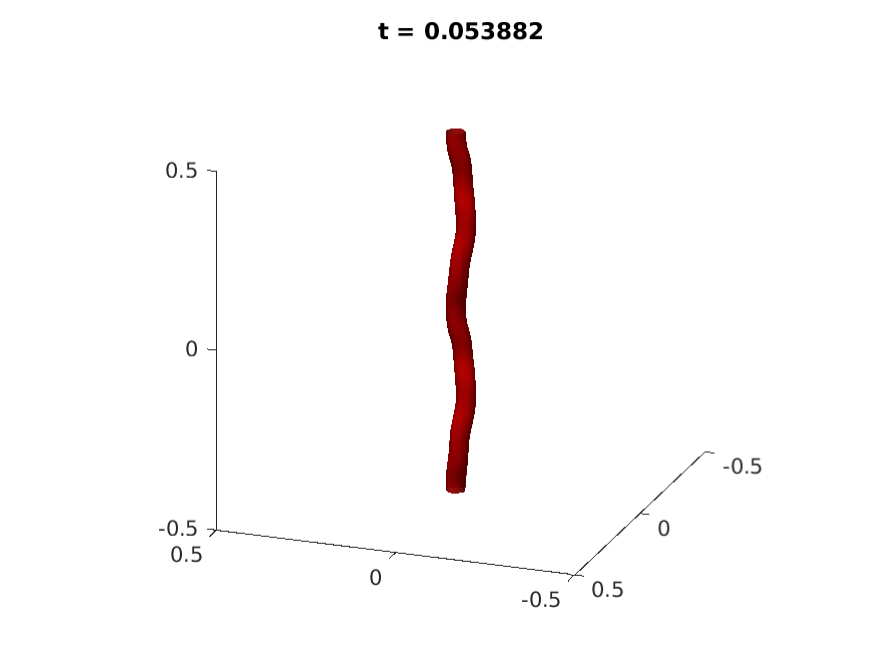} \\
	\includegraphics[width=.24\textwidth]{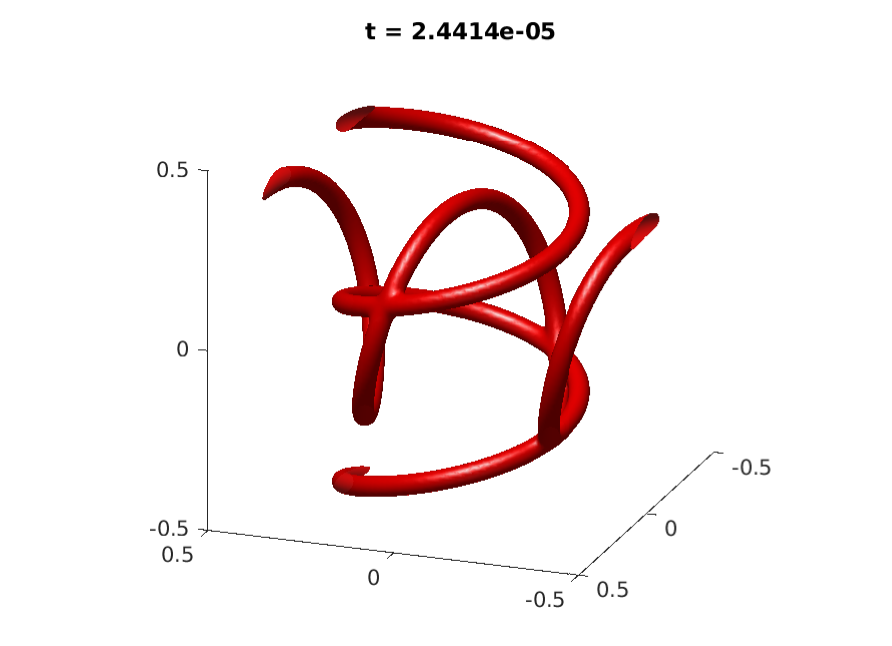}
	\includegraphics[width=.24\textwidth]{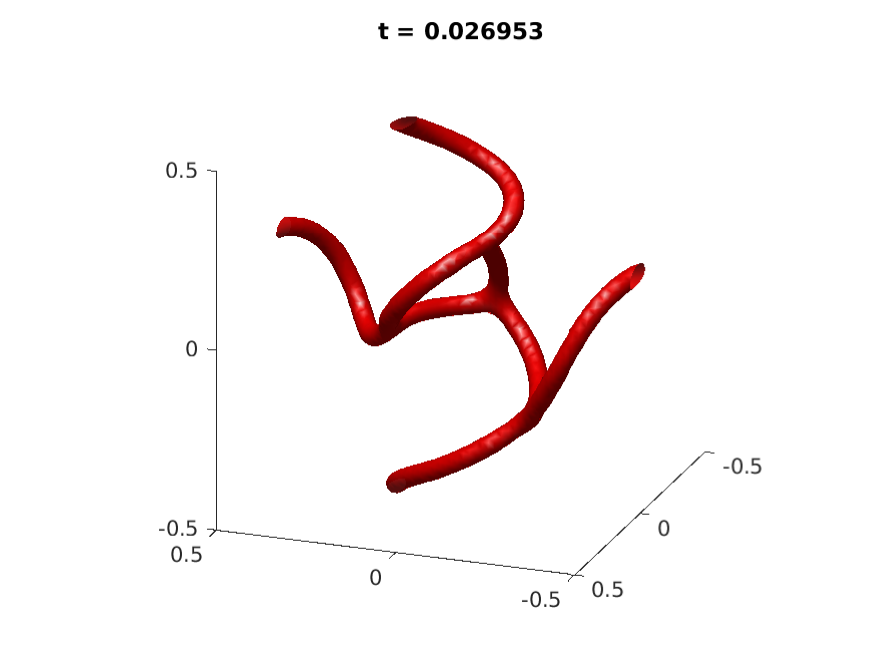}
	\includegraphics[width=.24\textwidth]{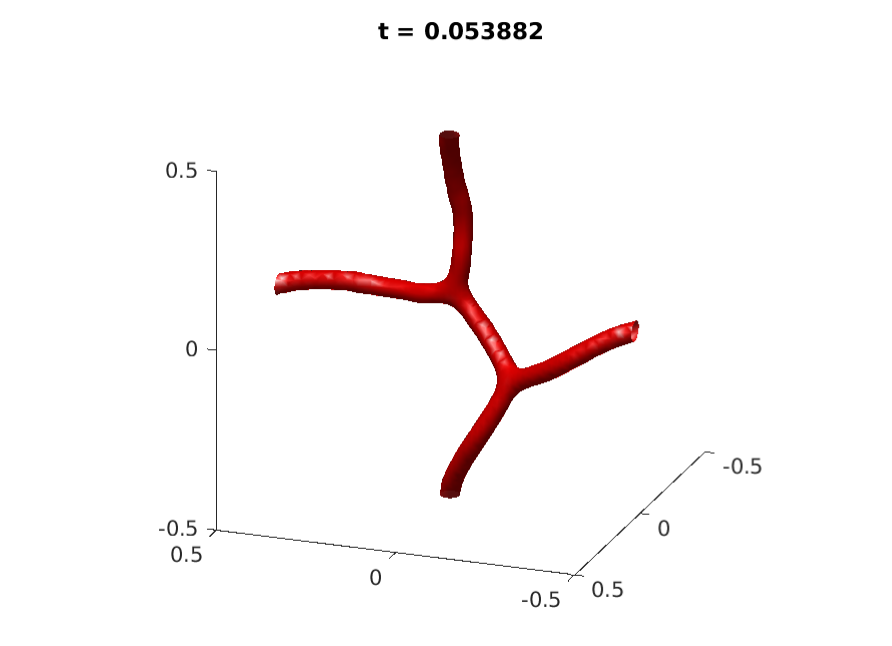}
	\includegraphics[width=.24\textwidth]{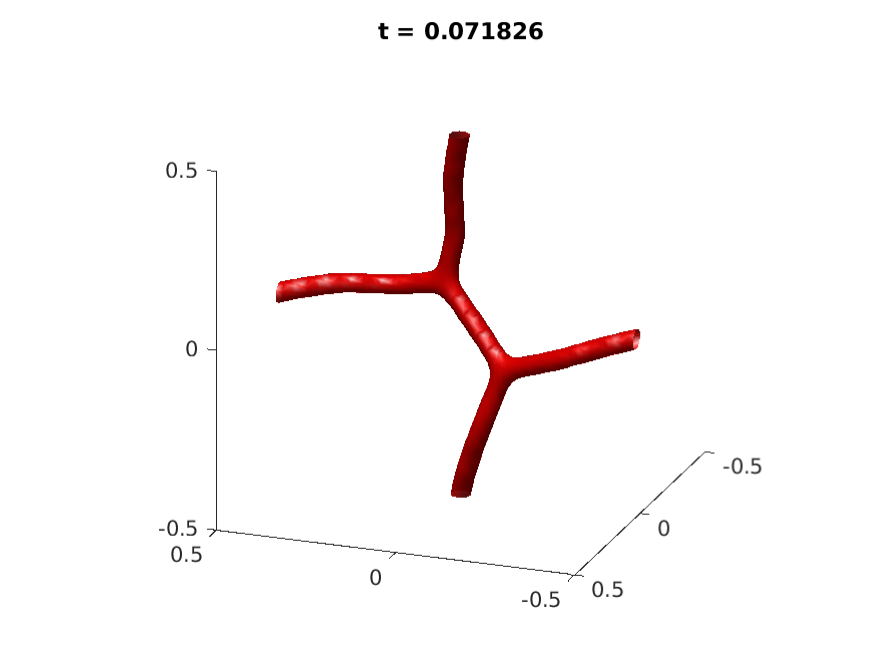} \\

\caption{Evolution of two tubular sets with periodic boundary conditions}
\label{fig_test_tubular3D}
\end{figure}

\section{Conclusion and remarks}
In this paper, we have introduced and studied a Cahn-Hilliard-type energy coupled with a phase-field Willmore-type energy.
A key aspect of our contribution is the use in the Cahn-Hilliard-type energy of a specific
potential that combines a smooth well at  $s=0$ and an obstacle at
$s=1/4$.
Our variational model admits a regular minimizer connecting these two phases, and the
 Willmore term not only promotes the stability of this regular profile along the gradient flow, it also imposes a penalty on the  $s=1/4$ phase which facilitates the approximation of the surface tension energy of non-oriented  interfaces. The asymptotic analysis of our model in dimension $1$, and in higher dimensions for radial functions, 
has enabled us to justify more rigorously these properties and the relevance of the new potential. \medskip

A surprising and theoretically unjustified result is that the
$L^2$  gradient flow of our energy is, at least numerically, a fairly good approximation of  the mean curvature flow of codimension $1$ or $2$ surfaces, which allows for the presence of triple junctions. A promising next step would be to justify this numerical observation by adapting classical asymptotic expansion techniques to our context. Additionally, it would be interesting to analyze the profiles of the solution at triple junctions and to identify more precisely the profile observed in the case of tubular sets, which likely does not correspond to the profile
$-q'$ mentioned in the introduction.

\section*{Acknowledgements}
The authors acknowledge support from the "France 2030'' funding ANR-23-PEIA-0004 ("PDE-AI") and from the French National Research Agency (ANR) under grants ANR-19-CE01-0009-01 (project MIMESIS-3D) and ANR-24-CE40-XXX (project STOIQUES). Part of this
work was also supported by the LABEX MILYON (ANR-10-LABX-0070) of Université de Lyon,
within the program "Investissements d’Avenir" (ANR-11-IDEX- 0007) operated by the French National Research Agency (ANR), and by the European Union Horizon 2020 research and innovation
programme under the Marie Sklodowska-Curie grant agreement No 777826 (NoMADS).

\end{document}